\documentclass[11pt]{amsart}

\usepackage[OT1,T1]{fontenc}
\usepackage{times}
\usepackage{calligra}
\usepackage{aurical}

\usepackage{mathpazo}
\usepackage{amsthm}
\usepackage[mathscr]{euscript}

\usepackage{stmaryrd}
\usepackage{amssymb}
\usepackage{amsmath}
\usepackage{titletoc}
\usepackage{extarrows}
\usepackage{varioref}
\usepackage{bm}
\usepackage{amsthm}
\usepackage[all,cmtip]{xy}
\usepackage{tikz-cd}
\usepackage{color} 
\usepackage{soul}
\setulcolor{red}
\sethlcolor{yellow}
\setstcolor{blue}
\usepackage{amscd}
\usepackage{xfrac}    

\usepackage{amsfonts}
\usepackage[colorlinks,linkcolor=blue,citecolor=blue, pdfstartview=FitH]{hyperref}
\usepackage{graphicx}
\usepackage{geometry}
\usepackage{cite}

\numberwithin{equation}{section}

\theoremstyle{definition}
\newtheorem{thm}{Theorem}[section]
\newtheorem{theorem}[thm]{Theorem}

\newtheorem{lemma}[thm]{Lemma}
\newtheorem{corollary}[thm]{Corollary}
\newtheorem{proposition}[thm]{Proposition}

\newtheorem{remark}[thm]{Remark}
\newtheorem*{notation}{Notation}
\newtheorem{definition}[thm]{Definition}

\newtheorem{defn-thm}[thm]{Definition-Theorem}

\begin{document}
	\title{Picard Groups of Spectral Varieties and Moduli of Higgs Sheaves}

\author{Xiaoyu Su}
\address[Xiaoyu Su]{School of Mathematical Sciences and Key Laboratory of
	Mathematics and Information Networks (Ministry of Education),
	Beijing University of Posts and Telecommunications,
	No.~10 Xitucheng Road, Beijing 100876, China}
\email{xiaoyusu@bupt.edu.cn}
\urladdr{\url{https://orcid.org/0000-0002-9567-9153}}

\author{Bin Wang}
\address[Bin Wang]{Hetao Institute of Mathematics and Interdisciplinary
	Sciences (Shenzhen), No.~3 Binglang Road,
	Shenzhen, Guangdong 518017, China}
\email{wangbin@himis-sz.cn}
\urladdr{\url{https://orcid.org/0000-0001-5668-696X}}


	\begin{abstract}
		We study moduli spaces of  Higgs sheaves valued in line bundles and the associated Hitchin maps on surfaces. We first work out Picard groups of generic (very general) spectral varieties which holds for dimension of at least 2, i.e., a Noether--Lefschetz type theorem for spectral varieties. We then apply these Picard-group calculations to obtain a necessary and sufficient condition for the non-emptiness of generic Hitchin
fibers in the range where the spectral Noether--Lefschetz theorem applies. Then we move on to detect the geometry of the moduli spaces of Higgs sheaves as the second Chern class varies.
	\end{abstract}
	\maketitle

	\section{Introduction}
	
	Tanaka--Thomas \cite{TT18, TT20} constructed Vafa--Witten invariants via moduli spaces of Gieseker-(semi)stable Higgs sheaves on a complex smooth projective surface $X$ with Higgs fields valued in the canonical bundle $\omega_X$. Gholampour--Sheshmani--Yau \cite{GSY20L} considered a similar setting on a simply-connected surface, by replacing $\omega_X$ with a line bundle ${\mathcal{L}}$ on $X$ and developed the local Donaldson--Thomas invariants which can specialize to the Vafa--Witten invariants aforementioned. Roughly speaking, moduli spaces of ${\mathcal{L}}$-valued Higgs sheaves on a surface can be treated as moduli spaces of compactly supported 2-dimensional sheaves on the total space ${\operatorname {Tot}}({\mathcal{L}}):={\operatorname{Spec}}{\operatorname {Sym}}^{\scriptscriptstyle\bullet}_{X}({\mathcal{L}}^\vee)$.  Such moduli spaces, the central object studied in this paper, admit a natural $\mathbb{G}_m$-action and the associated Hitchin maps are equivariant under this action. The Vafa--Witten invariants and local Donaldson--Thomas invariants are constructed via the localization of virtual fundamental classes which depends on the equivariant (reduced) perfect obstruction theories chosen on the moduli spaces.
	
	Though invariants are defined using ``integration'' on the $\mathbb{G}_m$-fixed point locus which is contained in the Hitchin fiber over $0$, we instead focus on generic fibers of the associated Hitchin maps. Like the curve case, the properties of generic fibers, if nonempty, can reflect many interesting geometric structures of the moduli spaces. But unlike the curve case, there maybe irreducible components which cannot be detected by generic fibers.
	
	Let $k$ be an algebraically closed field. \textit{In this introduction, we assume that $k=\mathbb{C}$ to make statements simpler.} Let $X$ be a smooth projective variety over $k$ with $\dim X\geq 2$ and ${\mathcal{L}}$ be a very ample line bundle over $X$. 
	\begin{definition}
		An ${\mathcal{L}}$-valued Higgs sheaf of rank $r$ on $X$ consists of $(\mathcal E,\theta)$ with $\mathcal E$ being a torsion free sheaf on $X$ of generic rank $r$, $\theta:\mathcal E\rightarrow \mathcal E\otimes \mathcal L$ being an $\mathcal O_{X}$-linear map. If there is no ambiguity, we will simplify by omitting "${\mathcal{L}}$-valued". When we want to emphasize that $\mathcal E$ is locally free, we call $(\mathcal E,\theta)$ a Higgs bundle. 
	\end{definition}
	Throughout the paper, we choose the stability condition given by ${\mathcal{L}}$ or any ample line bundle which is numerically equivalent to a positive scalar of ${\mathcal{L}}$ (hence they give the same stability condition). According to \cite{Yo91}, one can construct the moduli spaces of Gieseker semistable Higgs sheaves (resp. Higgs bundles), denoted by ${{\operatorname{Higgs}}}^\textrm{tf,ss}$ (reps. ${{\operatorname{Higgs}}}^\textrm{bun,ss}$  ). By definition, there is an inclusion ${{\operatorname{Higgs}}}^\textrm{bun,ss}\hookrightarrow{{\operatorname{Higgs}}}^\textrm{tf,ss}$. 
	\begin{definition}
		We call the affine space $\bm A:=\oplus_{i=1}^{r}\mathrm H^{0}(X,{\mathcal{L}}^{\otimes i})$ the Hitchin base. The following characteristic polynomial map:
		\[
		h :{{\operatorname{Higgs}}}^\textrm{tf,ss}\rightarrow\bm A,\ ({\mathcal{E}},\theta)\mapsto \text{char.poly.}(\theta)
		\]
		is called the Hitchin map.  
	\end{definition}
\noindent By Yokogawa \cite[Page 495-498, Theorem 5.10 and Corollary 5.12]{Yo93C}, the Hitchin map is projective.
 
	The first goal of the paper is to look at generic fibers of the Hitchin map. As the curve case, ``spectral varieties", which are divisors in ${\operatorname {Tot}}({\mathcal{L}})$, play an important role in the study of generic Hitchin fibers. But for later use, it is better to treat them as divisors in the projective completion, ${\operatorname{Proj}}_X {\operatorname {Sym}}_{{\mathcal{O}}_X}^{\scriptscriptstyle\bullet}({\mathcal{L}}^{\vee}\oplus {\mathcal{O}}_{X})$, of ${\operatorname {Tot}}({\mathcal{L}})$. We put
	\[
	\pi: {\operatorname{Proj}}_X {\operatorname {Sym}}_{{\mathcal{O}}_X}^{\scriptscriptstyle\bullet}({\mathcal{L}}^{\vee}\oplus {\mathcal{O}}_{X})=:Y\rightarrow X
	\]
	as the natural projection.
	\begin{definition}\label{def:spectral var}
		Spectral schemes  are divisors in the linear system $|\pi^*{\mathcal{L}}^{r}\otimes {\mathcal{O}}(r)|$ contained in ${\operatorname {Tot}}({\mathcal{L}})$,  where ${\mathcal{O}}(-1)$ is the relative tautological line bundle. They are parameterized by the Hitchin base $\bm A$. And for a closed point $s\in\bm A$, we denote the corresponding spectral scheme by $X_s$.
	\end{definition}
	Since we mainly consider generic spectral schemes which are varieties, we then call generic ones spectral varieties. As in the curve case, the spectral construction identifies generic Hitchin fibers with moduli spaces of rank-one torsion-free sheaves on the corresponding smooth spectral varieties. Their determinant--ideal decomposition relates these fibers to Picard schemes and Hilbert schemes of points (see Section~3 and Proposition \ref{prop:bundle over spectral}).	Hence, to study the generic fibers of the Hitchin map, we need to know Picard groups of generic members of smooth spectral varieties. Since they are divisors of $Y$ and ${\operatorname{Pic}}(Y)\cong {\operatorname{Pic}}(X)\oplus \mathbb{Z}$, the calculation of ${\operatorname{Pic}}(X_s)$ is revealed to be a Noether--Lefschetz type problem.

	Though ${\mathcal{L}}$ is very ample, spectral varieties are not ample divisors on $Y$. We can not apply the classical Noether--Lefschetz theory directly. But the situation is not far away from the classical ones as the line bundle $\pi^{*}{\mathcal{L}}^{r}\otimes {\mathcal{O}}_{Y/X}(r)$ is big and base-point free for any $r>0$, see Proposition \ref{prop:bigness of tauto}.  In characteristic $0$, Ravindra--Srinivas \cite{RS06,RS09} systematically deal with the Noether--Lefschetz type problem with the ample line bundle being replaced by a big and base-point free line bundle which is a nontrivial generalization of those in \cite{SGA2} and \cite[\S IV]{Har70}. We follow their strategy to work out Picard groups of (very general) generic spectral varieties.

	\begin{theorem}(see Theorem \ref{thm:eff lef for formal line bundle}, Theorem \ref{thm:NL for surfaces})
		 Let $X_s$ be a smooth spectral variety, and $\pi_s:X_s\rightarrow X$ is the restriction of $\pi:Y\rightarrow X$. Then $\pi^{*}_s:{\operatorname{Pic}}(X)\rightarrow{\operatorname{Pic}}(X_{s})$ is an isomorphism if
		\begin{enumerate}
			\item when $\dim X\ge 3$.
			\item $X_s$ is very general, i.e., $s$ lies in the complement of countably many proper closed subvarieties of $\bm A$ and $r\ge 4$, when $\dim X=2$.
		\end{enumerate}		
	\end{theorem}
	Since spectral varieties are quite special, we want to point it out that under certain assumptions, the above isomorphism holds over positive characteristics for $\dim X\ge 3$.  We will address this to Subsection \ref{subsec:higher dim}.

		We now assume that $\dim X=2$ and $k=\mathbb{C}$. For a fixed choice of Chern classes $(r,c_1,c_2)\in \oplus_{i=0}^{2}\mathrm{H}^{2i}(X,\mathbb Q)$, we denote the moduli space of Gieseker semistable Higgs sheaves by ${\operatorname{Higgs}}_{r,c_1,c_2}^{\text{
				tf,ss}}$ which admits a $\mathbb{G}_m$-action. 

		\begin{theorem}(see Theorem \ref{thm:criterion of existence}, Corollary \ref{cor:genericsmooth})\label{thm:intro vw1}
			Assume that
$\pi_s^*:{\operatorname{Pic}}(X)\to{\operatorname{Pic}}(X_s)$
is an isomorphism for a very general $s\in\bm A$. Generic fibers of the Hitchin map
			\[
			h:{{\operatorname{Higgs}}}_{r,c_{1},c_{2}}^\textrm{tf,ss}\rightarrow \bm{A}
			\]
			are nonempty,  if and only if the following equations:
			\begin{align*}
				&r\delta=c_1+\frac{r(r-1)}{2}c_1({\mathcal{L}})\\
				&(r-1)c_1^{2}-2rc_2=\frac{r^2(r^2-1)}{12}c_1({{\mathcal{L}}})^2-2r\pi_{s*}[\mathfrak D].
			\end{align*}
			have solutions for $\delta\in \text{Im}(c_1:{\operatorname{Pic}}(X)\rightarrow\mathrm{H}^{2}(X,\mathbb Q))$  and effective $[\mathfrak D]\in \mathrm {H}^{4}(X_s,\mathbb Z)$. Here $X_s$ is a general spectral surface and effective means $\int_{X_s}[\mathfrak D]\geq 0$.

 If the generic fibers are non-empty, the Hitchin map is surjective by its projectiveness. Moreover, in this case, the Hitchin map $h$ is generically smooth. i.e. there exist a non-empty Zariski open subset $W\subset \bm A$, such that $h|_{h^{-1}W}:h^{-1}W\to W$ is smooth.
		\end{theorem}

		Following \cite{TT20}, for a line bundle ${\mathcal{N}}$ on $X$, we denote  by ${\operatorname{Higgs}}_{r,{\mathcal{N}},c_2}^{\textrm{
				tf,ss}}$  the corresponding moduli space  of Gieseker semistable Higgs sheaves with fixed determinant ${\mathcal{N}}$ and tracefree Higgs fields. In \cite{TT20}, Tanaka and Thomas take ${\mathcal{L}}=\omega_X$. When semistability coincides with stability, they localize the virtual fundamental class to the fixed point components and define the Vafa--Witten invariant as 
		$$\mathsf{VW}_{r,c_1,c_2}=\int_{[({\operatorname{Higgs}}_{r,\mathcal N,c_{2}}^\textrm{tf,ss})^{\mathbb G_m}]^\textrm {vir}}\frac{1}{e(N^\textrm{vir})},$$
		where $c_1=c_1({\mathcal{N}})$, $N^\textrm{vir}$ is the virtual normal bundle. It is conjectured that, for some $s\in\mathbb{C}$, the following generating series (see \cite{JK22} for a modified and precised version):
		\[
		Z_r(X,{\mathcal{N}})=q^{-s}\sum_{n\in\mathbb Z}\mathsf{VW}_{r,{\mathcal{N}},n}\cdot q^{n}
		\]
		has certain modularity. In general, Tanaka--Thomas \cite{TT18} define Vafa--Witten invariants via the moduli space of Joyce--Song pairs which are equal to the aforementioned ones when semistability coincide with the stability. 
		Irreducible components in the fixed point locus $({\operatorname{Higgs}}_{r,{\mathcal{N}},c_2}^{\text{tf,ss}})^{\mathbb{G}_m}$ with zero Higgs fields are called the instanton branches and the others are called monopole branches. 
		
		\begin{definition}
			We define 
			$$c_2^{\text{g.bun}}:=\frac{(r-1)}{2r}c_1^2-\frac{r(r^2-1)}{24}c_1(\mathcal L)^2.$$ 
			which is a rational number (depends on $r,c_1,\mathcal L$).
		\end{definition}
		Then we can state our main theorem on the moduli space ${\operatorname{Higgs}}^{\text{tf,ss}}_{r,c_1,c_2}$.
		\begin{theorem}[see Theorem \ref{thm:graded piece are simple}, Theorem \ref{thm:small c2}]\label{thm:main intro} We assume that the Equation:
			\[
			r\delta=c_1+\frac{r(r-1)}{2}c_1(\mathcal{L})
			\]
			has a solution $\delta \in \text{Im}(c_1:{\operatorname{Pic}}(X)\rightarrow\mathrm{H}^{2}(X,\mathbb Q)) $.   Consequently, $c_2^{\text{g.bun}}$  is an integer. Then:

			\begin{enumerate}   
				\item If $c_2<c_2^{\text{g.bun}}$, then ${\operatorname{Higgs}}^{\text{tf},ss}_{r,c_1,c_2}$ is empty.
				\item If $c_2=c_2^{\text{g.bun}}$, then for $({\mathcal{E}},\theta)\in {\operatorname{Higgs}}^{\text{tf},ss}_{r,c_1,c_2}$,  the only possible Harder--Narasimhan type of $\mathcal E$ is that each graded piece is of generic rank $1$. 
			 
				\item Assume in addition that
$\pi_s^*:{\operatorname{Pic}}(X)\xrightarrow{\cong}{\operatorname{Pic}}(X_s)$
for a very general $s\in\bm A$.
Then, for $c_2\gg 0$, the moduli space ${\operatorname{Higgs}}_{r,c_1,c_2}^{\text{tf,ss}}$ is reducible.
			\end{enumerate}
		\end{theorem}
		
		The above theorem shows that the moduli spaces are generally reducible for $c_2\gg 0$. We can show the connectedness of the moduli spaces of Higgs sheaves of rank 2 on degree $d\geq 5$ hyper-surfaces in $\mathbb P^3$ for $c_2\gg 0$.
		\begin{remark}
			The above theorem shows that Vafa--Witten invariants for $c_2<c_2^{\textrm{g.bun}}$ vanish supposing that the Gieseker semistability coincides with stability. This bounds the possible pole order of $Z_r(X,{\mathcal{N}})$.
		\end{remark}

		Let us now indicate how this relates to previous work. Over $\mathbb{C}$, the notion of Infinitesimal Noether--Lefschetz Condition was first proposed by Joshi \cite{Jos95} to solve Noether--Lefschetz problem on a smooth projective complex variety with an ample linear system. As mentioned before, our calculation of Picard groups follows Ravindra and Srinivas' strategy for big and base-point free linear systems. 
		Since the linear system of spectral varieties is quite special, we adapt their method to prove Noether-Lefschetz type theorem in positive characteristics under certain assumptions.

		We should also mention that in positive characteristics, Ji \cite{Ji24} proved Noether--Lefschetz type theorems, up to $p$-torsions, for a linear system with sufficient ampleness for normal threefolds. In Theorem \ref{thm:NL for surfaces}, if we take $X$ to be $\mathbb{P}^{2}$, then spectral surfaces will be isomorphic to hypersurface in $\mathbb{P}^{3}$, which then is a special case of Deligne\cite{Del73} where the Picard groups of almost all complete intersections are shown to be $\mathbb{Z}$ for any characteristics.
		
		Gholampour--Sheshmani--Yau\cite{GSY20L,GSY20} work on the virtual fundamental classes of general nested Hilbert scheme of points and computed their contribution to the Vafa--Witten invariants in rank two cases. Gholampour--Thomas \cite{GT20, GT20II}, Thomas \cite{Tho20}, Laarakker\cite{Laa20, Laa21} have extensively studied virtual fundamental classes of nested Hilbert schemes and their contributions to Vafa--Witten invariants as the vertical components. In particular, Thomas \cite{Tho20} showed that if the rank is prime then only the nested Hilbert schemes contribute to the Vafa--Witten invariants. Jiang--Kool \cite{JK22} proved the modularity for the generating series of Vafa--Witten invariants on K3 surfaces for prime ranks. Our results show that starting from generic fibers, we can detect many geometric properties of the moduli spaces and give a bound on the possible leading exponent of the generating series. However such a method also has its drawback. As implied in the Theorem \ref{thm:main intro}, there are irreducible components which cannot be detected via generic fibers.

		We now close this section by briefly describing how the paper is organized. In the second section,  we prove our results on Picard groups. In the third section, we apply the calculation of Picard groups to study generic fibers. In the fourth section, we prove Theorem \ref{thm:main intro}. In the fifth section, we show the connectedness of the moduli spaces of Higgs sheaves of rank 2 on degree $d\geq 5$ hyper-surfaces in $\mathbb P^3$.\vspace{5pt}

        \noindent\textbf{Acknowledgements.}
        We are grateful to Qizheng Yin for many helpful discussions and suggestions. We also want to thank Michael McBreen, Conan Nai Chung Leung, Xucheng Zhang, Quan Shi and Muyao Zou for useful discussions. We would also like to thank the anonymous referee for suggesting a substantially shorter and more transparent proof of Lemma~\ref{lem:olympic} based on a simple splitting argument.
        
        X. Su is partially supported by NSFC (No.12301056) and the Fundamental Research Funds for the Central Universities. B. Wang conducted part of this work at the Steklov Institute of Mathematics (Moscow), supported by the Ministry of Science and Higher Education of the Russian Federation (agreement no. 075-15-2019-1614 ). He was supported by funding from the Department of Mathematics at the Chinese University of Hong Kong, the General Research Fund (Project code: 14307022), and the Early Career Scheme (Project code: 24307121) of Michael McBreen.
        
		\section{Picard Groups of Spectral Varieties}
		In this section, we calculate the Picard group of a generic (very general) spectral varieties. These types of problems are often referred to as Noether--Lefschetz type problems. The classical Noether–Lefschetz theorem asserts that for a very general divisor $S$ in the linear system $|\mathcal O_{\mathbb P^3}(d)|$ of degree $\geq 4$ in $\mathbb P^3$ over the complex numbers, the restriction map from the Picard group on $\mathbb P^3$ to $S$ is an isomorphism. In our case, we consider the Picard group of spectral varieties vary in a big and base-point free linear system.

		Let $k$ be an algebraically closed field of characteristic 0. Let $X$ be a smooth projective variety over $k$ with $\dim X\ge 2$ and ${\mathcal{L}}$ be a very ample line bundle over $X$.

		Recall that in the introduction, we denote the total space of the line bundle ${\mathcal{L}}$ by ${\operatorname {Tot}}({\mathcal{L}})$, its projective completion by
		$\pi:Y:=\mathbf P({\mathcal{L}}^\vee\oplus{\mathcal{O}}_X)\to X$ and $Y={\operatorname {Tot}}({\mathcal{L}})\sqcup D_{\infty}$. 
		As we specify the line bundle ${\mathcal{L}}$, we may denote ${\mathcal{O}}_{\mathbf P({\mathcal{L}}^{\vee}\oplus {\mathcal{O}})/X}(1)$ by ${\mathcal{O}}_{Y/X}(1)$ for simplicity.

		\begin{remark}
			It is straightforward to see that $\pi^*{\mathcal{L}}\otimes{\mathcal{O}}_{Y/X}(1)|_{D_\infty}$ is trivial, hence it is not ample. 
		\end{remark}

		\subsection{Vanishing Properties}
		In this subsection, we will prove several cohomological vanishing results for later use.
	
		\begin{proposition}\label{prop:bigness of tauto}
			The line bundle $\pi^*{\mathcal{L}}\otimes{\mathcal{O}}_{Y/X}(1)$ is big and base-point free on $Y$. 
		\end{proposition}

		\begin{proof}
		Since ${\mathcal{L}}$ is globally generated, and ${\mathcal{L}}\otimes \pi_*{\mathcal{O}}_{Y/X}(1)\cong {\mathcal{O}}_X\oplus{\mathcal{L}}$, we have the surjective evaluation map: 
			\[e:H^0(X,{\mathcal{O}}_X\oplus{\mathcal{L}})\otimes{\mathcal{O}}_X\twoheadrightarrow {\mathcal{O}}_X\oplus{\mathcal{L}}.\] Applying $\pi^*$, we have the factorization of the evaluation map
			\[	\begin{tikzcd} H^{0}(Y,\pi^{*}{\mathcal{L}}\otimes{\mathcal{O}}_{Y/X}(1))\otimes_k{\mathcal{O}}_{Y}\arrow[r,"{\text{ev}_{\pi^*{\mathcal{L}}\otimes{\mathcal{O}}_{Y/X}(1)}}"]&\pi^{*}{\mathcal{L}}\otimes{\mathcal{O}}_{Y/X}(1)\\
			\pi^*(H^{0}(X,{\mathcal{L}}\otimes \pi_*{\mathcal{O}}_{Y/X}(1))\otimes_k {\mathcal{O}}_{X})\arrow[r,twoheadrightarrow,"{\pi^*e}"]\arrow[u,"="]&
			\pi^*{\mathcal{L}}	\otimes_{{\mathcal{O}}_Y}		\pi^*\pi_*{\mathcal{O}}_{Y/X}(1)\arrow[u,twoheadrightarrow]
			\end{tikzcd},\]
			which shows that $\pi^{*}{\mathcal{L}}\otimes{\mathcal{O}}_{Y/X}(1)$ is base-point free. Hence $\pi^{*}{\mathcal{L}}\otimes{\mathcal{O}}_{Y/X}(1)$ is nef. To show that $\pi^{*}{\mathcal{L}}\otimes{\mathcal{O}}_{Y/X}(1)$ is big, we just have to check that the intersection number $(\pi^{*}{\mathcal{L}}\otimes{\mathcal{O}}_{Y/X}(1))^{\dim X+1}>0$.  But $\sigma (X)$ is a zero divisor of $\pi^{*}{\mathcal{L}}\otimes{\mathcal{O}}_{Y/X}(1)$, and $\pi^{*}{\mathcal{L}}\otimes{\mathcal{O}}_{Y/X}(1)|_{\sigma(X)}\cong {\mathcal{L}}$. Then $(\pi^{*}{\mathcal{L}}\otimes{\mathcal{O}}_{Y/X}(1))^{\dim X+1}=({\mathcal{L}})^{\dim X}>0$ which follows from the bigness of ${\mathcal{L}}$ on $X$.   	
		\end{proof}
	
		\begin{lemma}\label{lem:can bund formula}
			Let $\omega_{Y}$ be the canonical line bundle of $Y$, then:
			\[
			\omega_{Y}\cong \pi^{*}(\omega_{X}\otimes{\mathcal{L}}^{\vee})\otimes{\mathcal{O}}_{Y/X}(-2)
			\]
		\end{lemma}
		\begin{proof}
			By the relative Euler exact sequence:
			\[
			0\to \Omega^1_{Y/X}\to \pi^*({\mathcal{L}}^\vee\oplus{\mathcal{O}}_X)\otimes{\mathcal{O}}_{Y/X}(-1)\to {\mathcal{O}}_{Y}\to 0.
			\]
			we can calculate that $\omega_{Y}\cong \det(\Omega_{Y/X}^1)\otimes\pi^{*}\omega_{X}\cong \pi^{*}(\omega_{X}\otimes{\mathcal{L}}^{\vee})\otimes{\mathcal{O}}_{Y/X}(-2)$.
		\end{proof}
		
		\begin{lemma}\label{lem: a vanishing theorem}
			$H^{i}(Y,\pi^{*}{\mathcal{L}}^{-n}\otimes{\mathcal{O}}_{Y/X}(-n))=0$ for $i<3$ and any $n\geq 1$.
		\end{lemma}
		
		\begin{proof}

			Recall $\pi^{*}{\mathcal{L}}\otimes{\mathcal{O}}_{Y/X}(1)\cong{\mathcal{O}}_{\mathbf P({\mathcal{L}}\otimes{\mathcal{O}})/X}(1)$ and $Y=\mathbf P({\mathcal{L}}^\vee\oplus{\mathcal{O}})\cong \mathbf P({\mathcal{L}}\oplus {\mathcal{O}})$ is a projective bundle on $X$, thus by the direct image formula of projective bundles, $R\pi_{*}(\pi^{*}{\mathcal{L}}^{-n}\otimes{\mathcal{O}}_{Y/X}(-n))\cong \oplus_{\ell=1}^{n-1}{\mathcal{L}}^{-\ell}[-1]$ for $n>1$ and $0$ for $n=1$.
			Thus we have $$H^{i}(Y,\pi^{*}{\mathcal{L}}^{-n}\otimes{\mathcal{O}}_{Y/X}(-n))\cong H^{i-1}(X,\oplus_{\ell=1}^{n-1}{\mathcal{L}}^{-\ell}).$$ By the Kodaira vanishing theorem, we have $H^{i}(Y,\pi^{*}{\mathcal{L}}^{-n}\otimes{\mathcal{O}}_{Y/X}(-n))=0$ for $i<3$. 
		\end{proof}
		By Proposition \ref{prop:bigness of tauto}, this also follows straightforward from  Kawamata--Viehweg Vanishing theorem since we are in characteristic 0.

For notational convenience, put
\[
\mathcal A:=\pi^{*}\mathcal L\otimes\mathcal O_{Y/X}(1),
\qquad
\mathcal W:=\mathcal A^{\otimes r}
 =\pi^{*}\mathcal L^{r}\otimes\mathcal O_{Y/X}(r).
\]
Let
\[
\phi:=\phi_{|\mathcal A|}\colon
Y\rightarrow \mathbb P\bigl(H^{0}(Y,\mathcal A)^{\vee}\bigr)\ \text{and}\  
g:=\phi_{|\mathcal W|}\colon
Y\rightarrow \mathbb P\bigl(H^{0}(Y,\mathcal W)^{\vee}\bigr).
\]
Set
\[
U_r:=
\operatorname{Im}\!\left(
\operatorname{Sym}^{r}H^{0}(Y,\mathcal A)
\rightarrow H^{0}(Y,\mathcal W)
\right).
\]
The inclusion $U_r\subset H^{0}(Y,\mathcal W)$ induces a linear projection
\[
\rho\colon
\mathbb P\bigl(H^{0}(Y,\mathcal W)^{\vee}\bigr)
\dashrightarrow
\mathbb P(U_r^{\vee}).
\]
Since the sub-linear system $U_r$ is base-point free, the center of this projection is disjoint from $g(Y)$. Moreover, after identifying
$\mathbb P(U_r^{\vee})$ with the linear span of the image of the $r$-fold Veronese morphism, we have
\[
\rho\circ g=\nu_r\circ\phi.
\]
Here, $\nu_r$ is the $r$-fold Veronese embedding. Thus, more precisely, the morphism $\nu_r\circ\phi$ is obtained from $g$ by a linear projection.

Let $Z:=g(Y)$, let $Z_1:=\phi(Y)$, and put
$\widetilde Z:=\nu_r(Z_1)$. Denote by
$\widetilde{\mathcal H}$ the line bundle on $\widetilde Z$ corresponding, under the isomorphism
$Z_1\xrightarrow{\sim}\widetilde Z$, to
$\mathcal O_{Z_1}(1)$. Then
\[
\widetilde{\mathcal H}^{\otimes r}
 \cong \mathcal O_{\widetilde Z}(1).
\]
The restriction $\rho_Z:=\rho|_Z\colon Z\to\widetilde Z$ is a morphism. Define
\[
\mathcal H:=\rho_Z^{*}\widetilde{\mathcal H}.
\]
It follows that
\[
g^{*}\mathcal H\cong\mathcal A
\qquad\text{and}\qquad
\mathcal H^{\otimes r}
 \cong
 \mathcal O_{\mathbb P(H^{0}(Y,\mathcal W)^{\vee})}(1)|_Z.
\]
In particular, $\mathcal H$ is ample and globally generated.

	The following proposition and its corollary will only be used when $X$ is a surface.
	\begin{proposition}\label{prop:global gen}
		Let $X$ be a smooth projective surface. For $r>3$, we have $g_{*}\omega_{Y}\otimes{\mathcal{O}}_{\mathbf P(W)}(1)=g_*(\omega_Y\otimes {\mathcal{W}})$ is Castelnuovo--Mumford $0$-regular with respect to the ample and globally generated line bundle ${\mathcal{H}}$ on $Z$.
	\end{proposition}
	Recall that a coherent sheaf ${\mathcal{F}}$ on $Z$ is $m$-regular with respect to ${\mathcal{H}}$, if $H^{q}(Z,{\mathcal{F}}\otimes{\mathcal{H}}^{\otimes(m-q)}))=0$ for $q\geq 1$.
	\begin{proof}
		We first prove that $Rg_{*}\omega_{Y}\cong g_{*}\omega_{Y}$. Recall that for $f : B\to C$ be a proper morphism between varieties and ${\mathcal{F}}$ a coherent sheaf on $B$.
		The following are equivalent, see \cite[Proposition 2.69]{KM98}:
		\begin{itemize}
			\item $H^q(B, {\mathcal{F}}\otimes f^*{\mathcal{O}}_C(H)) =0$ for $H$ sufficiently ample,
			\item $R^qf_*{\mathcal{F}}=0$.
		\end{itemize}
		We take $B=Y, C=Z, {\mathcal{F}}=\omega_{Y}$, and $H={\mathcal{H}}^{\otimes \ell} $ for $\ell\gg 0$. Then by Lemma \ref{lem: a vanishing theorem}, we have $R^{q}g_{*}\omega_{Y}=0$ for $q=1,2,3$. Since $\dim Y=3$, we have $Rg_*\omega_Y=g_*\omega_Y$.

		Since $r>3$, and $g^{*}{\mathcal{H}}\cong \pi^{*}{\mathcal{L}}\otimes{\mathcal{O}}_{Y/X}(1)$, then again by Lemma \ref{lem: a vanishing theorem}, for $q=1,2,3$ we have 
		\begin{small}
			\[ H^q(Z,g_*(\omega_Y\otimes {\mathcal{W}})\otimes {\mathcal{H}}^{\otimes(-q)})
			=H^{3-q}(Y,\pi^*{\mathcal{L}}^{\otimes q-r}\otimes{\mathcal{O}}_{Y/X}(q-r))=0.
			\]
		\end{small}
		And the zero-regularity of $g_*(\omega_Y\otimes_Y{\mathcal{W}})$ follows.
	\end{proof}

	In particular, we have:
	\begin{corollary}\label{cor:surjection we need}
		For any $\ell\ge 0$, we have the following surjection:
		\[H^{0}(Z,g_{*}(\omega_{Y}\otimes{\mathcal{W}}))\otimes H^{0}(Z,{\mathcal{O}}_Z(\ell))\twoheadrightarrow H^{0}(Z,g_{*}(\omega_{Y}\otimes {\mathcal{W}})\otimes {\mathcal{O}}_Z(\ell)).\]
	\end{corollary}
	\begin{proof}
		By the $0$-regularity of $g_{*}(\omega_{Y}\otimes{\mathcal{W}})$ with respect to ${\mathcal{H}}$, and the Mumford's theorem, see \cite[Chapter 5, Lemma 5.1]{FGA} or \cite[Theorem 1.8.5]{Laz04}, we have:
		\[
		H^{0}(Z,g_{*}(\omega_{Y}\otimes {\mathcal{W}}))\otimes H^{0}(Z,{\mathcal{H}}^{\ell})\rightarrow H^{0}(Z,g_{*}(\omega_{Y}\otimes {\mathcal{W}})\otimes{\mathcal{H}}^{\ell})
		\]
		is surjective for any $\ell\ge 0$. Since ${\mathcal{H}}^{r}={\mathcal{O}}_Z(1)$, we are done.
	\end{proof}

\subsection{Higher Dimension Case}\label{subsec:higher dim}

{In this subsection, we introduce the modified Grothendieck's Lefscehtz condition by Ravindra--Srinivas\cite{RS06} and prove the Noether--Lefschetz theorem for $\dim X\ge 3$. On one hand this is for the self-containedness of the paper, on the other hand, since the spectral divisors are rather special, all the proofs still hold when ${\operatorname{char}}(k)=p>0$ by assuming that $ H^i(X,{\mathcal{L}}^{-\ell})=0, \forall i<\dim X,\ell\in\mathbb{Z}_{>0}$. We point it out that by the above assumption of vanishing, Lemma \ref{lem: a vanishing theorem} also holds.}

First, we check the smoothness of a generic spectral variety. 
Instead, we regard the Hitchin base as the affine chart of the complete
linear system $|{\mathcal{W}}|$ where 
\[
\mathcal W=\pi^*\mathcal L^r\otimes\mathcal O_{Y/X}(r),
\]
on which the leading coefficient is nonzero and normalized to $1$.
The morphism induced by $|\mathcal W|$ restricts to an immersion on
$Y\setminus D_\infty=\operatorname{Tot}(\mathcal L)$. Hence this
linear system is separable on $\operatorname{Tot}(\mathcal L)$, and
Bertini's theorem implies that a general member is integral and smooth there.
Moreover, every monic spectral divisor is disjoint from $D_\infty$.

Let $X_{s}$ be a generic smooth spectral variety with $s\in \bm A$, our goal in this section is to show that the natural map $\pi^{*}:{\operatorname{Pic}}(X)\rightarrow {\operatorname{Pic}}(X_{s})$ is an isomorphism when $\dim X\ge 3$. Considering the following exact sequence:
\[
0\rightarrow \mathbb{Z}[D_{\infty}]\rightarrow{\operatorname{Pic}}(Y)\rightarrow{\operatorname{Pic}}(U)\cong\pi^{*}{\operatorname{Pic}}(X)\rightarrow 0,
\]
we only need to prove the following exact sequence for generic $s$:
\begin{equation}\label{eq:key exact seq}
	0\rightarrow \mathbb{Z}[D_{\infty}]\rightarrow{\operatorname{Pic}}(Y)\rightarrow{\operatorname{Pic}}(X_{s})\rightarrow 0.
\end{equation}
Since we also consider positive characteristics and also for the completeness of the paper, we slightly modified the proofs in \cite{RS06} to our cases. 

Let us denote the formal completion of $Y$ along $X_{s}$ by $\widehat{Y}_{s}$  and the $\ell$-th thickening of $X_s$ by $X_{s,\ell}$. Then $\widehat Y_s=\displaystyle\lim_{\rightarrow}X_{s,\ell}$ in the category of locally ringed spaces. We denote the defining ideal $X_{s}$ by ${\mathcal{I}}_s\cong{\mathcal{W}}^{-1}$. One has the exact sequence:
\[0\to {\mathcal{I}}_{s}^{\ell}\to {\mathcal{O}}_Y\to {\mathcal{O}}_{X_{s,\ell}}\to 0.\] 

\begin{proposition} \label{prop:inj of formal pic}If $\dim X\geq 3$, ${\operatorname{Pic}}(\widehat Y_s)\cong {\operatorname{Pic}}(X_s)$. If $\dim X=2$, the natural map ${\operatorname{Pic}}(\widehat{Y}_{s})\rightarrow {\operatorname{Pic}}(X_{s})$ is an injection.
\end{proposition}
\begin{proof}
	
	One has the following exact sequence:
	\[
	0\rightarrow{\mathcal{I}}^{m}_s/{\mathcal{I}}_s^{m+1}\rightarrow{\mathcal{O}}^{\times}_{X_{s,m+1}}\rightarrow {\mathcal{O}}^{\times}_{X_{s,m}}\rightarrow 0
	\]
	
	Then we have $$H^i(X_s,\sfrac{{\mathcal{I}}^m}{{\mathcal{I}}^{m+1}})\cong  H^i(X_s,{\mathcal{W}}^{-m}|_{X_s}) $$ Since ${\mathcal{W}}^{-1}|_{X_s}=\pi^*{\mathcal{L}}^{-r}\otimes {\mathcal{O}}_{Y/X}(-r)|_{X_s}\cong \pi_{s}^*{\mathcal{L}}^{-r}$,  if $\dim X\geq 3$, $H^i(X_s,\sfrac{{\mathcal{I}}^m}{{\mathcal{I}}^{m+1}})=0$ for $i=1,2$ and $m\geq 1$. Then by \cite[Expos\'e XI.1]{SGA2}, We get the isomorphism. If $\dim X=2$, similarly, we get $H^i(X_s,\sfrac{{\mathcal{I}}^m}{{\mathcal{I}}^{m+1}})=0$ for $i=1$, thus ${\operatorname{Pic}}(\widehat{Y}_{s})\rightarrow {\operatorname{Pic}}(X_{s})$ is an injection.
\end{proof}

Now let us recall the modified Grothendieck's Lefschetz conditions introduced in \cite[Definition 1]{RS06}, which is weaker than that in\cite[Expos\'e X.2]{SGA2} and fits into our cases well. It is this weaker Lefschetz condition that helps explain why we can not have the results as in \cite[Theorem 3.1.8]{SGA2}.
\begin{definition}[\cite{RS06}]\label{def:lef cond}
	Let $T$ be a smooth projective variety, and $D$ an effective divisor in $T$. We put $\widehat{T}$ the formal completion of $T$ along $D$. We say the pair $(T,D)$ satisfies the weak Lefschetz condition, denoted by $\text{Lef}^w(T,D)$ if for any open neighborhood $V$ of $D$, and any locally free coherent sheaf ${\mathcal{F}}$ on $V$, there is an open subset $V'\subset V$, such that the natural map: $H^{0}(V',{\mathcal{F}}|_{V'})\rightarrow\ H^{0}(\widehat{T},\widehat{\mathcal{F}})$ is an isomorphism. Here $\widehat{{\mathcal{F}}}$ is the completion of ${\mathcal{F}}$ along $D$.
	
	We say the pair $(T,D)$ satisfies the weak effective Lefschetz condition, which we denote by $\text{Leff}^w(T,D)$, if it satisfies $\text{Lef}^w(T,D)$ and in addition, for all locally free coherent sheaf ${\EuScript{F}}$ on $\widehat{T}$, there exist an open neighbourhood $V$ of $D$ and a locally free coherent sheaf ${\mathcal{F}}$ on $V$ such that $\widehat{{\mathcal{F}}}\cong {\EuScript{F}}$.
\end{definition}

Notice that for any open neighborhood of $X_{s}$ and any locally free coherent sheaf ${\mathcal{F}}$ on $V$, it can be extended to a reflexive sheaf on $Y$.
The following proposition shows that the pair $(Y,X_{s})$ satisfies the weak Lefschetz condition in Definition \ref{def:lef cond}:
\begin{proposition}\label{prop:mod lef cond}
For any reflexive sheaf ${\mathcal{N}}$ on $Y$ which is locally free in an open neighborhood of $X_s$, there exists an integer $m$ such that the natural map:
\[
H^{0}(Y,{\mathcal{N}}(mD_{\infty}))\rightarrow H^{0}(\widehat{Y}_{s},\widehat{{\mathcal{N}}})
\]
is an isomorphism.
\end{proposition}
\begin{proof}
Since $X_{s}\cap D_{\infty}=\emptyset$, $H^{0}(\widehat{Y}_{s},\widehat{{\mathcal{N}}})\cong H^{0}(\widehat{Y}_{s},\widehat{{\mathcal{N}}(\ell D_{\infty})})$ for any $\ell\in\mathbb{Z}$. Recall that, as a divisor, $X_s\simeq r(D_\infty+\pi^*L)$ and $X_s-iD_\infty$ is ample for $0<i<r$. In fact, $D_\infty+(1+\epsilon)\pi^*L$ is ample for any $\epsilon>0$, i.e., $r(D_\infty+\pi^*L)+m\pi^*L$ is ample for all $m>0$. This is because ${\mathcal{O}}_Y(r(D_\infty+\pi^*L)+m\pi^*L)$ is the pull back of relative ${\mathcal{O}}(1)$ on $\mathbf P({\mathcal{L}}^m\oplus\cdots\oplus{\mathcal{L}}^{r+m})/X$
and ${\mathcal{L}}^m\oplus\cdots\oplus {\mathcal{L}}^{r+m}$ is an ample vector bundle. 

The map in this proposition is induced by first considering the exact sequence for each thickening $X_{s,n}$ (the sequence is exact because we assume ${\mathcal{N}}$ to be locally free along $X_s$)
\[\xymatrix{ 0\ar[r] &{\mathcal{N}}(mD_\infty-nX_s)\ar[r]& {\mathcal{N}}(mD_\infty)\ar[r]^{t_n\ \ \ \ \ \ \ \ \ }&{\mathcal{N}}(mD_\infty)|_{X_{s,n}}\cong {\mathcal{N}}|_{X_{s,n}}\ar[r]& 0
},
\]
then taking the inverse limit $\displaystyle\lim_{\leftarrow}H^0(t_n)$ (see \cite[Chapter 8, Corollary 8.2.4]{FGA}). To prove the proposition, we have to show that for $n$ sufficiently large $H^0(t_n)$ is both injective and surjective (see \cite[8.2.5.2]{FGA} or \cite[Ch. 0, 13]{EGAIII-1}). By Serre duality, ${\operatorname{Ext}}^{i}({\mathcal{F}},\omega_Y)\cong H^{\dim Y-i}(Y,{\mathcal{F}})^{\vee}$ for any coherent sheaf on $Y$. Since ${\mathcal{N}}$ is reflexive, we have:
\[
H^0(Y,{\mathcal{N}}(mD_\infty-nX_s))=H^{\dim Y}(Y,\omega_Y\otimes {\mathcal{N}}^\vee(nX_s-mD_\infty))^\vee
\]
Via the spectral sequence:
\[
0\rightarrow H^{1}(Y,{\mathcal{E}} xt^0({\mathcal{F}},\omega_Y))\rightarrow{\operatorname{Ext}}^1({\mathcal{F}},\omega_Y)\rightarrow H^0(Y,{\mathcal{E}} xt^1({\mathcal{F}},\omega_Y))\rightarrow
\]
we have the inclusion (and also by the reflexivity of ${\mathcal{N}}$):
\[
H^1(Y,{\mathcal{N}}(mD_\infty-nX_s))\hookrightarrow H^{\dim Y-1}(Y,\omega_Y\otimes {\mathcal{N}}^\vee(nX_s-mD_\infty))^\vee
\]
Hence it suffices to show the vanishing of $H^{\dim Y-i}(Y,\omega_Y\otimes {\mathcal{N}}^\vee(nX_s-mD_\infty))=0, i=0,1$.
Since $nX_s-mD_\infty=m(X_s-D_\infty)+(n-m)X_s$, $(X_s-D_\infty)$ is ample and $X_s$ is nef, then by the Fujita vanishing theorem \cite[Theorem (1)]{Fuj83} (see also \cite[Remark 1.4.36]{Laz04}) for $m$ sufficiently large, and all $n>m$, we have the desired vanishing of cohomology, which completes the proof.
\end{proof}

\begin{proposition}\label{prop:ker of pic res}
We have the following exact sequence:
\[
0\rightarrow \mathbb{Z}D_{\infty}\rightarrow{\operatorname{Pic}}(Y)\rightarrow {\operatorname{Pic}}(\widehat{Y}_{s}).
\]
Since ${\operatorname{Pic}}(\widehat{Y}_{s})\rightarrow{\operatorname{Pic}}(X_{s})$ is injective, we also have:
\[
0\rightarrow \mathbb{Z}D_{\infty}\rightarrow{\operatorname{Pic}}(Y)\rightarrow {\operatorname{Pic}}(X_{s}).
\]
\end{proposition}
\begin{proof}
It is clear that $D_{\infty}$ is trivial when restricts to $\widehat{Y}_{s}$. Hence the map
\[
\mathbb Z[D_{\infty}]\rightarrow {\operatorname{Pic}}(Y)
   \rightarrow {\operatorname{Pic}}(\widehat{Y}_s)
\]
is zero. Put $U:=Y\setminus D_{\infty}$. By the localization sequence for the smooth variety $Y$, we have an
exact sequence
\[
\mathbb Z[D_{\infty}]
   \rightarrow {\operatorname{Pic}}(Y)
   \rightarrow {\operatorname{Pic}}(U)
   \rightarrow 0.
\]
It is therefore enough to prove that the restriction map
\[
{\operatorname{Pic}}(U)\rightarrow {\operatorname{Pic}}(X_s)
\]
is injective.

Let $\mathcal M$ be a line bundle on $U$ such that
\[
\mathcal M|_{X_s}\cong\mathcal O_{X_s}.
\]
Extending $\mathcal M$ to a rank-one reflexive sheaf on $Y$ and taking
its double dual, we may regard $\mathcal M$ as a line bundle on $Y$,
since $Y$ is smooth. By Proposition~\ref{prop:inj of formal pic}, the completion
$\widehat{\mathcal M}$ is trivial on $\widehat{Y}_s$. Choose a
nowhere-vanishing formal section
\[
\widehat{\sigma}\in
H^0(\widehat{Y}_s,\widehat{\mathcal M}).
\]
By Proposition~\ref{prop:mod lef cond}, there exists an integer $m$ such that the natural
map
\[
H^0\bigl(Y,\mathcal M(mD_{\infty})\bigr)
   \rightarrow
H^0\bigl(\widehat{Y}_s,\widehat{\mathcal M}\bigr)
\]
is an isomorphism. Thus $\widehat{\sigma}$ lifts to a global section $\sigma\in H^0\bigl(Y,\mathcal M(mD_{\infty})\bigr)$. Let $D:=(\sigma)_0$ be its zero divisor. Since the restriction of $\sigma$ to
$\widehat{Y}_s$ is invertible, we have ${\operatorname{Supp}}(D)\cap X_s=\varnothing$.

We claim that every prime divisor $\Gamma\subset Y$ disjoint from
$X_s$ is equal to $D_{\infty}$. Indeed, let
\[
g:Y\rightarrow
Z\subset\mathbb P\bigl(H^0(Y,\mathcal W)^\vee\bigr)
\]
be the morphism induced by the complete linear system $|\mathcal W|$,
and let $H_s$ be the hyperplane corresponding to the section defining
$X_s$. Then
$X_s=g^{-1}(Z\cap H_s)$. Hence, if $\Gamma\cap X_s=\varnothing$, then $g(\Gamma)\cap H_s=\varnothing$.
If $g(\Gamma)$ were positive-dimensional, then, being a projective
subvariety of a projective space, it would meet every hyperplane.
Equivalently, a positive-dimensional projective variety cannot be
contained in the affine open
$\mathbb P(H^0(Y,\mathcal W)^\vee)\setminus H_s$. Therefore
$g(\Gamma)$ is a point.

On the other hand, by the construction of $g$, its restriction
\[
g|_{Y\setminus D_{\infty}}:
Y\setminus D_{\infty}\rightarrow Z\setminus\{o\}
\]
is an isomorphism, while $D_{\infty}$ is contracted to the point $o$.
Thus every positive-dimensional subvariety contracted by $g$ is
contained in $D_{\infty}$. It follows that 
$\Gamma\subset D_{\infty}$. Since both $\Gamma$ and $D_{\infty}$ are prime divisors, we obtain
$\Gamma=D_{\infty}$, proving the claim.

Applying the claim to every prime component of $D$, we conclude that $D=\ell D_{\infty}$ for some integer $\ell\geq 0$. Therefore the section $\sigma$ induces
an isomorphism
\[
\mathcal M(mD_{\infty})
   \cong \mathcal O_Y(D)
   \cong \mathcal O_Y(\ell D_{\infty}),
\]
and hence $\mathcal M
   \cong \mathcal O_Y\bigl((\ell-m)D_{\infty}\bigr)$. Restricting to $U=Y\setminus D_{\infty}$ gives $\mathcal M|_U\cong\mathcal O_U$. Thus ${\operatorname{Pic}}(U)\to{\operatorname{Pic}}(X_s)$ is injective, and consequently we obtain
the exact sequence
\[
0\rightarrow \mathbb Z[D_{\infty}]
\rightarrow {\operatorname{Pic}}(Y)
\rightarrow {\operatorname{Pic}}(\widehat{Y}_s).
\]
\end{proof}

\begin{proposition}
For the pair $(Y,X_{s})$, $\operatorname{Leff}^w(Y,X_{s})$ holds.
\end{proposition}
\begin{proof}
By \cite[Chapter IV, Theorem 1.5]{Har70} (also see \cite[5.2.4]{EGAIII-1} and \cite[8.4.3]{FGA}), since ${\mathcal{W}}=\pi^{*}{\mathcal{L}}^{r}\otimes{\mathcal{O}}_{Y/X}(r)$ restricts to a very ample line bundle in an open neighbourhood of $X_{s}$, for any locally free coherent formal sheaf ${\EuScript{F}}$ on $\widehat{Y}_s$, we have the exact sequence:
\[{\mathcal{O}}_{\widehat Y_s}(-m_1)^{\oplus M_1}\xrightarrow{\widehat \phi}{\mathcal{O}}_{\widehat Y_s}(-m_2)^{\oplus M_2}\twoheadrightarrow {\EuScript{F}}\to 0.\]
For notation ease, we simply write $\widehat{{\mathcal{W}}^{m}|_{Y_{s}}}$ by ${\mathcal{O}}_{\widehat Y_s}(m)$.

By Corollary \ref{prop:mod lef cond}, $\operatorname{Lef}^w(Y,X_s)$ holds and $$\widehat \phi\in {\mathcal{H}} om({\mathcal{O}}_{\widehat Y_s}(-m_1)^{\oplus M_1},{\mathcal{O}}_{\widehat Y_s}(-m_2)^{\oplus M_2})\cong \widehat{{\mathcal{H}} om}({\mathcal{O}}_{Y}(-m_1)^{\oplus M_1},{\mathcal{O}}_{Y}(-m_2)^{\oplus M_2})$$ is algebrizable by $\phi\in\Gamma(Y, {\mathcal{H}} om_{Y}({\mathcal{O}}(-m_1)^{\oplus M_1},{\mathcal{O}}(-m_2)^{\oplus M_2}\otimes{\mathcal{O}}_Y(m_3 D_\infty)))$. Then we have $\widehat{\operatorname{CoKer}(\phi)}\cong {\EuScript{F}}$. For any $y\in X_s$, we have ${\mathcal{I}}_{s}\subset \mathfrak m_y$. Then completion along ${\mathcal{I}}_{s}$ and then completion along ${\widehat{\mathfrak m_y}}^{{\mathcal{I}}_s}$, is equal to completing directly at $\mathfrak m_y$. 

This means that ${\widehat{\operatorname{Coker}(\phi)}}^{\mathfrak m_y}$ is locally free at each point $y\in X_s$. By faithful flatness of the completion along a maximal ideal, $\operatorname{Coker}(\phi)$ is locally free after being localized at each closed point $y\in X_s$. Thus  $\operatorname{Coker}(\phi)$ is locally free over a neighbourhood $U$ of $X_s$.\end{proof}

\begin{theorem}\label{thm:eff lef for formal line bundle}
We have the following exact sequence:
\[
0\rightarrow\mathbb{Z}D_{\infty}\to {\operatorname{Pic}}(Y)\xrightarrow{c} {\operatorname{Pic}}(\widehat Y_s)\to 0.
\]
in particular,  if $\dim X\geq 3$,  we have $\pi_{s}^{*}:{\operatorname{Pic}}(X)\rightarrow {\operatorname{Pic}}(X_{s})$ is an isomorphism provided that $X_{s}$ is smooth.
\end{theorem}
\begin{proof}
By Proposition \ref{prop:inj of formal pic} and Proposition \ref{prop:ker of pic res}, we only need to show that $c$ is surjective. For any formal line bundle ${\EuScript{M}}$ on $\widehat Y_s$, by $\operatorname{Leff}^w(Y,X_s)$ there is an open neighbourhood $V$ of $X_s$ and an invertible sheaf ${\mathcal{M}}_V$ on $V$ such that $\widehat {\mathcal{M}}_V\cong {\EuScript{M}}$. One check that ${\mathcal{M}}_V$ can always be extend to a line bundle ${\mathcal{M}}$ over $Y$ provided $Y$ is smooth (see \cite[Corollary 3.4]{HeinLec}, extending ${\mathcal{M}}_V$ to a coherent sheaf and taking the double dual, which is relexive of rank 1, hence a line bundle). Thus $c({\mathcal{M}})=\widehat{{\mathcal{M}}}\cong{\EuScript{M}}$ and $c$ is surjective. 

If $\dim X\ge 3$, ${\operatorname{Pic}}(\widehat{Y}_{s})\rightarrow {\operatorname{Pic}}(X_{s})$ is an isomorphism, thus we conclude that $\pi_{s}^{*}:{\operatorname{Pic}}(X)\rightarrow {\operatorname{Pic}}(X_{s})$ is an isomorphism.	\end{proof}

\subsection{Surface Case} 

In this section, we prove a similar result when $\dim X=2$ and $\mathrm {char}(k)=0$. We
 follow \cite{RS09}. But to make the paper self-contained, we adapt their methods to the setting of spectral varieties to use the property that $g_{*}(\omega_{Y}\otimes\pi^{*}{\mathcal{L}}^{\otimes r}\otimes {\mathcal{O}}_{Y/X}(r))$ is $0$-regular.

\begin{notation} $\underline{{\operatorname{Pic}}}_{X/S}^{\natural}$ means the relative Picard functor,  $\underline{{\operatorname{Pic}}}_{X/S}^{\natural,{\operatorname{\acute{e}t}}}$ the \'etale sheafification of $\underline{{\operatorname{Pic}}}_{X/S}^{\natural}$ and $\underline{{\operatorname{Pic}}}_{X/S}$ means the scheme it is represented by (if representable). If $S={\operatorname{Spec}} k$, we may omit $S$. ${\operatorname{Pic}}(X)$ means the Picard group of $X$, which is isomorphic to $H^1(X,{\mathcal{O}}_X^\times)$.     
\end{notation}

As before, let ${\mathcal{L}}$ be an ample line bundle on $X$ and we put $Y=\mathbf P({\mathcal{L}}^{\vee}\oplus {\mathcal{O}})$ with $\pi:Y\to X$ the natural projection. Since ${\mathcal{L}}^{r}$ is very ample, we have $Z=g(Y)$ consisting a unique singularity $o=\phi(D_{\infty})$ where $g$ is induced by the complete linear system $|\pi^*{\mathcal{L}}^{r}\otimes {\mathcal{O}}_{Y/X}(r)|=|{\mathcal{W}}|$. 

Let $X_s\subset Y$ be a smooth spectral surface defined by $s\in |\pi^*{\mathcal{L}}^{r}\otimes {\mathcal{O}}_{Y/X}(r)|$, which does not intersect with $D_\infty$. One can always view $X_s$ as a very ample divisor in $Z$ defined by a global section of ${\mathcal{O}}_Z(1)$. The main goal of this subsection is to prove the following theorem:
\begin{theorem}\label{thm:iso for surfaces}[${\operatorname{char}} =0, \dim X=2$ case] Assuming that $r\ge 4$, then for a very general $s\in |\pi^*{\mathcal{L}}^{r}\otimes {\mathcal{O}}_{Y/X}(r)|$ the map $\pi^{*}:\underline{{\operatorname{Pic}}}_X\cong \underline{{\operatorname{Pic}}}_{X_s}$ is an isomorphism.
\end{theorem}
We start with the following lemma on formal Picard groups.
\begin{lemma}
For every smooth spectral surface $X_s$, the natural maps give
\[
{\operatorname{Pic}}(X)\xrightarrow{\sim}{\operatorname{Pic}}(\widehat Y_s)
\hookrightarrow{\operatorname{Pic}}(X_s).
\]
\end{lemma}

\begin{proof}
By Theorem~\ref{thm:eff lef for formal line bundle},
\[
{\operatorname{Pic}}(\widehat Y_s)
\cong {\operatorname{Pic}}(Y)/\mathbb Z[D_\infty]
\cong {\operatorname{Pic}}(X).
\]
The injectivity of the second map follows from
Proposition~\ref{prop:inj of formal pic}.
\end{proof}

Let ${\mathcal{Y}}$ contained in $\mathbf P_Y:=Y\times \mathbf P(H^{0}(Y,\pi^{*}{\mathcal{L}}^{r}\otimes{\mathcal{O}}(r))^\vee)$ be the universal family of divisors parameterized by $\mathbf P(H^{0}(Y,\pi^{*}{\mathcal{L}}^{r}\otimes{\mathcal{O}}(r))^\vee)$. We put $p=\bm p_Y|_{{\mathcal{Y}}}:{\mathcal{Y}}\to Y$ and $q:=\bm p_{\mathbf P}|_{{\mathcal{Y}}}:{\mathcal{Y}}\to \mathbf P(H^{0}(Y,\pi^{*}{\mathcal{L}}^{r}\otimes{\mathcal{O}}(r))^\vee)$:
\[
\begin{tikzcd}
&{\mathcal{Y}}\ar[ldd, "p"']\ar[rdd,"q"]& {\qquad}{\qquad}\ldots\supset {\mathcal{Y}}_{s,\ell}\ldots \supset {\mathcal{Y}}_s\ar[rdd]\\
\\
Y&& {\small\mathbf P(H^{0}(Y,\pi^{*}{\mathcal{L}}^{r}\otimes{\mathcal{O}}(r))^\vee)}&\ldots\supset  \small{{\operatorname{Spec}}(\frac{{\mathcal{O}}}{\mathfrak m_{s}^{\ell-1}})}\ldots\ni s
\end{tikzcd}
\]
By the construction, ${\mathcal{Y}}\subset \mathbf P_Y$ is relatively very ample to $p_Y$. Thus $R^i{p_Y}_*{\mathcal{O}}_{\mathbf P_Y}(-{\mathcal{Y}})=0$ for all $i\geq 1$. For a closed point $s\in\bm A\subset \mathbf P(H^{0}(Y,\pi^{*}{\mathcal{L}}^{r}\otimes{\mathcal{O}}_{Y/X}(r))^\vee)$, we write ${\mathcal{Y}}_{s}:=q^{-1}(s)$  which is isomorphic to $X_{s}$ under the projection $p$. 

Similar as before, we put ${\mathcal{Y}}_{s,\ell}$ the $\ell$-th thickening of ${\mathcal{Y}}_{s}$, and $\widehat{{\mathcal{Y}}}_{s}$ the formal completion along ${\mathcal{Y}}_{s}$. We denote the maximal ideal of $s$ by $\mathfrak m_{s}$, thus the defining ideal of ${\mathcal{Y}}_{s}$ is $q^{*}\mathfrak m_{s}=\mathfrak m_{s}{\mathcal{O}}_{{\mathcal{Y}}}$. We denote the defining ideal of $X_{s}$ in $Y$ by ${\mathcal{I}}_{s}\cong {\mathcal{W}}^{-1}$. Then we have $p^{*}{\mathcal{I}}_{s}\subset q^{*}\mathfrak m_{s}$. In particular, we have a map of infinitesimal thickenings $ {\mathcal{Y}}_{s,\ell}\to X_{s,\ell}$ induced by $p$. As a result, we have the following commutative diagram:
\[
\xymatrix{&{\operatorname{Pic}}(\widehat{Y}_{s})\ar[dd]\ar[r]&{\operatorname{Pic}}(X_{s,\ell})\ar[dd]\ar[rd]\\
{\operatorname{Pic}}(Y)\ar[rd]\ar[ru]&&&{\operatorname{Pic}}(X_{s})\\
&{\operatorname{Pic}}(\widehat{{\mathcal{Y}}}_{s})\ar[r]&{\operatorname{Pic}}({\mathcal{Y}}_{s,\ell})\ar[ru]
}.
\]

Following \cite{RS09}, we introduce the infinitesimal Noether--Lefschetz condition.
\begin{definition}
We say the pair $(Y,X_{s})$ satisfies the $\ell$-th infinitesimal Noether--Lefschetz condition, denoted by $\text{INL}_{\ell}$ if the following map is an isomorphism:
\[
\text{Image}({\operatorname{Pic}}(X_{s,\ell})\rightarrow{\operatorname{Pic}}(X_{s}))\rightarrow\text{Image}({\operatorname{Pic}}({\mathcal{Y}}_{s,\ell})\rightarrow{\operatorname{Pic}}(X_{s})).
\]
Similarly, we say the pair $(Y,X_{s})$ satisfies the formal Noether--Lefschetz condition, denoted by $\text{FNL}$ if the following map is an isomorphism:
\[
\text{Image}({\operatorname{Pic}}(\widehat{Y}_{s})\rightarrow{\operatorname{Pic}}(X_{s}))\rightarrow\text{Image}({\operatorname{Pic}}(\widehat{{\mathcal{Y}}}_{s})\rightarrow{\operatorname{Pic}}(X_{s})).
\]
\end{definition}

We restrict the universal family to get a smooth family, without causing ambiguity, we still denote it by 
\begin{equation}\label{eq:smooth family}
q:{\mathcal{Y}}_{U}\rightarrow U 
\end{equation}
where $q$ is smooth with connected fibers. Since $q$ is flat, projective with integral geometric fibers, by \cite[Theorem 9.4.8]{FGA}, $\underline{{\operatorname{Pic}}}_{q}^{\natural,{\operatorname{\acute{e}t}}}$, the \'etale sheaf associated with the relative Picard functor $\underline{{\operatorname{Pic}}}^{\natural}_{q}$, is representable by a scheme $\underline{{\operatorname{Pic}}}_q$ which is separated and locally of finite type over $U$.

By the base change property of relative Picard scheme, for any closed point $\xi\in U$, $\underline{{\operatorname{Pic}}}_{q}|_{\xi}=\underline{{\operatorname{Pic}}}_{X_\xi/\xi}$ is smooth. It means that each fiber of $\underline{{\operatorname{Pic}}}_{q}\rightarrow U$ is smooth.

Let $\operatorname{HP}$ be the set of Hilbert polynomials of line bundles on ${\mathcal{Y}}_{U}$. $\Phi\subset \operatorname{HP}$ be a finite subset. Denote by $\underline{{\operatorname{Pic}}}_q^\Phi\subset \underline{{\operatorname{Pic}}}_q$ be the components with Hilbert polynomials in $\Phi$, then $\underline{{\operatorname{Pic}}}_q^\Phi$ is of finite type over $U$. Denote by $U^\Phi\subset U$ the open subset of $U$ on which $\underline{{\operatorname{Pic}}}_q^\Phi$ is flat (hence smooth). This subset is non-empty because it contains the generic point of $U$. The intersection of $U^\Phi$ for $\Phi$ covering $\operatorname{HP}$ is a very general subset $\bm V$ of $U$.  

\begin{lemma}
If for any $s\in \bm V\subset U$ and any $m_1>m_2\in {\mathbb {Z}}_{>0}$ , we have $\pi^{*}:{\operatorname{Pic}}(\widehat{{\mathcal{Y}}}_{s,m_1})\rightarrow{\operatorname{Pic}}(\widehat{{\mathcal{Y}}}_{s,m_2})$ is surjective.
\end{lemma}
\begin{proof}
For $s\in \bm V\subset U$, we have $\underline{{\operatorname{Pic}}}_{q}$ is smooth over $s$. By the smoothness, we have the following surjective morphisms, \begin{footnotesize}
	$${\operatorname{Hom}}({\operatorname{Spec}}\widehat{{\mathcal{O}}}_{U,s}, \underline{{\operatorname{Pic}}}_{q})\twoheadrightarrow {\operatorname{Hom}}({\operatorname{Spec}}({\mathcal{O}}_{U,s}/\mathfrak m_{s}^{m_1}),\underline{{\operatorname{Pic}}}_{q})\twoheadrightarrow {\operatorname{Hom}}({\operatorname{Spec}}({\mathcal{O}}_{U,s}/\mathfrak m_{s}^{m_2}),\underline{{\operatorname{Pic}}}_{q}).$$
	
\end{footnotesize} Consider the following exact sequence, coming from low-degree terms in the Leray spectral sequence \cite[Chapter III, Theorem 1.18]{Mil80} for ${\mathbb {G}}_m$ relative to $q_{S}$ with $S= {\operatorname{Spec}}\widehat{{\mathcal{O}}}_{U,s}$, (see also \cite[Chapter 9, (9.2.11.5)]{FGA})
$$0\to {\operatorname{Pic}}(S)\to {\operatorname{Pic}}({\mathcal{Y}}_S)\to \underline{{\operatorname{Pic}}}^{\natural,{\operatorname{\acute{e}t}}}_{q}(S)\to H^2(S_{\operatorname{\acute{e}t}},{q_S}_*{\mathbb {G}}_m).$$
We have ${q_S}_*{\mathbb {G}}_m={\mathbb {G}}_{m,S}$. Thus the previous sequence becomes
$$0\to {\operatorname{Pic}}(S)\to {\operatorname{Pic}}({\mathcal{Y}}_S)\to \underline{{\operatorname{Pic}}}^{\natural,{\operatorname{\acute{e}t}}}_{q}(S)\to H^2(S_{{\operatorname{\acute{e}t}}},{\mathbb {G}}_m).$$
By \cite[Chapter III, Theorem 3.9]{Mil80}, $S$ is strictly local, then $H^2(S_{\operatorname{\acute{e}t}},{\mathbb {G}}_m)=0$, we have $\underline{{\operatorname{Pic}}}_q^{\natural,{\operatorname{\acute{e}t}}}(S)\cong {\operatorname{Pic}}({\mathcal{Y}}_S)$ for $S= {\operatorname{Spec}}\widehat{{\mathcal{O}}}_{U,s}$. Since ${{\mathcal{O}}}_{U,s}/\mathfrak m_{s}^\ell$ for $\ell\in{\mathbb {Z}}_{>0}$ is an Artin local ring with algebraically closed residue field, we have $\underline{{\operatorname{Pic}}}_q^{\natural,{\operatorname{\acute{e}t}}}(S)\cong {\operatorname{Pic}}({\mathcal{Y}}_S)$ for $S={\operatorname{Spec}}{{\mathcal{O}}}_{U,s}/\mathfrak m_{s}^\ell$ for $\ell\in{\mathbb {Z}}_{>0}$. And the previous surjections show that $\pi^{*}:{\operatorname{Pic}}({\mathcal{Y}}_{s,m_1})\rightarrow{\operatorname{Pic}}({\mathcal{Y}}_{s,m_2})$ is surjective for any $s\in \bm V$ and $m_1>m_2\in {\mathbb {Z}}_{>0}$.
\end{proof}

	{Combining with Proposition \ref{prop:ker of pic res}, }it is easy to deduce the following corollary,
	\begin{corollary}\label{cor:FNL holds then}
For any closed point $s\in \bm V$, if the pair $(Y,X_{s})$ satisfies the $\text{FNL}$, then we have $\pi^{*}:{\operatorname{Pic}}(X)\rightarrow{\operatorname{Pic}}(X_{s})$ is an isomorphism.
\end{corollary}

We here collect cohomological results which follows from the fact that $p:{\mathcal{Y}}\rightarrow Y$ is a projective bundle. 
\begin{lemma}\cite[Lemma 1 and  \S 2.1]{RS09}\label{lem:direct image comp}
\begin{enumerate}
	\item $Rp_{*}(q^{*}\mathfrak m_{s})\cong {\mathcal{I}}_{s}$;
	\item $p_{*}(q^*\mathfrak m_{s}^{\ell})={\mathcal{I}}_{s}^{\ell}$ for $\ell\ge 1$;
	\item $0\rightarrow{\mathcal{O}}_{X_{s,\ell}}\rightarrow p_{*}{\mathcal{O}}_{{\mathcal{Y}}_{s,\ell}}\rightarrow R^{1}p_{*}(q^{*}\mathfrak m_{s}^{\ell})\rightarrow 0$ .
	\item $R^jp_*q^*\mathfrak m^\ell_s=0$ for $j\geq 2$ and $\ell\in{\mathbb {Z}}_{>0}$. 
\end{enumerate}
\end{lemma}

The following proposition is crucial for us to prove the $\text{FNL}$ for $s\in \bm V\subset U$.
\begin{proposition}\label{prop:inj for every step}
\[
H^{2}(X_{s}, {\mathcal{I}}_{s}^{\ell}/{\mathcal{I}}_{s}^{\ell+1})\rightarrow H^{2}({\mathcal{Y}}_{s},q^{*}\mathfrak m_{s}^{\ell}/q^{*}\mathfrak m_{s}^{\ell+1})
\]
is injective for all $\ell>0$.
\end{proposition}

\begin{proof}
Consider the following commutative diagram:
\[
\xymatrix{
	0\ar[r]&{\mathcal{I}}_{s}^{\ell+1}\ar[r]\ar[d]&{\mathcal{I}}_{s}^{\ell}\ar[r]\ar[d]&{\mathcal{I}}_{s}^{\ell}/{\mathcal{I}}_{s}^{\ell+1}\ar[r]\ar[d]&0\\
	0\ar[r]&q^{*}\mathfrak m_{s}^{\ell+1}
	\ar[r]&q^{*}\mathfrak m_{s}^{\ell}\ar[r]&q^{*}\mathfrak m^{\ell}/q^{*}\mathfrak m_{s}^{\ell+1}\ar[r]&0}.
\]
Recall that ${\mathcal{I}}_{s}\cong{\mathcal{W}}^{-1}=\pi^{*}{\mathcal{L}}^{-r}\otimes{\mathcal{O}}_{Y/X}(-r)$. By Lemma \ref{lem: a vanishing theorem},  $H^{i}(Y, {\mathcal{I}}_{s}^{\ell})=0$ for all $\ell>0, i=1,2$.
Thus we have:
\begin{footnotesize}
	\[
	\xymatrix{
		0\ar[r]&H^{2}(X_{s},{\mathcal{I}}_{s}^{\ell}/{\mathcal{I}}_{s}^{\ell+1})\ar[r]\ar[d]&H^{3}(Y,{\mathcal{I}}_{s}^{\ell+1})\ar[r]\ar[d]&H^{3}(Y,{\mathcal{I}}_{s}^{\ell})\ar[d]\\
		\ar[r]&H^{2}({\mathcal{Y}}_{s},q^{*}\mathfrak m_{s}^{\ell}/q^{*}\mathfrak m_{s}^{\ell+1})\ar[r]&H^{3}({\mathcal{Y}},q^{*}\mathfrak m_{s}^{\ell+1})\ar[r]&H^{3}({\mathcal{Y}},q^{*}\mathfrak m_{s}^{\ell})
	}.
	\]
\end{footnotesize}

\noindent Then we only need to prove that $H^{3}(Y,{\mathcal{I}}_{s}^{\ell})\rightarrow H^{3}({\mathcal{Y}},q^{*}\mathfrak m_{s}^{\ell})$ is injective for all $\ell>0$. Since $p_{*}(q^{*}m_{s}^{\ell})={\mathcal{I}}^{\ell}_s$, and considering that the Leray spectral sequence for $p:{\mathcal{Y}}\rightarrow Y$, $H^{3}(Y,{\mathcal{I}}_{s}^{\ell})\rightarrow H^{3}({\mathcal{Y}},q^{*}\mathfrak m_{s}^{\ell})$ is the map $E_{2}^{3,0}\rightarrow H^{3}({\mathcal{Y}},q^{*}\mathfrak m^{\ell})$, thus to show the injectivity, we only need to show the following differential vanishes for all $\ell>0$:
\[
H^{1}(Y,R^{1}p_{*}q^{*}\mathfrak m_{s}^{\ell})=E_{2}^{1,1}\rightarrow E_{2}^{3,0}=H^{3}(Y,{\mathcal{I}}_{s}^{\ell})
\]

By Ravindra--Srinivas \cite[\S 2.1, \S 2.2]{RS09},it suffices to prove that the evaluation map:
\[
 H^{0}(Y,\omega_{Y}\otimes{\mathcal{W}})\otimes {\operatorname {Sym}}^{(\ell -1)}H^{0}(Y,{\mathcal{W}})\rightarrow H^{0}(Y,\omega_{Y}\otimes{\mathcal{W}}^{\otimes \ell})
\]
is surjective for all $\ell>0$. 
By Proposition~\ref{prop:global gen}, $g_*(\omega_Y\otimes{\mathcal{W}})$ is $0$-regular
with respect to the ample and globally generated line bundle
${\mathcal{H}}$. Iterating Mumford's theorem
\cite[Theorem~1.8.5(ii),(iii)]{Laz04} gives a surjection
\[
H^0(Z,g_*(\omega_Y\otimes{\mathcal{W}}))
\otimes H^0(Z,{\mathcal{H}})^{\otimes r(\ell-1)}
\twoheadrightarrow
H^0(Z,g_*(\omega_Y\otimes{\mathcal{W}})\otimes{\mathcal{H}}^{r(\ell-1)})
\]
Since $g^*{\mathcal{H}}^r\cong{\mathcal{W}}$, we have 
$$H^0(Z,g_*(\omega_Y\otimes{\mathcal{W}})\otimes{\mathcal{H}}^{r(\ell-1)})\cong H^0(Y,\omega_Y\otimes{\mathcal{W}}^\ell).$$
Since $g^*{\mathcal{H}}^r\cong{\mathcal{W}}$, this surjection factors through the
preceding symmetric-power multiplication map. Hence that
map is surjective.\qed
\end{proof}

\begin{remark}
    In fact, after diagram-chasing and Serre duality, the map $$ H^{0}(Y,\omega_{Y}\otimes{\mathcal{W}})\otimes {\operatorname {Sym}}^{(\ell -1)}H^{0}(Y,{\mathcal{W}})\rightarrow H^{0}(Y,\omega_{Y}\otimes{\mathcal{W}}^{\otimes \ell})$$ should be a shifted (up to some constants depends on $\ell$) of the standard evaluation map. Thus in positive characteristic cases, the map can be non-surjective because some of the coefficients may be divided by the characteristic.
\end{remark}

\begin{proposition}\label{prop:FNL}
When $r>3$, then for any closed point $s\in \bm V$ (i.e. over $s$ the relative Picard variety is smooth), we have the pair $(Y,X_{s})$ satisfies $\text{INL}_{m}$ for all $m>0$ and then it satisfies $\text{FNL}$.
\end{proposition}

\begin{proof}
For any $m>0$, we have:
\[    \xymatrix{
	0\ar[r]& {\mathcal{I}}_{s}^{m}/{\mathcal{I}}_{s}^{m+1}\ar[r]\ar[d]|{\beta_{m+1}}&{\mathcal{O}}^{\times}_{X_{s,m+1}}\ar[r]\ar[d]|{p^{\flat}_{m+1}}&{\mathcal{O}}^{\times}_{X_{s,m}}\ar[r]\ar[d]|{p^{\flat}_{m}}&0\\
	0\ar[r]& q^*\mathfrak m_{s}^{m}/q^*\mathfrak m_{s}^{m+1}\ar[r]&{\mathcal{O}}^{\times}_{{\mathcal{Y}}_{s,m+1}}\ar[r]&{\mathcal{O}}^{\times}_{{\mathcal{Y}}_{s,m}}\ar[r]&0
}.\]
Since $H^{1}(X_{s},{\mathcal{I}}_{s}^{m}/{\mathcal{I}}_{s}^{m+1})=0$, we have the following exact sequence:
\begin{footnotesize}\begin{equation}
		\xymatrix@C=3ex{
			0\ar[r]&{\operatorname{Pic}}(X_{s,m+1})\ar[r]\ar[d]&{\operatorname{Pic}}(X_{s,m})\ar[r]\ar[d]& H^{2}(X_{s}, {\mathcal{I}}_{s}^{m}/{\mathcal{I}}_{s}^{m+1})\ar[d]\\
			0\ar[r]&{\operatorname{Pic}}({\mathcal{Y}}_{s,m+1})/H^{1}({\mathcal{Y}}_{s},q^{*}\mathfrak m_{s}^{m}/q^{*}\mathfrak m_{s}^{m+1})\ar[r]&{\operatorname{Pic}}({\mathcal{Y}}_{s,m})\ar[r]& H^{2}({\mathcal{Y}}_{s},q^{*}\mathfrak m_{s}^{m}/q^{*}\mathfrak m_{s}^{m+1})
		}
\end{equation}\end{footnotesize}

\noindent	Notice that $s\in \bm V$, we have the surjection ${\operatorname{Pic}}({\mathcal{Y}}_{s,m+1})\twoheadrightarrow {\operatorname{Pic}}({\mathcal{Y}}_{s,m})$, which implies that the map ${\operatorname{Pic}}({\mathcal{Y}}_{s,m})\rightarrow H^{2}({\mathcal{Y}}_{s},q^{*}\mathfrak m_{s}^{m}/q^{*}\mathfrak m_{s}^{m+1})$ is a zero map. By the injectivity of $$H^{2}(X_{s}, {\mathcal{I}}_{s}^{m}/{\mathcal{I}}_{s}^{m+1})\hookrightarrow H^{2}({\mathcal{Y}}_{s},q^{*}\mathfrak m_{s}^{m}/q^{*}\mathfrak m_{s}^{m+1})$$ for all $m>0$, we know that ${\operatorname{Pic}}(X_{s,m})\to H^{2}(X_{s}, {\mathcal{I}}_{s}^{m}/{\mathcal{I}}_{s}^{m+1})$ is a zero map. Thus ${\operatorname{Pic}}(X_{s,m+1})\to{\operatorname{Pic}}(X_{s,m})$ is surjective. As a result we have $\text{Image}({\operatorname{Pic}}(X_{s,m})\rightarrow{\operatorname{Pic}}(X_{s}))={\operatorname{Pic}}(X_{s})$ and the $\text{FNL}$ holds.
\end{proof}

To conclude, we have:
\begin{theorem}\label{thm:NL for surfaces}
When $r\ge 4$, for $s$ in the very general subset $\bm V\subset \bm A$, we have the isomorphism $\pi^{*}:{\operatorname{Pic}}(X)\rightarrow{\operatorname{Pic}}(X_{s})$.
\end{theorem}
Notice that, if $X=\mathbb{P}^{2}$, ${\mathcal{L}}={\mathcal{O}}(1)$, and $r=3$, then $X_{s}$ is a cubic surface whose Picard number is different from the Picard number of $X$. However, if the canonical line bundle of $X$ is sufficiently ample, we expect that the theorem holds for smaller $r$.

\paragraph{\textbf{"Bigness" is Necessary}}
For simplicity, we give an example using a complex elliptic surface. Now we assume $f:X\rightarrow \mathbb{P}^1$ is a relatively minimal non-isotrivial elliptic surface with a section and irreducible singular fibers.  We assume that $e:=\chi({\mathcal{O}}_X)\ge 4$ and fix the rank $r\ge 2$. We put
\[
\Sigma:=\{x\in C| f^{-1}(x) \;\text{is singular}\}.
\]
Then by Kodaira's canonical bundle formula, 
\[
(R^1f_*{\mathcal{O}}_X)^\vee\cong{\mathcal{O}}_{\mathbb P^1}(e),
\qquad
\omega_X\cong f^*{\mathcal{O}}_{\mathbb P^1}(e-2).
\]
Thus $\omega_X$ is nef but not big: it is pulled back from a curve and has self-intersection zero.

Since $f_*{\mathcal{O}}_X={\mathcal{O}}_{\mathbb P^1}$, the projection formula identifies
\[
H^0(X,\omega_X^i)
\cong
H^0\bigl(\mathbb P^1,{\mathcal{O}}_{\mathbb P^1}(i(e-2))\bigr).
\]
Consequently, a degree-$r$ spectral surface for $\omega_X$ is obtained by base change from the corresponding spectral curve:
\[
X_s\cong X\times_{\mathbb P^1}C_s,
\]
where $C_s$ lies in the total space of
${\mathcal{O}}_{\mathbb P^1}(e-2)$.

For general $s$, the curve $C_s$ is smooth and its projection to $\mathbb P^1$ is unramified over $\Sigma$. The base change is then smooth above $\Sigma$ because this projection is etale there, and it is smooth away from $\Sigma$ because $f$ is smooth there. Hence $X_s$ is smooth. On the other hand,
\[
H^1(X_s,{\mathcal{O}}_{X_s})
\cong
\bigoplus_{i=0}^{r-1}H^1(X,\omega_X^{-i}).
\]
The Leray spectral sequence gives an injection
\[
H^1\bigl(\mathbb P^1,{\mathcal{O}}_{\mathbb P^1}(2-e)\bigr)
\hookrightarrow H^1(X,\omega_X^{-1}),
\]
and the source has dimension $e-3>0$. Since $r\geq2$,
\[
h^1(X_s,{\mathcal{O}}_{X_s})>h^1(X,{\mathcal{O}}_X).
\]
Therefore $\dim\underline{{\operatorname{Pic}}}_{X_s}^0>
\dim\underline{{\operatorname{Pic}}}_X^0$. In particular, without the bigness assumption, pullback on Picard schemes can fail to be an isomorphism even when the spectral surface is smooth.

\section{Generic Fibers of Hitchin Systems of Higgs sheaves}
In the following sections, we take $k=\mathbb{C}$, and let $X$ be a smooth  complex projective surface.  We consider the moduli space of Gieseker (semi-)stable Higgs sheaves on $X$ and the associated Hitchin maps. To be more precise, we use the Auslander--Buchsbaum formula from homological algebra to identity generic fibers of such Hitchin systems with products of Picard schemes and Hilbert schemes of points of general spectral surfaces. Then by the previous calculation of Picard groups of very general spectral surfaces to obtain a necessary and sufficient condition for the non-emptyness of generic fibers. 

\subsection{Moduli space of Higgs sheaves and the $\mathbb G_m$-action}
Following \cite{GSY20L}, we start from more general Higgs sheaves than that in \cite{TT18}.
\begin{definition}
An ${\mathcal{L}}$-valued Higgs sheaf of rank $r$ on $X$ consists of $(\mathcal E,\theta)$ with $\mathcal E$ being a torsion free sheaf on $X$ of generic rank $r$, $\theta:\mathcal E\rightarrow \mathcal E\otimes \mathcal L$ being an $\mathcal O_{X}$-linear map. If no ambiguity caused, we will drop "${\mathcal{L}}$-valued" for simplicity. When we want to emphasize that $\mathcal E$ is locally free, we call $(\mathcal E,\theta)$ a Higgs bundle. 
\end{definition}

\begin{remark}Though $\mathcal L$-valued Higgs sheaves can be defined (and also of its importance) for general line bundles, to detect generic fibers of $h^{\text{tf}}$, we here assume $\mathcal L$ to be very ample. 
\end{remark}

\begin{remark} In higher dimensional case, Simpson in \cite{Sim92} introduce the concept of Higgs bundles valued in the cotangent bundle $\Omega_X^1$. For vector bundle valued Higgs sheaves, we refer the readers to \cite{CN20,GGPN24,HL24,HLM24,SS24} and references therein. 
\end{remark}

Let us recall the stability conditions. As mentioned in the introduction,  we choose the polarization $H$ such that $H$ is numerically equivalent to $\lambda c_1({\mathcal{L}})$ for some $\lambda>0$. For a torsion-free sheaf
$\mathcal E$, write
\[
\mu_H(\mathcal E)
=
\frac{c_1(\mathcal E)\cdot H}{\operatorname{rk}(\mathcal E)}
\]
and let $p_{H,\mathcal E}$ denote its reduced Hilbert polynomial.
Slope and Gieseker stability for a Higgs sheaf
$(\mathcal E,\theta)$ are defined by testing only
$\theta$-invariant proper subsheaves. As usual,
\[
\text{slope stable}
\Longrightarrow
\text{Gieseker stable}
\Longrightarrow
\text{Gieseker semistable}
\Longrightarrow
\text{slope semistable}.
\]
If $\gcd(r,c_1(\mathcal E)\cdot H)=1$, all four conditions
coincide; if $\gcd(r,\chi(\mathcal E))=1$, Gieseker
semistability coincides with Gieseker stability.

With the help of the Geometric Invariant Theory, the coarse moduli spaces
of Higgs sheaves with fixed Chern classes are constructed. In the curve case, Nitsure\cite{N91} constructed the moduli space of (semi-)stable $\mathcal L$-valued Higgs sheaves over an algebraically closed field and showed the properness of Hitchin maps. In more general setting (smooth projective scheme over a universally Japanese ring with some boundedness assumption), Yokogawa\cite{Yo91}(see also \cite{Yo93C}) constructed the moduli space of Gieseker (semi)stable Higgs sheaves. 

\begin{definition}
We denote by ${\mathrm{Higgs}}_{r,c_{1},c_{2}}^\textrm{tf,ss/s}$ (resp. $\mathrm{Higgs}_{r,c_{1},c_{2}}^\textrm{bun,ss/s}$) the coarse moduli space of Gieseker (semi)stable torsion-free $\mathcal L$-valued Higgs sheaves (resp. Higgs bundles)  with $\operatorname{rank} (\mathcal{E})=r, c_{1}(\mathcal{E})=c_{1}\in \mathrm {Im}(\text{NS}(X)\to H^2(X,\mathbb Q)), c_2(\mathcal{E})=c_{2}\in H^4(X,\mathbb Z)\cong \mathbb Z$. 
We call the affine space $\bm A:=\oplus_{i=1}^{r}H^{0}(X,\mathcal L^{\otimes i})$ the Hitchin base. The following characteristic polynomial map
\[
h^\text{tf} :{\operatorname{Higgs}}_{r}^\text{tf,ss} \rightarrow\bm A,\ (\mathcal E,\theta)\mapsto \text{char.poly.}(\theta)
\]
is called the Hitchin map which is projective by Yokogawa \cite[Corollary 5.12]{Yo93C}. 

\end{definition}
For later use, let
\[
p:R^{ss}\longrightarrow
{\operatorname{Higgs}}^{\mathrm{tf,ss}}_{r,c_1,c_2}
\]
be Yokogawa's framed parameter scheme and its good quotient.
\begin{definition}
Consider the natural map
$${\mathrm{Higgs}}_{r,c_{1},c_{2}}^\textrm{tf,ss} \to \mathrm {Pic}^{c_1}_X\times \Gamma(X,\mathcal L),\quad (\mathcal E,\theta)\mapsto (\det \mathcal E,\mathrm {tr}(\theta)).$$ Let $\mathcal{N}\in \text{Pic}(X)$ be a line bundle over $X$, then we denote the inverse image over $(\mathcal N,0)$ by  ${\mathrm{Higgs}}_{r,\mathcal{N},c_{2}}^\text{tf,ss}$ (resp. $\mathrm{Higgs}_{r,\mathcal{N},c_{2}}^\text{bun,ss} $). Hence  ${\mathrm{Higgs}}_{r,\mathcal{N},c_{2}}^\text{tf,ss}$ (resp. $\mathrm{Higgs}_{r,\mathcal{N},c_{2}}^\text{bun,ss} $) is the moduli space of semistable $\mathcal L$-valued Higgs sheaves  (resp. Higgs bundles) $(\mathcal{E},\theta)$ where $\det\mathcal{E}=\mathcal{N}$ and $\theta$ is traceless. We denote the corresponding Hitchin maps by $h_{\mathcal{N}}^\text{tf}$(resp. $h_{\mathcal{N}}$ ). 
\end{definition}
\begin{remark}
As pointed out by Tanaka--Thomas\cite{TT20}, ${\mathrm{Higgs}}_{r,\mathcal{N},c_{2}}^\text{tf,ss}$ is a better candidate to define the Vafa--Witten invariants.  
\end{remark}

When Gieseker semistability coincides with stability,  the $\mathbb G_m$-action on the moduli of Higgs sheaves is also used by Tanaka--Thomas\cite{TT20} to localize the virtual fundamental class to the fixed point components and define the Vafa--Witten invariant as 
$$\mathsf{VW}_{r,c_1,c_2}=\int_{[(\mathrm{Higgs}_{r,\mathcal N,c_{2}}^\textrm{tf,ss})^{\mathbb G_m}]^\textrm {vir}}\frac{1}{e(N^\textrm{vir})},$$
where $N^\textrm{vir}$ is the virtual normal bundle, see \cite[(4.2)]{TT20}. The definition below follows from \cite[\S 1.6]{TT20}.
\begin{definition}
Let $(X,\mathcal L)$ be a smooth polarized surface, $\mathcal{N}$ be a line bundle over $X$ and $\mathrm{Higgs}_{r,c_1,c_{2}}^\text{tf,ss}$, $\mathrm{Higgs}_{r,\mathcal{N},c_{2}}^\text{tf,ss} $ (resp. $\mathrm{Higgs}_{r,c_1,c_{2}}^\text{bun,ss} $, $\mathrm{Higgs}_{r,\mathcal{N},c_{2}}^\text{bun,ss} $) be the moduli space of Gieseker semistable  $\mathcal L$-valued Higgs sheaves  (resp. Higgs bundles) $(\mathcal{E},\theta)$. If we fix $\mathcal N$, we mean $\det\mathcal{E}$ is isomorphic to $\mathcal{N}$ and $\theta$ is traceless. Consider the components of the $\mathbb G_m$ fixed points ${(\mathrm{Higgs}_{r,c_1,c_{2}}^\text{tf,ss} )}^{\mathbb G_m}$ and ${(\mathrm{Higgs}_{r,\mathcal{N},c_{2}}^\text{tf,ss} )}^{\mathbb G_m}$,
\begin{enumerate}
	\item we call the components consists of $\mathbb G_m$ fixed points with zero Higgs fields as the instanton branch and  
	\item we call the components consists of $\mathbb G_m$ fixed points with non-zero Higgs fields as the monopole branch. 
\end{enumerate}
We call the fiber $h^{-1}(0)$ and $h_{\mathcal{N}}^{-1}(0)$ the global nilpotent cone.
\end{definition}
\begin{remark} We now make some brief remarks on the moduli spaces:
\begin{enumerate}
	\item In our case, we assume that our Higgs field lies in $\Gamma(X,\mathcal End(\mathcal E)\otimes \mathcal L)$, where $\mathcal L$ coincide with our polarization which is an ample line bundle.
	\item Tanaka--Thomas \cite{TT20} uses $\omega_X$ instead of a general very ample $\mathcal{L}$ here. Thus if we assume that our surface has very ample canonical bundle which also gives the polarization, then we get the moduli as in \cite{TT20}.
	\item If we give the natural $\mathbb{G}_m$-action on the Hitchin base, then the Hitchin map is $\mathbb G_m$-equivariant.	
	Then positiveness of the weights shows that the $\mathbb G_m$-fixed points $(\mathrm{Higgs}_{r,c_1,c_{2}}^\textrm{tf,ss})^{\mathbb G_m}$ and  $(\mathrm{Higgs}_{r,\mathcal N,c_{2}}^\textrm{tf,ss})^{\mathbb G_m}$ are contained in the global nilpotent cone. 
	\item If the Higgs field is zero, then the stability condition for the Higgs sheaves coincide with the stability condition of underlying torsion-free sheaves. Thus the instanton branch $F_\textrm{ins}$  is isomorphic to the corresponding moduli of semistable torsion-free sheaves. 
\end{enumerate}
\end{remark}\vspace{10pt}

\subsection{The generic fibers of Hitchin maps} 	As the curve case, we detect the generic fiber of the Hitchin map
$$h:{\mathrm{Higgs}}_{r,c_{1},c_{2}}^\textrm{tf,ss}\to \bm A,\quad (\mathcal E,\theta)\mapsto \textrm {char.poly.}(\theta)$$
by ``spectral surfaces" (see Definition \ref{def:spectral var}) and we can identify generic Hitchin fibers with Picard varieties of corresponding smooth spectral varieties.

Given a Higgs sheaf $(\mathcal{E},\theta)$, with Hitchin image $s\in \bm A$. We see that the Higgs field $\theta\in \Gamma(X,\mathcal End(\mathcal E)\otimes \mathcal L)$ is equivalent as an $\pi_{s*}\mathcal O_{X_s}$-module structure $\pi_{s*}\mathcal O_{X_s}\to \mathcal End(\mathcal E)$, which can be treated as a coherent sheaf on the spectral variety $X_{s}$.  Let $\pi_{s}:X_{s}\rightarrow X$ be the projection, which is finite flat of degree $r=\mathrm{rank} (\mathcal{E})$. Since $X_s/X$ is an affine morhphism, we have the following equivalence of categories (see \cite[Proposition 12.5]{GW20}) given by 
$$\pi_{s*}:\textrm {q-Coh}(X_s)\leftrightarrows \textrm {q-Coh}(\pi_{s*}\mathcal O_{X_s}\textrm {-mod}):\mathrm {tilde}_{X_s/X}.$$
Since $\pi_s$ is also finite, we see that the previous equivalent of categories restrict to the subcategories of coherent sheaves. In particular, we have

\begin{proposition}\label{prop:bundle over spectral}
Let $\mathcal{M}$ be a coherent sheaf over $X_{s}$. If $X_{s}$ is smooth and integral, then ${\pi_{s}}_*\mathcal{M}$ is a locally free sheaf of rank $r$ if and only if $\mathcal{M}$ is an invertible sheaf on $X_s$. 
\end{proposition}
By Bertini's theorem, see Section~\ref{subsec:higher dim}, generic $X_s$ is smooth and integral.
\begin{proof}
Over a Noetherian regular local ring $R$, a finitely generated module $M$ satisfies the following Auslander--Buchsbaum formula:
\[
\text{Projdim}_R M+\text{Depth}_R M=\dim R
\]	
Thus by calculating the depth, we know that ${\pi_{s}}_*\mathcal{M}$ is locally free if and only if the projective dimension of $\mathcal{M}$ is $0$. As a result, it is locally free of rank one.    	
\end{proof}
\begin{remark}
As we can see the above lemma also holds in higher dimensional cases.
\end{remark}

Let $X_s$ be a generic spectral surface. Hence it is smooth and connected.
\begin{lemma}
The canonical line bundle $\omega_{X_{s}}$ of $X_s$ is isomorphic to $\pi_{s}^*(\omega_{X}\otimes\mathcal{L}^{r-1})$.
\end{lemma}
\begin{proof}
By Lemma \ref{lem:can bund formula}, we have:
\[
\omega_{Y}\cong \pi^{*}(\omega_{X}\otimes\mathcal{L}^{\vee})\otimes\mathcal{O}_{Y/X}(-2).
\]
Then by the adjunction formula, we have:
\[
\omega_{X_s}\cong (\omega_{Y}\otimes \mathcal{O}_{Y}(X_s))|_{X_s}\cong \pi_{s}^*(\omega_{X}\otimes\mathcal{L}^{r-1}).
\]
The last isomorphism ``$\cong$'' is due to that $X_s$ lies in the linear system $\left\vert\pi^{*}\mathcal{L}^{r}\otimes \mathcal{O}_{Y/X}(r)\right\vert$ and $\mathcal{O}_{Y/X}(1)|_{X_s}$ is trivial.
\end{proof}

\begin{lemma} The Chern character of the cotangent bundle on the spectral surface satisfies the formula 
$$\textrm{ch}(\Omega_{X_s}^{1})=\pi_{s}^*(\textrm{ch}(\Omega_X^1)+\textrm{ch}(\mathcal{L}^{\vee})-\textrm{ch}(\mathcal{L}^{\vee})^{r}).$$
\end{lemma}
\begin{proof}
Recall the following exact sequences:
\[0\to \Omega_{Y/X}^1\to \pi^*(\mathcal{L}^\vee\oplus\mathcal{O}_X)(-1)\to \mathcal{O}_{Y}\to 0, \]
\[0\to \pi^*\Omega_X^1\to \Omega_{Y}^1\to \Omega_{Y/X}^1\to 0,\]
\[0\to \pi^*\mathcal{L}^{\otimes (-r)}(-r)|_{X_s}\to \Omega_{Y}^1|_{X_s}\to \Omega_{X_s}^1\to 0.\]
Notice that $\mathcal{O}(1)$ defines an infinite divisor, its restriction to spectral varieties is trivial. We have:
\[\begin{split}\textrm{ch}(\Omega_{X_s}^1)=&i_s^*\textrm{ch}(\Omega_{Y}^1)-\textrm{ch}(\pi^*_s\mathcal{L}^{\otimes (-r)})\\
	=&\pi^*_s\textrm{ch}(\Omega_X^1)+i_s^*\textrm{ch}(\Omega_{Y/X}^1)-\pi_{s}^*\textrm{ch}(\mathcal{L}^\vee)^{r}\\
	=&\pi^*_s\textrm{ch}(\Omega_X^1)+\pi^*_s(\textrm{ch}(\mathcal{L}^\vee)-\textrm{ch}(\mathcal{L}^\vee)^{r})
\end{split}\]

\end{proof}
\begin{lemma} Let $\text{Td}(X_s)$ be the Todd class of $X_s$. Then
$$\begin{array}{rcl}
	\mathrm{Td}(X_s)&=&1-\frac{1}{2}\pi_{s}^*(c_1(\omega_{X})+(r-1)c_1(\mathcal{L}))\\
	&&+\frac{1}{12}\pi_{s}^{*}(c_1(\omega_{X})^2+(2r-1)(r-1)c_1(\mathcal{L})^2\\&&+3(r-1)c_1(\omega_{X}). c_1(\mathcal{L})+c_2(TX)).
\end{array}$$
\end{lemma}
\begin{proof}
We first calculate $c_2(TX_s)$($=c_2(\Omega_{X_s}^1)$).
\begin{align*}
	c_{2}(TX_s)&=2+c_1(\Omega_{X_s}^1)+\frac{1}{2}c_1(\Omega_{X_s}^1)^{2}-\textrm{ch}(\Omega_{X_s}^1)\\
	&=\pi_{s}^*(r(r-1)c_1(\mathcal{L})^{2}+(r-1)c_1(\omega_{X}). c_1(\mathcal{L})+c_2(TX))
\end{align*}
Then we have:
\begin{align*}
	\mathrm{Td}(X_s)&=1-\frac{1}{2}c_1(\omega_{X_{s}})+\frac{1}{12}(c_1(\omega_{X_{s}})^{2}+c_2(TX_s))\\
	&=1-\frac{1}{2}\pi_{s}^*(c_1(\omega_{X})+(r-1)c_1(\mathcal{L}))\\&+\frac{1}{12}\pi_{s}^{*}(c_1(\omega_{X})^{2}+(2r-1)(r-1)c_1(\mathcal{L})^{2}+3(r-1)c_1(\omega_{X}). c_1(\mathcal{L})+c_2(TX)).
\end{align*}
\end{proof}

        Now let $(\mathcal{E},\theta)$ be a Higgs sheaf over $X$ with image $s$ in $\bm A$ under the Hitchin map $h^{\text{tf}}$. Then $(\mathcal{E},\theta)$ can be identified with torsion-free sheaf of rank $1$ on $X_{s}$  denoted by $\mathcal{M}\otimes \mathcal{I}_{\mathfrak D}$. Here $\mathcal{I}_{\mathfrak D}$ is the ideal sheaf of a $0$-dimensional closed subscheme $\mathfrak D$ of finite length $\mathrm{length}{\mathfrak D}$ on $X_{s}$. First recall the following result on Chern classes of torsion-free sheaves of rank 1 on a surface.
        \begin{lemma}
        	$\textrm{ch}(\mathcal{M}\otimes \mathcal{I}_{\mathfrak D})=\textrm{ch}(\mathcal{M})-[\mathfrak D]$.
        \end{lemma}
        \begin{proof}
        	We denote by $i:\mathfrak D\rightarrow X$ the closed immersion. Again by Grothendieck--Riemann--Roch, we have $c_{1}(i_{*}\mathcal{O}_{\mathfrak D})=0,c_{2}(i_{*}\mathcal{O}_{\mathfrak D})=-[\mathfrak D]$, then $c_2(\mathcal I_\mathfrak D)=[\mathfrak D]$ and $\textrm{ch}(\mathcal{I}_{\mathfrak D})=1+c_1(\mathcal{I}_{\mathfrak D})+\frac{(c_1^2-2c_2)}{2}=1-[\mathfrak D]$. As a result, $\textrm{ch}(\mathcal{M}\otimes \mathcal{I}_{\mathfrak D})=\textrm{ch}(\mathcal{M})-[\mathfrak D]$.      
        \end{proof}
        Let us fix $(r,c_{1},c_{2})$ as before. Now we can prove the first theorem on the nonemptyness of generic fibers of the Hitchin map.
        \begin{theorem}\label{thm:criterion of existence} Assume that $\pi_s^*:\mathrm {Pic}(X)\to \mathrm {Pic}(X_s)$ is an isomorphism for a very general $s\in\bm A$, the generic fibers of Hitchin map
        	\[
        	h:{\mathrm{Higgs}}_{r,c_{1},c_{2}}^\textrm{tf,ss}\rightarrow \bm{A}
        	\]
        	are nonempty if and only if the following equations 
        	\begin{equation}\label{eq:generic fib eq1}
        		r\delta=\textrm{ch}_1(\mathcal{E})+\frac{r(r-1)}{2}c_1(\mathcal{L})\end{equation}
        	\begin{equation}\label{eq:generic fib eq2}
        		(r-1)c_1(\mathcal{E})^{2}-2rc_2(\mathcal{E})=\frac{r^2(r^2-1)}{12}c_1({\mathcal{L}})^2-2r\pi_{s*}[\mathfrak D]
        	\end{equation}
        	have solutions for $\delta\in Im({\operatorname{Pic}}(X)\rightarrow H^2(X,\mathbb Q))$  and effective $[\mathfrak D]\in H^4(X_{s})$.
        \end{theorem}
   
        \begin{proof}
        	Since $\pi_{s*}(\mathcal{M}\otimes\mathcal{I}_{\mathfrak D})=\mathcal{E}$. By Grothendieck--Riemann--Roch, we have:
        	\begin{align*}
        		{\pi_{s}}_*(\textrm{ch}(\mathcal{M}\otimes\mathcal{I}_{\mathfrak D}).\mathrm{Td}(X_{s}))=\textrm{ch}(\mathcal{E}).\mathrm{Td}(X),\\
        		{\pi_{s}}_*(\mathrm{Td}(X_{s}))=\textrm{ch}({\pi_{s}}_*\mathcal{O}_{X_{s}}).\mathrm{Td}(X).
        	\end{align*}
        	
        	From the above two equations, we have:
        	\begin{equation}\label{eq: GRR}
        		{\pi_{s}}_*((c_{1}(\mathcal{M})+\frac{c_{1}(\mathcal{M})^{2}}{2}-[\mathfrak D]).\mathrm{Td}(X_{s}))=(\textrm{ch}(\mathcal{E})-\textrm{ch}({\pi_{s}}_*\mathcal{O}_{X_{s}})).\mathrm{Td}(X).
        	\end{equation}
        	By our previous computation of Todd classes of $X$ and $X_s$, the left hand side of (\ref{eq: GRR}) becomes
        	\begin{align*}&\pi_{s*}((c_{1}(\mathcal{M})+\frac{1}{2}c_{1}(\mathcal{M})^{2}-[\mathfrak D]).\mathrm{Td}(X_{s}))\\=& \pi_{s*}c_{1}(\mathcal{M})+\frac{1}{2}\pi_{s*}c_{1}(\mathcal{M})^{2}-\pi_{s*}[\mathfrak D]-\frac{1}{2}{\pi_{s}}_*c_{1}(\mathcal{M}). (c_1(\omega_{X})+(r-1)c_1(\mathcal{L})) ,
        	\end{align*}
        	and the right hand side of (\ref{eq: GRR}) becomes
        	{\footnotesize\begin{align*}
        			(&\textrm{ch}(\mathcal{E})-\textrm{ch}({\pi_{s}}_*\mathcal{O}_{X_{s}})).\mathrm{Td}(X)\\=&(\textrm{ch}(\mathcal{E})-\textrm{ch}({\pi_{s}}_*\mathcal{O}_{X_{s}})). (1-\frac{1}{2}c_1(\omega_{X}))\\
        			=&(\textrm{ch}_1(\mathcal{E})+\frac{r(r-1)}{2}c_1(\mathcal{L}))+\textrm{ch}_2(\mathcal{E})-\frac{r(r-1)(2r-1)}{12}c_1(\mathcal{L})^{2}-\frac{1}{2}(\textrm{ch}_1(\mathcal{E})+\frac{r(r-1)}{2}c_1(\mathcal{L})). c_1(\omega_{X}).
        	\end{align*}}
        	
        	We may assume that $c_1(\mathcal{M})=\pi_{s}^*\delta$, where $\delta$ is a first chern class of certain line bundle on $X$. Now we get the following:
        	\begin{align*}
        		&r\delta=\textrm{ch}_1(\mathcal{E})+\frac{r(r-1)}{2}c_1(\mathcal{L})\\
        		&r\delta^{2}-r\delta. (c_1(\omega_{X})+(r-1)c_1(\mathcal{L}))-2\pi_{s*}[\mathfrak D]=\\
        		&2\textrm{ch}_2(\mathcal{E})-\frac{r(r-1)(2r-1)}{6}c_1(\mathcal{L})^{2}-(\textrm{ch}_1(\mathcal{E})+\frac{r(r-1)}{2}c_1(\mathcal{L})). c_1(\omega_{X})
        	\end{align*}
        	We can reduce to the following equations:
        	
        	\begin{align*}
        		&r\delta=c_1+\frac{r(r-1)}{2}c_1(\mathcal{L})\\
        		&(r-1)c_1^{2}-2rc_2=\frac{r^2(r^2-1)}{12}c_1({\mathcal{L}})^2-2r\pi_{s*}[\mathfrak D],
        	\end{align*}
        	for  $\delta\in {\operatorname{NS}}(X)$, $[\mathfrak D]\in \mathbb{Z}_{\geq 0}$.
        
            Conversely, consider the direct image $\pi_{s*}({\mathcal{I}}_D\otimes \pi_s^*{\mathcal{N}})$ where $c_1({\mathcal{N}})=\delta$, its Chern numbers are as expected. Since the spectral surface is integral,  $\pi_{s*}({\mathcal{I}}_D)\otimes {\mathcal{N}}$ is stable as a Higgs sheaf.
        \end{proof}

        Motivated by Equation~\eqref{eq:generic fib eq2}, we introduce
        the following number.
        
        \begin{definition}\label{def:minimal c2}
        We define
        \[
        c_2^{\mathrm{g.bun}}
        :=
        \frac{r-1}{2r}c_1^2
        -
        \frac{r(r^2-1)}{24}c_1(\mathcal L)^2.
        \]
        It is a rational number depending only on $r$, $c_1$, and
        $\mathcal L$.
        \end{definition}
        
        \begin{remark}\label{rmk:meaning-gbun}
        With the notation of Theorem~\ref{thm:criterion of existence},
        Equation~\eqref{eq:generic fib eq2} can be rewritten as
        \[
        c_2
        =
        c_2^{\mathrm{g.bun}}
        +
        \operatorname{length}(\mathfrak D).
        \]
        Thus $c_2^{\mathrm{g.bun}}$ is the value of the second Chern
        class corresponding to the spectral construction with
        $\mathfrak D=\varnothing$, namely to a line bundle on the
        spectral surface and hence to a Higgs bundle on $X$.
        This explains the notation ``g.bun''.
        
        In general, the number $c_{2}^{\mathrm{g.bun}}$ defined above is a
        rational number. Under the divisibility assumption
        \[
        r\delta
        =
        c_{1}+\frac{r(r-1)}{2}c_{1}(\mathcal L)
        \]
        for some $\delta\in \operatorname{Im}
        \bigl(c_{1}\colon {\operatorname{Pic}}(X)\to H^{2}(X,\mathbb Q)\bigr)$, it is
        an integer. In fact, choose
        $\mathcal N\in\operatorname{Pic}(X)$ with
        $c_1(\mathcal N)=\delta$. For any smooth spectral surface $X_s$,
        the length-zero construction
        $\pi_{s*}\pi_s^*\mathcal N$
        is a vector bundle with first Chern class $c_1$. The preceding
        formula with $D=\varnothing$ gives
        $c_2(\pi_{s*}\pi_s^*\mathcal N)=c_2^{\mathrm{g.bun}}$.
        Hence $c_2^{\mathrm{g.bun}}\in\mathbb Z$.
    \end{remark}

    \begin{remark} 
    If the generic fibers are non-empty, the Hitchin map is surjective by the projectiveness proved in \cite[Theorem 5.10 and Corollary 5.12]{Yo93C}. 
    \end{remark}

\subsection{On Generic Fibration}
Though the above theorem might only hold for very general $s\in \bm{A}$ that $\pi_s^*:{\operatorname{Pic}}(X)\rightarrow {\operatorname{Pic}}(X_s)$ is an isomorphism, we can still describe generic fibers. To do this, we first construct a relative version of Beauville--Narasimhan--Ramanan Correspondence \cite{BNR}, i.e.,  a family of Higgs sheaves can be viewed as a family of torsion free sheaves of rank $1$ on the relative spectral surface assuming smoothness and integrality.

We denote the universal family of spectral surfaces over the Hitchin base:
\[
q_{\bm A}:\mathcal Y_{\bm A}\longrightarrow\bm A
\]
which is obtained by restricting the
universal divisor of Section~2 to the Hitchin base. In what follows, we will use the notation $X_s$ to denote a spectral surface for a geometric point $s\in \bm A$ and we will use the notation $\mathcal Y_T$ when we consider the family of spectral surface over an $\bm A$-scheme $T$. \footnote{This is to follow the convention in Section 2, where we have to distinguish two differential kinds of formal neighborhoods of spectral surfaces, hence we use ${\mathcal{Y}}$ to denote the universal family.}

Now for any $\bm A$ scheme $T/\bm A$, we have a finite morphism:
\[
\pi:\mathcal Y_{T}\rightarrow X\times T.
\]
Then $\pi^{*}\pi_X^*{\mathcal{L}}$ is a $q$-ample line bundle on ${\mathcal{Y}}_{T}$, where $\pi_X:X\times T\rightarrow X$ is the projection to $X$. We remark here that we have isomorphism of $\mathcal O_{X\times T}$-modules
$$ \pi_*\mathcal O_{\mathcal Y_T}\cong \mathcal O_{X\times T}\oplus\pi_X^*\mathcal L^{-1}\oplus \cdots \oplus \pi_X^*\mathcal L^{\otimes(-r+1)}.
$$
Here $\pi_*\mathcal O_{\mathcal Y_T}$ is the quotient affine $\mathcal O_{X\times T}$-algebra of $\mathrm {Sym}^\textrm{\tiny{\textbullet}}_{X\times T}(\pi_X^{*}\mathcal L^{-1})$ defined by the universal closed inclusion $\mathcal Y_T\subset {\operatorname {Tot}}(\pi_X^*({\mathcal{L}}))$. 

By the construction of the Higgs moduli space in \cite{Yo91}, a more natural topological invariant is the Hilbert polynomial. We recall the corresponding moduli functor as follows. 

\begin{definition}
    Let $P(x)$ be a numerical polynomial. For any $k$-scheme $T$, we define the Higgs moduli functor as
$$\mathrm{Higgs}^{\mathrm {tf,ss},\natural}_{r,P(x)} (T)=\left\{ ({\mathcal{E}},\Theta)\ \left| \ \begin{array}{c}
	{\mathcal{E}}\; \text{is a coherent sheaf on}\ X\times T \ \textrm{and } \\
	\text{satisfies the following conditions}\ (1) \ \text{and}\ (2),\\ \Theta:\mathcal E\to \mathcal E\otimes \pi_X^*\mathcal L\ \textrm{is a}\ \mathcal O_{X\times T}\textrm{-linear map}. 
\end{array} 
\right.  \right\} /\sim
$$
\begin{itemize}
	\item [(1)] $\mathcal E$ is flat over $T$ and is torsion-free of rank $r$ over $X\times\{t\}$ for any geometric point $t\in T$.
	\item[(2)]  the Hilbert polynomial of ${\mathcal{E}}\otimes_{{\mathcal{O}}_T} k(t)$ with respect to $\pi_X^*{\mathcal{L}}|_{X\times\{t\}}$ is $P(x)$ and $(\mathcal E,\Theta )$ is Higgs semi-stable with respect to $\pi_X^*{\mathcal{L}}|_{X\times\{t\}}$.
\end{itemize}
Here two families of Higgs sheaves $(\mathcal E,\Theta)\sim(\mathcal E_1,\Theta_1)$ if and only if they become isomorphic as Higgs sheaves after tensoring one of them with a line bundle pulled back from $T$. And for any morphism $T_1 \xrightarrow{ f}T_2$ of $k$-schemes, we define the corresponding map $\mathrm{Higgs}^{\mathrm {tf,ss},\natural}_{r,P(x)} (f)$ to be the pull back of Higgs sheaves by the base change map $X\times{T_1}\to X\times T_2$.  Here the sign ``$\natural$" emphasize that this is a presheaf on $\mathrm {Sch}/k$. 
\end{definition}
Yokogawa\cite{Yo91} proved that the moduli functor is corepresented by a scheme $\mathrm{Higgs}^{\mathrm {tf,ss}}_{r,P(x)}$. In this case, we also have the Hitchin map
$$h:{\mathrm{Higgs}}_{r,P(x)}^\textrm{tf,ss}\to \bm A,\quad (\mathcal E,\theta)\mapsto \textrm {char.poly.}(\theta)$$
and the properness of $h$ is also proved by Yokogawa in \cite{Yo93C}.

Let $U$ be the open subset of $\bm A$ parameterize smooth spectral surfaces, we consider the functor 
$$U\times_{\bm A}{\operatorname{Higgs}}^{\textrm{tf,ss},\natural}_{r,P(x)}.$$
For any $U$-scheme $T$, $U\times_{\bm A}{\operatorname{Higgs}}^{\textrm{tf,ss},\natural}_{r,P(x)}(T)$ consists of flat families of Higgs sheaves on $X\times T$ such that the image of $(\mathcal E,\Theta)$ under the Hitchin map lies in $U$.

\begin{definition}
    For any $T\in \text{Sch}/U$, we put:
\[
\mathrm{M}^{s,\natural}_{q_U,1,P(x)} (T)=\{{\mathcal{F}}|{\mathcal{F}}\; \text{is a coherent sheaf on}\ {\mathcal{Y}}_U\times_{U}T \ \text{and satisfies}\ \eqref{eq:flatness}\;\ \text{and}\;\eqref{eq:stability in family}\}/\sim
\]
\begin{enumerate}
	\item\label{eq:flatness} ${\mathcal{F}}$ is flat over $T$ and for all geometric point $t$ of $T$,  ${\mathcal{F}}\otimes_{{\mathcal{O}}_T} k(t)$ on $\mathcal Y_t$ is a torsion free sheaf of rank $1$.
	\item \label{eq:stability in family}  The Hilbert polynomial of ${\mathcal{F}}\otimes_{{\mathcal{O}}_T} k(t)$ with respect to $(\pi^{*}\pi_X^*{\mathcal{L}})|_{{\mathcal{Y}}_t}$ is $P(x)$ and is stable with respect to $(\pi^{*}\pi_X^*{\mathcal{L}})|_{{\mathcal{Y}}_t}$.
\end{enumerate}

Here ${\mathcal{F}}\sim{\mathcal{F}}'$ if and only if there is a line bundle ${\mathcal{M}}$ on $T$ such that ${\mathcal{F}}\cong{\mathcal{F}}'\otimes q_T^*{\mathcal{M}}$ where $q_T:{\mathcal{Y}}_U\times_UT\rightarrow T$. And for any morphism $f:T_1 \xrightarrow{ f}T_2$ of $U$-schemes, we define the corresponding map $\mathrm{M}^{s,\natural}_{q_U,1,P(x)}(f)$ to be the pull back of sheaves.
\end{definition}
We now prove that this spectral functor is naturally isomorphic to the restriction of the Higgs moduli functor over $U$.

\begin{lemma}\label{lem:relative BNR}
	We have an isomorphism of functors	$$\mathrm{BNR}_U:\mathrm{M}_{q_U,1,P(x)}^{s,\natural}\xrightarrow{\cong} U\times_{\bm A}{\operatorname{Higgs}}^{\textrm{tf,ss},\natural}_{r,P(x)},$$ 
	given by push-forward of the natural projection $\pi: {\mathcal{Y}}_U\rightarrow X\times U$. Hence, an isomorphism of the corresponding moduli spaces
    $$\mathrm{BNR}_U:\mathrm{M}_{q_U,1,P(x)}^{s}\xrightarrow{\cong} U\times_{\bm A}{\operatorname{Higgs}}^{\textrm{tf,ss}}_{r,P(x)}.$$
\end{lemma}    
\begin{proof}
Simpson \cite[Lemma 6.5]{Sim94II} (see also  \cite[\S 1]{Yo91}) shows that a Higgs sheaf sheaf on $X\times T$ over $T$ is the same thing as an $\mathcal O_{X\times T}$-coherent $\mathrm {Sym}^\textrm{\tiny{\textbullet}}_{X\times T}(\pi_X^{*}\mathcal L^{-1})$-module on $X\times T$.  Moreover, by the Cayley--Hamilton theorem, the ideal of the spectral surface in $\mathrm {Sym}^\textrm{\tiny{\textbullet}}_{X\times T}(\pi_X^{*}\mathcal L^{-1})$ annihilate ${\mathcal{E}}$ and thus $\mathcal E$ is indeed a $\pi_*\mathcal O_{\mathcal Y_T}$ module. 
Hence, pushforward along
$\pi:\mathcal Y_T\to X\times T$
identifies coherent sheaves on $\mathcal Y_T$ with coherent
Higgs sheaves having the prescribed characteristic polynomial.
This correspondence preserves flatness over $T$ and commutes
with base change.

On every geometric fiber, integrality of $\mathcal Y_t$ and
$\deg(\pi_t)=r$ give
\[
\operatorname{rk}(\mathcal E_t)
=
r\,\operatorname{rk}(\mathcal F_t),
\qquad
\mathcal E_t=\pi_{t*}\mathcal F_t,
\]
and torsion-freeness corresponds under pushforward.
Thus the rank-$r$ torsion-free Higgs sheaves under consideration
correspond precisely to rank-one torsion-free spectral sheaves.
Higgs-invariant subsheaves correspond to spectral subsheaves,
and the projection formula shows that Hilbert polynomials are
preserved. Since rank-one torsion-free sheaves on an integral
surface are Gieseker stable, the corresponding Higgs sheaves
are stable as well.

Finally, the correspondence respects tensoring by line bundles
pulled back from $T$. It therefore induces the asserted
isomorphism of functors and hence of their coarse moduli spaces.
\end{proof}

	We see that the restricted Hitchin map $h|_{U}:h^{-1}(U)\to U$ should be understood as the moduli, $\mathrm{M}^{s,\natural}_{q_U,r,P(x)}$, of families of torsion free sheaves of rank $1$. Altman and Kleiman\cite{AK79,AK80} show that the generic rank $1$ moduli functor (which is the \'etale sheafification of $\mathrm{M}^{s,\natural}_{q_U,1,P(x)}$) is represented by a scheme. Moreover, the projectivity of the moduli space is reproved by Koll\'ar via a different method in \cite[Section 6]{Kol90} and he \cite[Lemma 6.13]{Kol90}  also shows that 
	\begin{lemma} Let $f: X\to S$ be smooth and let $F$ be a sheaf on $X$, flat over $S$. Assume that for every $s \in  S$ the restriction $F|_{X_s}$ is torsion
		free of rank one. Then the double dual of $F$ is locally free. 
	\end{lemma}
	
	By this, we can show that the double dual coincides with the determinant and have the following interpretation of $\mathrm{M}^{s,\natural}_{q_U,1,P(x)}$.
	
	\begin{lemma}

For each $d\in \mathbb Z_{\geq 0}$, let
$\operatorname{Hilb}^{d,\natural}_{q_U}$ denote the relative Hilbert functor of zero-dimensional subschemes of length $d$, and let
$\operatorname{Pic}^{P(x)+d,\natural}_{q_U}$ denote the relative Picard functor of line bundles with Hilbert polynomial $P(x)+d$. Then there is a natural isomorphism of functors
\[
F^\natural_{\mathrm{MPH}}\colon
\mathrm{M}^{s,\natural}_{q_U,1,P(x)}
\xrightarrow{\ \sim\ }
\bigsqcup_{d\geq 0}
\left(
\operatorname{Pic}^{P(x)+d,\natural}_{q_U}
\times_U
\operatorname{Hilb}^{d,\natural}_{q_U}
\right),
\]
given by
\[
\mathcal F\mapsto
\left(
\det(\mathcal F),
\operatorname{Coker}
\bigl(
\mathcal F\otimes\det(\mathcal F)^\vee
\hookrightarrow \mathcal O_{\mathcal Y_T}
\bigr)
\right).
\]
Only the nonempty summands are understood to occur in the above disjoint union. 		Here $P(x)$ is the numerical polynomial defined by $P(n)=\chi(\mathcal F|_t\otimes \pi_t^*\mathcal L^{n})$ and ``MPH" means that the natrual isomorphism split the moduli functor $\mathrm{M}^{s,\natural}_{q,1,P(x)}$  into a product of Picard functor and Hilbert functor.

	\end{lemma}
	
	\begin{remark}\label{rmk:FMPH} Before the proof, we should give some remarks on the definition of the functor.
		\begin{itemize}
			\item [(1)] The determinant of a coherent sheaf. see \cite{KM76}, \cite[Theorem 4.11]{Mar78} and also \cite[15.122 tag 0FJI]{stacks-project}. Let $\mathcal F$ be a coherent sheaf with finite length locally free resolution on a scheme, one can define $\det (\mathcal F):=\bigotimes_{i}\ \det(\mathcal E_i)^{(-1)^i}$ for a finite length resolution $\mathcal E_\textrm{\tiny{\textbullet}}$ of $\mathcal F$ by locally free sheaves. The determinant is independent of the chosen complex $\mathcal E_\textrm{\tiny{\textbullet}}$ and compatible with base changes. In our case, by the smoothness of $\mathcal Y_T/T$ and \cite[Lemma 15.77.6, tag 09PC]{stacks-project}, our coherent sheaves admit finite length locally free resolutions. See also \cite[\S 8.1]{BLMNPS21} to define relatively perfect objects in a much more general setting.
			\item[(2)] The determinant and the double dual.
By Koll\'ar's \cite[6.13.\ Lemma]{Kol90},
$\mathcal F^{\vee\vee}$ is a line bundle on $\mathcal Y_T$.
We claim that there is a canonical isomorphism
$\det(\mathcal F)\cong \mathcal F^{\vee\vee}$
which is compatible with base change.

Consider the double-dual morphism $\mathrm{dd}_{\mathcal F}\colon
\mathcal F\rightarrow \mathcal F^{\vee\vee}$. 
By item~(1), $\mathcal F$ is perfect, while
$\mathcal F^{\vee\vee}$ is locally free by Koll\'ar's lemma.
Thus $\mathrm{dd}_{\mathcal F}$ may be regarded as a morphism
of perfect complexes. We use the determinant--divisor construction of
Knudsen--Mumford \cite[Chapter~II, pp.~47--49]{KM76}.
Recall that a good morphism
$\lambda\colon \mathcal E^\bullet\to\mathcal G^\bullet$
of perfect complexes determines a Cartier divisor
$\operatorname{Div}(\lambda)$ and a canonical isomorphism
$\det(\lambda)\colon
\det(\mathcal E^\bullet)
\bigl(\operatorname{Div}(\lambda)\bigr)
\xrightarrow{\sim}
\det(\mathcal G^\bullet)$. 
We check that $\mathrm{dd}_{\mathcal F}$ is good and that its
associated divisor vanishes. Let $y\in\mathcal Y_T$ be a point
with $\operatorname{depth}\mathcal O_{\mathcal Y_T,y}\leq 1$,
and let $t\in T$ be its image. Since
$q_T\colon\mathcal Y_T\to T$ is smooth, the depth formula for
flat local homomorphisms gives
\[
\operatorname{depth}\mathcal O_{\mathcal Y_T,y}
=
\operatorname{depth}\mathcal O_{T,t}
+
\operatorname{depth}\mathcal O_{\mathcal Y_t,y}.
\]
Hence $\operatorname{depth}\mathcal O_{\mathcal Y_t,y}\leq 1$.
The fiber $\mathcal Y_t$ is a smooth surface, so
$\mathcal O_{\mathcal Y_t,y}$ is a regular local ring of
dimension at most one. Since $\mathcal F_t$ is torsion free of
rank one, $(\mathcal F_t)_y$ is therefore free of rank one.
As $\mathcal F$ is flat over $T$, the fiberwise criterion for
freeness implies that $\mathcal F$ is locally free at $y$.
Consequently, $\mathrm{dd}_{\mathcal F}$ is an isomorphism at every point of depth at most one. It follows that $\mathrm{dd}_{\mathcal F}$ is a good morphism
in the sense of Knudsen--Mumford and that
$\operatorname{Supp}(\mathrm{dd}_{\mathcal F})$ contains no
points of depth one. Since 
$\operatorname{Supp}
\bigl(\operatorname{Div}(\mathrm{dd}_{\mathcal F})\bigr)
\subseteq
\operatorname{Supp}(\mathrm{dd}_{\mathcal F})$,
we obtain $\operatorname{Div}(\mathrm{dd}_{\mathcal F})=0$.
The Knudsen--Mumford construction therefore gives a canonical
isomorphism
\[
\alpha_{\mathcal F}\colon
\det(\mathcal F)
\xrightarrow{\cong}
\det(\mathcal F^{\vee\vee})
=
\mathcal F^{\vee\vee}.
\]
Moreover, this construction is compatible with base change.

Composing the double-dual map with
$\alpha_{\mathcal F}^{-1}$ gives a natural morphism
$\mathcal F
\rightarrow
\det(\mathcal F)$. On every geometric fiber it is the usual inclusion of a torsion-free rank-one sheaf into its reflexive hull, and hence is injective. Thus we obtain the natural inclusion
$\mathcal F\otimes\det(\mathcal F)^\vee
\hookrightarrow
\mathcal O_{\mathcal Y_T}$.

\item[(3)] We can restrict the inclusion ${\mathcal{F}}\otimes(\det{\mathcal{F}})^{\vee}\hookrightarrow {\mathcal{O}}_{\mathcal Y_T}$ over a geometric point $t\in T$ and obtain an exact sequence
\[
0\rightarrow
\mathcal F_t\otimes\det(\mathcal F_t)^\vee
\rightarrow
\mathcal O_{\mathcal Y_t}
\rightarrow
\mathcal O_{D_t}
\rightarrow 0.
\] Here $D_t\subset\mathcal Y_t$ is a zero-dimensional closed subscheme. 			 Let $d=\operatorname{length}(D_t)=
c_2(\mathcal F_t)$. Since tensoring a zero-dimensional sheaf with a line bundle does not change its Hilbert polynomial, the preceding exact sequence gives
\[
P_{\det(\mathcal F_t)}(x)-P_{\mathcal F_t}(x)
=
P_{\det(\mathcal F_t)\otimes\mathcal O_{D_t}}(x)
=
d.
\]
In particular, fixing $P(x)$ does not fix the Hilbert polynomial of the determinant. For a family over a connected base $T$, the Hilbert polynomial of the line bundle $\det(\mathcal F_t)$ is locally constant; hence $d$ is constant on $T$. This explains why the target of $F^\natural_{\mathrm{MPH}}$ must be decomposed according to $d$.

		\end{itemize}
	\end{remark}
	
	\begin{proof} 
        By Remark~\ref{rmk:FMPH}, a family $\mathcal F$ over a connected $U$-scheme $T$ determines an ideal $I_D=
        \mathcal F\otimes\det(\mathcal F)^\vee\hookrightarrow\mathcal O_{\mathcal Y_T}$ whose fibers define zero-dimensional subschemes of a constant length $d$. The quotient $\mathcal O_D$ is flat over $T$: indeed, both $\mathcal F\otimes\det(\mathcal F)^\vee$ and $\mathcal O_{\mathcal Y_T}$ are flat over $T$, and the inclusion remains injective on every geometric fiber. Hence
        $\mathcal F\mapsto \bigl(\det(\mathcal F),D\bigr)$ defines a natural transformation
        \[
        \mathrm M^{s}_{q_U,1,P(x)}
        \rightarrow
        \bigsqcup_{d\geq0}
        \left(
        \operatorname{Pic}^{P(x)+d,\natural}_{q_U}
        \times_U
        \operatorname{Hilb}^{d,\natural}_{q_U}
        \right).
        \]
        Its inverse is given by the tensor product. So the two constructions are mutually inverse and compatible with base change. 
    \end{proof}
	Indeed, the natural isomorphism $F_{MPH}^{\natural}$ passes to the \'etale shefification.
	\begin{corollary}\label{coro:FMPH} 
       We have an isomorphism of the corresponding moduli spaces
        \[
        F_{\mathrm{MPH}}\colon
        \mathrm{M}^{s}_{q_U,1,P(x)}
        \xrightarrow{\cong }
        \bigsqcup_{d\geq 0}
        \left(
        \operatorname{Pic}^{P(x)+d}_{q_U}
        \times_U
        \operatorname{Hilb}^{d}_{q_U}
        \right).
        \] which represent the \'etale sheafification of the previous functors.
        Since $\mathrm{M}^{s}_{q_U,1,P(x)}$ is of finite type, only finitely many of these summands are nonempty.	
	\end{corollary}
	
	\begin{remark} For rank $r\ge 2$, taking double dual might not induce a morphism from moduli spaces of torsion-free sheaves with fixed chern classes to moduli spaces of vector bundles (Cf. \cite{Artamkin91}, \cite{GL94} and \cite{O96}). In \cite{O96}, the locus which parametrize non-locally free sheaves in the moduli of torsion free sheaves is called boundary. O'Grady \cite[Theorem A]{O96} shows that, in most cases, the boundary is non-empty. Thus for $c_2>\!\!\!>0$ such that the moduli space of torsion free sheaves is irreducible,  we can get a flat family of torsion free sheaves with special fibers not locally free. However, taking double dual will decrease the $c_2$ for torsion free sheaves, and keep the $c_2$ for bundles.
	\end{remark}

	By the previous corollary, we also have:
	\begin{corollary}\label{cor:genericsmoothcoarser}
		There exists an open subset $V$ of $\bm A$ parameterizing smooth spectral surfaces, such that the restriction of the Hitchin map $h|_{V}:V\times_{\bm A}{\operatorname{Higgs}}^{\textrm{tf,ss},\natural}_{r,P(x)}\rightarrow V$	is smooth.
	\end{corollary}

        \begin{proof}
By Lemma~\ref{lem:relative BNR} and
Corollary~\ref{coro:FMPH}, every geometric fiber over $U$
is a finite disjoint union of products of Picard schemes
and Hilbert schemes of points on smooth surfaces.
These fibers are smooth, since we work over $\mathbb C$.
By generic flatness, after shrinking $U$ to a nonempty
open subset $V$, the morphism $h^{-1}(V)\to V$ is flat.
Being of finite presentation with smooth geometric fibers,
it is therefore smooth.
\end{proof}

	In what follows, we consider our moduli of Higgs bundles with refined invariants $(r,c_1,c_2)$. In fact, this moduli space with refined invariants are union of connected componentes in $\mathrm {Higgs}_{r,P(x)}^\mathrm {tf,ss}$. As pointed out in \cite{Sasha11}, we have the following lemma.
	\begin{lemma} The rational Chern classes $c_1,c_2$ are locally constant in a flat family of coherent sheaves on a smooth projective algebraic surface. That is, let $X$ be a smooth projective algebraic surface over $\mathbb C$, $\mathcal F$ is a flat family of torsion free coherent sheaves on $X\times T$, then $c_1(\mathcal F|_t)\in H^2(X,\mathbb Q)$ and $c_2(\mathcal F|_t)\in H^4(X,\mathbb Q)$ are locally constant with respect to $t\in T$.   
	\end{lemma}
	\begin{proof} Let $\mathcal L$ be an ample line bundle on $X$. Then we have the Hilbert polynomial of $\mathcal F|_t$ on $X$ given by 
		{\small
			\begin{align*}
				\chi(\mathcal F|_t\otimes \mathcal L^n)=\frac{rc_1(\mathcal L)^2}{2}\cdot n^2+r(\mu_\mathcal L(\mathcal F|_t)+\mu_{\mathcal L}(T_X))\cdot n+\chi (\mathcal F|_t).
		\end{align*}}
		
		\noindent The Hilbert polynomial is locally constant in a flat family of coherent sheaves, thus $\mu_{\mathcal L}(\mathcal F_t)$ is locally constant with respect to $t$. By this, if we take $\mathcal L=L_1,...,L_n$ such that $L_1,...,L_n$ are ample line bundles on $X$ and $\mathbb R$-linear generate the $\mathrm {NS}(X)_{\mathbb R}$ (Cf. \cite[Positivity I, Example 1.3.14]{Laz04} ), then $c_1(\mathcal F_t)\in H^2(X,\mathbb Q)$ is locally constant with respect to $t$. From this, since $\chi (\mathcal F|_t)$ is locally constant, thus $c_2(\mathcal F_t)\in H^4(X,\mathbb Q)$ is also constant with respect to $t$. \end{proof}

	By this, ${\mathrm{Higgs}}_{r,c_{1},c_{2}}^\textrm{tf,ss}$ is a union of connected componentes in $\mathrm {Higgs}_{r,P(x)}^\mathrm {tf,ss}$. Thus we have the following corollary.

    \begin{corollary}\label{cor:genericsmooth}
        Assume that $\pi_s^*:\operatorname{Pic}(X)\xrightarrow{\sim}
        \operatorname{Pic}(X_s)$
        for very general $s\in\bm A$, and assume that the generic
        fiber of      $h:\operatorname{Higgs}^{\mathrm{tf,ss}}_{r,c_1,c_2}
        \rightarrow \bm A$
        is nonempty. Then there exists a nonempty Zariski open subset $W\subset\bm A$ over which the spectral construction induces an isomorphism
\[
h^{-1}(W)
\cong
\operatorname{Pic}_X^\delta
\times\operatorname{Hilb}^{n}_{q_W},
\qquad
\delta=\frac{c_1}{r}+\frac{r-1}{2}c_1(\mathcal L),
\quad
n=c_2-c_2^{\mathrm{g.bun}}.
\]
Here $\operatorname{Pic}_X^\delta$ parametrizes line bundles on $X$ with rational first Chern class $\delta$. In particular, $h^{-1}(W)\to W$ is smooth and projective.
    \end{corollary}
        \begin{remark}
            Notice that here we consider the Picard varities of $X$. Compare this with Corollary \ref{coro:FMPH} where we consider the product of relative Picard schemes with relative Hilbert schemes of points.
        \end{remark}
        \begin{proof}
        By Corollary \ref{cor:genericsmoothcoarser}, we can choose a nonempty open subset $W\subset\bm A$ parametrizing smooth spectral surfaces such that $h^{-1}(W)\to W$ is smooth. 
        
        Since we have
$\operatorname{Nm}_{\pi_W}\circ\pi_W^*=[r]$, the kernel of the relative pullback
$\pi_W^*:\operatorname{Pic}_X\times W
\to\operatorname{Pic}_{q_W}$
lie inside the finite \'etale group scheme
$\operatorname{Pic}_X[r]\times W$.
Its geometric fibers are trivial by the hypothesis, so $\pi_W^*$ is a monomorphism. Since $\operatorname{Pic}_X^\delta$ is proper and $\operatorname{Pic}_{q_W}$ is separated over $W$, the restriction of $\pi_W^*$ to
$\operatorname{Pic}_X^\delta\times W$ is a closed immersion.

Taking its product over $W$ with
$\operatorname{Hilb}^{n}_{q_W}$, the identifications of Lemma~\ref{lem:relative BNR} and Corollary~\ref{coro:FMPH}, together with the Chern class formulas of Theorem~\ref{thm:criterion of existence}, give a closed immersion
\[
\Phi_W:
\operatorname{Pic}_X^\delta
\times\operatorname{Hilb}^{n}_{q_W}
\rightarrow h^{-1}(W),
\]
given fiberwise by
$(N,D)\mapsto\pi_{s*}(\pi_s^*N\otimes I_D)$.

For very general $s\in W$, the hypothesis on $\pi_s^*$ and the same Chern class formulas imply that
$(\Phi_W)_s$ is surjective. Since $h^{-1}(W)\to W$
is smooth, the subset
\[
h\bigl(h^{-1}(W)\setminus\operatorname{Im}(\Phi_W)\bigr)
\]
is open in $W$ and contains no very general point,
so it is empty. Thus $\Phi_W$ is a surjective closed
immersion into the reduced scheme $h^{-1}(W)$, and hence an isomorphism. Smoothness holds by the choice of $W$
\end{proof}

		\section{Instanton and Monopole Branches} 
		
		In this section, we use Hitchin maps to study the geometric properties of the Higgs bundle moduli spaces and its $\mathbb G_m$-fixed points. We start from the following lemma.

		\begin{lemma}\label{lem:discriminant equality}
			For a torsion-free sheaf $\mathcal{E}$ on $X$, suppose it is endowed with the Harder--Narasimhan filtration of length $m$:
			\[
			0=\mathrm {HN}_0\mathcal{E}\subset \mathrm {HN}_1\mathcal{E}\subset\ldots\subset \mathrm {HN}_{m}\mathcal{E}=\mathcal{E},
			\]
			we put $\mathcal{E}_i=\mathrm{gr}_{i}^{\mathrm {HN}}\mathcal{E}=\frac{\mathrm {HN}_{i}\mathcal{E}}{\mathrm {HN}_{i-1}\mathcal{E}}$ and $r_i$ as the generic rank of $\mathcal{E}_i$ for $1\le i\le m$. Then we have:
			\[
			\frac{\Delta(\mathcal{E})}{r}=\sum_{i=1}^{m}\frac{\Delta(\mathcal{E}_i)}{r_i}-\sum_{1\le i<j\le m}\frac{r_ir_j}{r}(\frac{c_1(\mathcal{E}_i)}{r_i}-\frac{c_1(\mathcal{E}_j)}{r_j})^2.
			\]
		\end{lemma}
		\begin{proof}
			Notice that:
			\begin{align*}
				\frac{\Delta(\mathcal{E})}{2r}&=c_2(\mathcal{E})-\frac{r-1}{2r}c_1(\mathcal{E})^{2}\\       c_2(\mathcal{E})&=\sum_{i}c_2(\mathcal{E}_i)+\sum_{i<j}c_1(\mathcal{E}_i)c_1(\mathcal{E}_j)
			\end{align*}
			Hence we have:
			\begin{align*}
				\frac{\Delta(\mathcal{E})}{2r}&= \sum_{i}c_2(\mathcal{E}_i)+\sum_{i<j}c_1(\mathcal{E}_i)c_1(\mathcal{E}_j)-\frac{r-1}{2r}(\sum_{i}c_1(\mathcal{E}_i))^{2}\\
				&=\sum_{i}(\frac{\Delta(\mathcal{E}_i)}{2r_i}+\frac{r_i-1}{2r_i}c_1(\mathcal{E}_i)^{2})+\sum_{i<j}c_1(\mathcal{E}_i)c_1(\mathcal{E}_j)-\frac{r-1}{2r}(\sum_{i}c_1(\mathcal{E}_i))^{2}\\
				&=\sum_{i}\frac{\Delta(\mathcal{E}_i)}{2r_i}-\frac{1}{2r}\left(\sum_{i}\frac{\sum_{j\ne i}r_j}{r_i}c_{1}(\mathcal{E}_i)^2-2\sum_{i<j}c_1(\mathcal{E}_i)c_1(\mathcal{E}_j)\right)\\
				&=\sum_{i}\frac{\Delta(\mathcal{E}_i)}{2r_i}-\frac{1}{2r}\sum_{i<j}r_ir_j(\frac{c_1(\mathcal{E}_i)}{r_i}-\frac{c_1(\mathcal{E}_j)}{r_j})^2.
			\end{align*}
			The result follows.
		\end{proof}
		\begin{remark}
			As can be seen from the proof, the lemma holds for general filtrations.
		\end{remark}

		\begin{lemma}\label{lem:inequality from Higgs}
			Let $(\mathcal E,\theta)$ be a Gieseker-semistable Higgs sheaf. We have the inequalities $0<\frac{c_1(\mathcal{E}_i).c_1(\mathcal{L})}{r_i}-\frac{c_1(\mathcal{E}_{i+1}).c_1(\mathcal{L})}{r_{i+1}}\le c_1(\mathcal{L})^{2}$. 
		\end{lemma}
		
		\begin{proof}
			By the Gieseker semistablity of $(\mathcal{E},\theta)$, the induced map by the Higgs field:
			\[
			\mathrm{HN}_i\mathcal{E}\rightarrow \frac{\mathcal{E}}{\mathrm{HN}_i\mathcal{E}}\otimes\mathcal{L}
			\]
			is nonzero. Then $\mu_{\min}\mathrm{HN}_i\le\mu_{\max}(\frac{\mathcal{E}}{\mathrm{HN}_i\mathcal{E}})+c_1(\mathcal{L})^{2}$. And the lemma follows.
		\end{proof}
		We recall the Hodge index theorem: for any divisors $D$ on $X$, 
		\[
		(D. c_1({{\mathcal{L}}}))^2\ge D^2\cdot c_1({\mathcal{L}})^{2}
		\]
		\begin{lemma}\label{lem:HN inequality}
			For the Harder--Narasimhan filtration, we have:
			\[
			\sum_{1\le i<j\le m}\frac{r_ir_j}{r}(\frac{c_1(\mathcal{E}_i)}{r_i}-\frac{c_1(\mathcal{E}_j)}{r_j})^2\le  \sum_{1\le i<j\le m}\frac{r_ir_j}{r}(i-j)^{2} c_1(\mathcal{L})^2.
			\]
		\end{lemma}
		\begin{proof}
			By Lemma \ref{lem:inequality from Higgs}, 
			\[
			0<\frac{c_1(\mathcal{E}_i).c_1({\mathcal{L}})}{r_i}-\frac{c_1(\mathcal{E}_{i+1}).c_1({\mathcal{L}})}{r_{i+1}}\le c_1(\mathcal{L})^{2}.
			\]
			Then combined with the Hodge index theorem:
			\begin{align*}
				\sum_{1\le i<j\le m}\frac{r_ir_j}{r}(\frac{c_1(\mathcal{E}_i)}{r_i}-\frac{c_1(\mathcal{E}_j)}{r_j})^2. c_1(\mathcal{L})^2&\le \sum_{1\le i<j\le m}\frac{r_ir_j}{r}(\frac{c_1(\mathcal{E}_i).c_1(\mathcal{L})}{r_i}-\frac{c_1(\mathcal{E}_{j}).c_1(\mathcal{L})}{r_{j}})^2\\
				&= \sum_{1\le i<j\le m}\frac{r_ir_j}{r} (\sum_{s=i}^{j-1}\frac{c_1(\mathcal{E}_s). c_1(\mathcal{L})}{r_s}-\frac{c_1(\mathcal{E}_{s+1}). c_1(\mathcal{L})}{r_{s+1}})^2\\
				&\le \sum_{1\le i<j\le m}\frac{r_ir_j}{r}(i-j)^{2} (c_1(\mathcal{L})^2)^{2}
			\end{align*}
			
			Hence we get the conclusion.
		\end{proof}
		\begin{definition}
			We denote by $^{GF}{\operatorname{Higgs}}^{\text{tf},ss}_{r,c_1,c_2}$ the union of irreducible components which dominate the Hitchin base, i.e., the closure of generic fibers.
		\end{definition}

        	Recall the number $c_2^{\mathrm{g.bun}}$ introduced in Definition~\ref{def:minimal c2} (see also \ref{rmk:meaning-gbun}), 
			$$c_2^{\text{g.bun}}=\frac{(r-1)}{2r}c_1^2-\frac{r(r^2-1)}{24}c_1(\mathcal L)^2.$$ 
            This is a rational number depending on $r$, $c_{1}$, and $\mathcal L$. Under the divisibility assumption appearing below, it is an integer; see the following remark. Here g.bun means that the Hitchin map of Higgs bundles has nonempty generic fibers.

		\begin{remark} For torsion free sheaf $\mathcal E$ on $X$, we have equality
			\begin{equation}\label{eq:c2c2gbun}c_2(\mathcal E)-c_2^{\text{g.bun}}=\frac{\Delta(\mathcal E)}{2r}+\frac{r(r^2-1)}{24}c_1(\mathcal L)^2.
			\end{equation}
		\end{remark}

		\begin{theorem}\label{thm:graded piece are simple}
            We Assume that $\pi_s^*:{\operatorname{Pic}}(X)\xrightarrow{\sim}{\operatorname{Pic}}(X_s)$ for a very general $s\in\bm A$, and that the equation
            \[
            r\delta
            =
            c_1+\frac{r(r-1)}{2}c_1(\mathcal L)
            \]
            has a solution $\delta\in
            \operatorname{Im}
            \bigl(c_1:{\operatorname{Pic}}(X)\to \mathrm H^2(X,\mathbb Q)\bigr)$. Then the following statements hold.
            
            \begin{enumerate}
            \item
            Suppose that $c_2>c_2^{\mathrm{g.bun}}$. \begin{enumerate}
                \item Let
            $(\mathcal E,\theta)$ be the polystable representative of a closed point ${}^{GF}{\operatorname{Higgs}}^{\mathrm{tf,ss}}_{r,c_1,c_2}$.
            Then the Harder--Narasimhan filtration of $\mathcal E$ has length
            $r$. 
            \item  Writing
            $\mathcal E_i
            :=
            \operatorname{gr}^{\mathrm{HN}}_i(\mathcal E)$,
            $1\leq i\leq r$,
            there are zero-dimensional subschemes
            $D_1\supseteq D_2\supseteq\cdots\supseteq D_r$
            such that
            $\mathcal E_i
            \cong
            \mathcal E_1^{\vee\vee}
            \otimes
            \mathcal L^{-(i-1)}
            \otimes
            \mathcal I_{D_i},
            \
            1\leq i\leq r$ ,and $
            c_1(\mathcal E_1)=\delta
            \ \text{in }\mathrm H^2(X,\mathbb Q)$.
            Moreover,
            \[
            \sum_{i=1}^r \mathrm{length}D_i
            =
            \frac{r(r^2-1)}{24}c_1(\mathcal L)^2
            +
            \frac{\Delta(\mathcal E)}{2r}.
            \]
            \end{enumerate}
            \item
            If $c_2\gg 0$, the moduli space ${\operatorname{Higgs}}_{r,c_1,c_2}^{\text{tf,ss}}$ is reducible.
            \end{enumerate}
        \end{theorem}
            
            \begin{proof}
            We first prove (1).
            Since $c_2>c_2^{\mathrm{g.bun}}$ and the divisibility equation has a solution, the generic fibers of
            the Hitchin map are nonempty by
            Theorem~\ref{thm:criterion of existence}.
            
            Let $C\subset {}^{GF}{\operatorname{Higgs}}^{\mathrm{tf,ss}}_{r,c_1,c_2}$ be an irreducible component containing $({\mathcal{E}},\theta)$.
            By Corollary~\ref{cor:genericsmooth}, after restricting the Hitchin
            base to a nonempty Zariski open subset $W\subset\bm A$, every point
            of $C\cap h^{-1}(W)$ is obtained from a spectral sheaf of the form
            $\pi_s^*\mathcal N\otimes\mathcal I_{\mathfrak D}$, $c_1(\mathcal N)=\delta$.
            Its pushforward is $\mathcal N\otimes\pi_{s*}\mathcal I_{\mathfrak D}$. Since
            \[
            \pi_{s*}\mathcal O_{X_s}
            \cong
            \bigoplus_{j=0}^{r-1}\mathcal L^{-j}
            \]
            and the quotient $\pi_{s*}\mathcal O_{X_s}/\pi_{s*}\mathcal I_{\mathfrak D}$ is zero-dimensional. Put $H:=c_1(\mathcal L)$, $\ell:=H^2$, the Harder--Narasimhan factors of the
            underlying sheaf at such a point have generic rank one and slopes
            \[
            \delta\cdot H,\quad
            \delta\cdot H-\ell,\quad
            \ldots,\quad
            \delta\cdot H-(r-1)\ell.
            \]
            In particular, the successive slope gaps are all equal to $\ell$.

            We claim that the above Harder--Narasimhan polygon is maximal
            among those arising from semistable Higgs sheaves with rank $r$
            and first Chern class $c_1$. Indeed, let
            $\lambda_1\geq\cdots\geq\lambda_r$ be the Harder--Narasimhan
            slopes of the underlying sheaf of any such Higgs sheaf, repeated
            according to their ranks. By Lemma~\ref{lem:inequality from Higgs},
            \[
            \lambda_a-\lambda_b\leq (b-a)\ell,
            \qquad 1\leq a<b\leq r.
            \]
            Hence, for every $1\leq k<r$,
            \[
            \operatorname{HNP}_{\mathcal F}(k)-k\mu_H(\mathcal F)
            =
            \frac{1}{r}
            \sum_{\substack{1\leq a\leq k\\ k<b\leq r}}
            (\lambda_a-\lambda_b)
            \leq
            \frac{\ell}{2}k(r-k).
            \]
            For the spectral sheaves above,
            \[
            \lambda_a=\delta\cdot H-(a-1)\ell,
            \qquad 1\leq a\leq r,
            \]
            and equality holds for every $k$. Thus the spectral
            Harder--Narasimhan polygon is maximal.
            
            On the framed parameter scheme $R^{ss}$, Shatz's upper semicontinuity \cite[Corollary following Proposition~10, p.~183]{Sha77}
            shows that the locus on which the Harder--Narasimhan polygon
            dominates $\operatorname{HNP}_{\mathrm{sp}}$ is closed and
            $\mathrm{SL}(V)$-invariant. It contains the inverse image of
            $C\cap h^{-1}(W)$. Since its image under the good quotient is closed; hence its image
            contains $C$. Being closed and invariant, this locus contains the
            closed orbit corresponding to $({\mathcal{E}},\theta)$. Therefore
            \[
            \operatorname{HNP}_{\mathcal E}
            \geq
            \operatorname{HNP}_{\mathrm{sp}}.
            \]
            The maximality proved above gives the reverse inequality, and hence
            \[
            \operatorname{HNP}_{\mathcal E}
            =
            \operatorname{HNP}_{\mathrm{sp}}.
            \]
            As a result, the Harder--Narasimhan filtration of
            $\mathcal E$ has $r$ graded pieces, each of generic rank one with $\mu_H(\mathcal E_i) -\mu_H(\mathcal E_{i+1}) =
            \ell$, $1\leq i<r$. 
            
            To show the second part of $(1)$, we next identify the double duals of the graded pieces and the
            induced maps between them. For $1\leq i<r$, the slope inequalities imply
            \[
            \theta(\operatorname{HN}_i\mathcal E)
            \subset
            \operatorname{HN}_{i+1}\mathcal E\otimes\mathcal L.
            \]
            Hence the Higgs field induces a morphism
            $\overline\theta_i:
            \mathcal E_i
            \rightarrow
            \mathcal E_{i+1}\otimes\mathcal L$.
            By the semistability of the Higgs sheaf, this morphism is nonzero. 
            Since $\mathrm {gr}^\text {HN}_i\mathcal E$ and $\mathrm {gr}^\text {HN}_{i+1}\mathcal E$ are rank-one
            torsion-free sheaves $\overline\theta_i$ extends uniquely to a
            nonzero morphism of line bundles
            $$\overline\theta_i^{\vee\vee}:
            \mathcal E_i^{\vee\vee}
            \rightarrow
            \mathcal E_{i+1}^{\vee\vee}\otimes\mathcal L.$$  The two line bundles
            have the same $H$-degree, Since $H$ is ample, the zero divisor of this morphism must vanish; hence it is an isomorphism. Equivalently, $\mathcal E_i^{\vee\vee}
            \cong
            \mathcal E_1^{\vee\vee}
            \otimes
            \mathcal L^{-(i-1)}$.

            For each $i$, write
            $\mathcal E_i
            \cong
            \mathcal E_i^{\vee\vee}\otimes\mathcal I_{D_i}$
            for a zero-dimensional subscheme $D_i\subset X$.
            Then the morphism $\overline\theta_i$ becomes
            $\mathcal I_{D_i}
            \hookrightarrow
            \mathcal I_{D_{i+1}}$.
            Therefore
            \[
            D_i\supseteq D_{i+1},
            \qquad
            1\leq i<r.
            \]
            
            Furthermore,
            \[
            c_1(\mathcal E)
            =
            \sum_{i=1}^r c_1(\mathcal E_i)
            =
            r\,c_1(\mathcal E_1)
            -
            \frac{r(r-1)}{2}H.
            \]
            Comparing this with
            $
            r\delta
            =
            c_1+\frac{r(r-1)}{2}H$
            gives
            $c_1(\mathcal E_1)=\delta$ in $\mathrm H^2(X,\mathbb Q)$.

            Finally, Lemma~\ref{lem:discriminant equality}, together with
            $\operatorname{rank}(\mathcal E_i)=1$, gives
            $$\frac{\Delta(\mathcal E)}{2r} =
            \sum_{i=1}^r \frac{\Delta(\mathcal E_i)}{2}
            -
            \frac{1}{2r}
            \sum_{1\leq i<j\leq r}
            \bigl(c_1(\mathcal E_i)-c_1(\mathcal E_j)\bigr)^2 =
            \sum_{i=1}^r \mathrm{length}D_i
            -
            \frac{\ell}{2r}
            \sum_{1\leq i<j\leq r}(j-i)^2.$$
            Using
            \[
            \sum_{1\leq i<j\leq r}(j-i)^2
            =
            \frac{r^2(r^2-1)}{12},
            \]
            we obtain
            \[
            \sum_{i=1}^r\mathrm{length}D_i
            =
            \frac{r(r^2-1)}{24}c_1(\mathcal L)^2
            +
            \frac{\Delta(\mathcal E)}{2r}.
            \]
            This proves (1).
            
            We now prove (2).
            For $c_2\gg0$, the generic fibers of the Hitchin map are
            nonempty, so ${}^{GF}{\operatorname{Higgs}}^{\mathrm{tf,ss}}_{r,c_1,c_2}$
            is a nonempty union of irreducible components. On the other hand, by the asymptotic existence results of
            Donaldson~\cite{Donaldson90}, Friedman~\cite{Fri98}, Gieseker--Li~\cite{GL94}, and
            O'Grady~\cite{O93,O96}, for $c_2\gg0$ there exists a slope-stable torsion-free sheaf $\mathcal F$ on $X$ with the prescribed Chern classes $(r,c_1,c_2)$. Hence $(\mathcal F,0)$
            is a stable Higgs sheaf and determines a point of
            ${\operatorname{Higgs}}^{\mathrm{tf,ss}}_{r,c_1,c_2}$. The underlying sheaf $\mathcal F$ is slope stable, so its
            Harder--Narasimhan type is $(r)$. By part~(1), this point cannot belong to ${}^{GF}{\operatorname{Higgs}}^{\mathrm{tf,ss}}_{r,c_1,c_2}$, the polystable representative of each of whose closed points has Harder--Narasimhan type $(1,\ldots,1)$. Thus ${}^{GF}{\operatorname{Higgs}}^{\mathrm{tf,ss}}_{r,c_1,c_2}$ is a nonempty proper union of irreducible components of ${\operatorname{Higgs}}^{\mathrm{tf,ss}}_{r,c_1,c_2}$. Therefore the latter moduli space is reducible.
            \end{proof}

        \begin{theorem}\label{thm:small c2}
			We assume that the Equation:
			\[
			r\delta=c_1+\frac{r(r-1)}{2}c_1(\mathcal{L})
			\]
			has a solution $\delta\in {\operatorname {CH}}_1(X)$. 
			\begin{enumerate}
				\item If $c_2<c_2^{\text{g.bun}}$, then ${\operatorname{Higgs}}^{\text{tf},ss}_{r,c_1,c_2}$ is empty.
				\item If $c_2=c_2^{\text{g.bun}}$, for Higgs sheaf $({\mathcal{E}},\theta)\in {\operatorname{Higgs}}^{\text{tf},ss}_{r,c_1,c_2}$, $\mathcal E$ is always locally free. And the Harder--Narasimhan type of $\mathcal E$ is that each graded piece $\mathcal E_i$ is locally free of rank $1$.  Moreover, $\mathcal{E}_{i}\cong \mathcal{E}_1^{\vee\vee}\otimes \mathcal{L}^{-(i-1)}$ for $i=1,\ldots, r$ with $c_1({\mathcal{E}}_1)=\delta$. 
			\end{enumerate}
		\end{theorem}
		\begin{proof} Consider the equality (\ref{eq:c2c2gbun}), we have
			$$c_2(\mathcal E)-c_2^{\text{g.bun}}=\frac{\Delta(\mathcal E)}{2r}+\frac{r(r^2-1)}{24}c_1(\mathcal L)^2.$$
			by Lemma \ref{lem:discriminant equality} :
			\[
			\frac{\Delta(\mathcal{E})}{r}=\sum_{i=1}^{m}\frac{\Delta(\mathcal{E}_i)}{r_i}-\sum_{1\le i<j\le m}\frac{r_ir_j}{r}(\frac{c_1(\mathcal{E}_i)}{r_i}-\frac{c_1(\mathcal{E}_j)}{r_j})^2.
			\]
			Thus 
			$$c_2(\mathcal E)-c_2^{\text{g.bun}}=\sum_{i=1}^{m}\frac{\Delta(\mathcal{E}_i)}{2r_i}+\frac{r(r^2-1)}{24}c_1(\mathcal L)^2-\sum_{1\le i<j\le m}\frac{r_ir_j}{2r}(\frac{c_1(\mathcal{E}_i)}{r_i}-\frac{c_1(\mathcal{E}_j)}{r_j})^2.$$
			By Lemma \ref{lem:HN inequality}, we have:
			\begin{equation}\label{eq:calculation 2}
				\sum_{1\le i<j\le m}\frac{r_ir_j}{2r}(\frac{c_1(\mathcal{E}_i)}{r_i}-\frac{c_1(\mathcal{E}_j)}{r_j})^2\le \frac{1}{2r}\sum_{1\le i<j\le m}r_ir_j(j-i)^2c_1({\mathcal{L}})^{2}.
			\end{equation}
			By the following Lemma \ref{lem:olympic} that we have:
			\[
			\sum_{1\le i<j\le m}r_ir_j(j-i)^2\le \frac{r^2(r^2-1)}{12},
			\]
			and the equality holds only if each $r_i=1$. So
			$$c_2(\mathcal E)-c_2^{\text{g.bun}}\geq\sum_{i=1}^{m}\frac{\Delta(\mathcal{E}_i)}{2r_i}$$
			and the equality holds only if each $r_i=1$. 
			
			By the Bogomolov inequality and the slope semistability of each Harder--Narasimhan factor, we have $\frac{\Delta(\mathcal{E}_i)}{2r_i}\geq 0$. From this, if $c_2(\mathcal E)-c_2^\mathrm {g.bun}<0$ then the moduli space is empty and we get (1). 
			
			Moreover, if  $c_2(\mathcal E)=c_2^\mathrm {g.bun}$, then each $r_i=1$ and each $\Delta(\mathcal E_i)=2c_2(\mathcal E_i)=2\mathrm{length}D_i=0$. Thus $D_i=\emptyset$. So $\mathcal E$ is locally free on $X$. And the Harder--Narasimhan type of $\mathcal E$ is that each graded piece $\mathcal E_i$ is locally free of rank $1$. The remaining part can be proved analogously to the first part of Theorem \ref{thm:graded piece are simple}.
		\end{proof}

        \begin{lemma}\label{lem:olympic}
        Let
        $r=r_1+\cdots+r_m$
        be an ordered partition of $r$, where $r_i\in\mathbb Z_{>0}$. Then
        \[
        \sum_{1\leq i<j\leq m} r_i r_j(j-i)^2
        \leq
        \frac{r^2(r^2-1)}{12}.
        \]
        Moreover, equality holds if and only if $m=r$ and $r_1=\cdots=r_r=1$.
        \end{lemma}
        
        \begin{proof}
        For an ordered partition $(r_1,\ldots,r_m)$ of $r$, set
        \[
        S(r_1,\ldots,r_m)
        :=
        \sum_{1\leq i<j\leq m}r_i r_j(j-i)^2.
        \]
        
        Suppose that $r_k>1$, and write $r_k=a+b$  with $a,b\in\mathbb Z_{>0}$. Consider the refinement obtained by
        splitting the $k$-th part into two consecutive parts:
        \[
        (r_1,\ldots,r_{k-1},a,b,r_{k+1},\ldots,r_m).
        \]
        Denote the values of $S$ before and after the splitting by
        $S_{\mathrm{old}}$ and $S_{\mathrm{new}}$, respectively. A direct
        comparison of the terms gives
        \[
        S_{\mathrm{new}}-S_{\mathrm{old}}
        ={}ab
        +b\sum_{p<k}r_p\bigl(2(k-p)+1\bigr)+a\sum_{q>k}r_q\bigl(2(q-k)+1\bigr)+\sum_{p<k<q}r_pr_q\bigl(2(q-p)+1\bigr).
        \]
        All terms on the right-hand side are nonnegative, and the first term
        satisfies $ab>0$. Hence
        \[
        S_{\mathrm{new}}>S_{\mathrm{old}}.
        \]
        Thus splitting any part greater than one strictly increases $S$.
        Repeating this procedure, we eventually obtain the ordered partition
        $(1,\ldots,1)$ with $r$ parts. Consequently,
        \[
        S(r_1,\ldots,r_m)
        \leq S(\underbrace{1,\ldots,1}_{r\text{ times}}),
        \]
        and equality can hold only when the original ordered partition is
        already $(1,\ldots,1)$.
        
        It remains to calculate the maximal value. We have
        \[
        S(\underbrace{1,\ldots,1}_{r\text{ times}})
        =\sum_{1\leq i<j\leq r}(j-i)^2=\sum_{d=1}^{r-1}(r-d)d^2=r\sum_{d=1}^{r-1}d^2-\sum_{d=1}^{r-1}d^3=\frac{r^2(r^2-1)}{12}.
        \]
        This proves both the inequality and the equality statement.
        \end{proof}
        
        		Following Simpson's notation on the $\mathbb G_m$-fixed Higgs fields in \cite{Sim92},
        		\begin{definition}\label{def:systems of hodge pair}
        			We define that a system of Hodge sheaves is a Higgs sheaf $(\mathcal {F},\theta)$ with a decomposition
        			$\mathcal {F}\cong \oplus\mathcal {F}_i$,
        			such that $\theta$ is decomposed as a direct sum of $\theta_i:\mathcal {F}_i\to \mathcal {F}_{i+1}\otimes \mathcal L$. 
        		\end{definition}
        		\begin{theorem}\label{thm:only nested hilbert}
        Assume that the generic fibers of the Hitchin map
        $h^{\mathrm{tf}}$ are nonempty. Let
        $x\in
        \left({\operatorname{Higgs}}^{\mathrm{tf},ss}_{r,c_1,c_2}\right)^{\mathbb G_m}$
        be a closed point fixed by the $\mathbb G_m$-action, and let
        $(\mathcal E,\theta)$ be the polystable representative of the
        $S$-equivalence class corresponding to $x$. Then
        $(\mathcal E,\theta)$ has the structure of a system of Hodge sheaves; namely, there is a decomposition 
        $\mathcal E=\bigoplus_{p\in\mathbb Z}\mathcal F_p$ such that
        $\theta(\mathcal F_p)
        \subset
        \mathcal F_{p+1}\otimes\mathcal L$ for every $p\in\mathbb Z$.
        
        Moreover, suppose that the underlying sheaf $\mathcal E$ has
        Harder--Narasimhan type $(1,\ldots,1)$. Write
        $\mathcal E_i
        :=
        \operatorname{gr}^{\mathrm{HN}}_i(\mathcal E)
        =
        \frac{\mathrm{HN}_i\mathcal E}
             {\mathrm{HN}_{i-1}\mathcal E}$, $1\leq i\leq r$. Assume that the morphisms induced by the Higgs field
        \[
        \overline{\theta}_i:
        \mathcal E_i
        \rightarrow
        \mathcal E_{i+1}\otimes\mathcal L
        \]
        are nonzero for all $1\leq i<r$. Then the
        system-of-Hodge-sheaves decomposition is compatible with the
        Harder--Narasimhan filtration. More precisely, there is an
        isomorphism
        \[
        \mathcal E\cong\bigoplus_{i=1}^{r}\mathcal E_i
        \]
        under which
        \[
        \theta=\bigoplus_{i=1}^{r-1}\theta_i,
        \qquad
        \theta_i:
        \mathcal E_i
        \rightarrow
        \mathcal E_{i+1}\otimes\mathcal L.
        \]
        In particular, the conclusion applies to the polystable
        representative of every closed point in each of the following
        $\mathbb G_m$-fixed loci:
        \begin{enumerate}
            \item
            \[
            \left(
            {}^{GF}{\operatorname{Higgs}}^{\mathrm{tf},ss}_{r,c_1,c_2}
            \right)^{\mathbb G_m},
            \qquad
            c_2>c_2^{\mathrm{g.bun}};
            \]
            \item
            \[
            \left(
            {\operatorname{Higgs}}^{\mathrm{tf},ss}_{r,c_1,c_2}
            \right)^{\mathbb G_m},
            \qquad
            c_2=c_2^{\mathrm{g.bun}}.
            \]
        \end{enumerate}
        
        		\end{theorem}
        	
        \begin{proof}
        Let $ x\in
        \left({\operatorname{Higgs}}^{\mathrm{tf},ss}_{r,c_1,c_2}\right)^{\mathbb G_m}$ be the fixed closed point, and let $(\mathcal E,\theta)$ be the
        polystable representative of the corresponding $S$-equivalence
        class. We first prove that $(\mathcal E,\theta)$ is a system of Hodge
        sheaves.
        
        We follow Simpson's argument in
        \cite[Section~4, p.~44]{Sim92}. Fix a single element
        $t\in k^\times$ of infinite order. Since $x$ is fixed by the
        $\mathbb G_m$-action, the Higgs sheaves $(\mathcal E,\theta)$ and
        $(\mathcal E,t\theta)$ determine the same point of the coarse moduli
        space. Moreover, $(\mathcal E,t\theta)$ is polystable. By the
        uniqueness of the polystable representative of an $S$-equivalence
        class up to isomorphism, there exists an isomorphism
        \[
        \phi_t:
        (\mathcal E,\theta)
        \xrightarrow{\sim}
        (\mathcal E,t\theta).
        \]
        Equivalently,
        \begin{equation}\label{eq:fixed-higgs-isomorphism}
        (\phi_t\otimes\mathrm{id}_{\mathcal L})\circ\theta
        =
        t\theta\circ\phi_t,
        \end{equation}
        as expressed by the commutative diagram
        \[
        \begin{tikzcd}
        \mathcal E
          \ar[r,"\theta"]
          \ar[d,"\phi_t"']
        &
        \mathcal E\otimes\mathcal L
          \ar[d,"\phi_t\otimes\mathrm{id}_{\mathcal L}"]
        \\
        \mathcal E
          \ar[r,"t\theta"']
        &
        \mathcal E\otimes\mathcal L.
        \end{tikzcd}
        \]
        We emphasize that only this single element $t$ and the corresponding
        isomorphism $\phi_t$ will be used. No compatibility between
        isomorphisms associated with different elements of $\mathbb G_m$ is
        required.
        
        The automorphism $\phi_t$ extends uniquely to an automorphism
        $\phi_t^{\vee\vee}$ of the vector bundle $\mathcal E^{\vee\vee}$.
        Since $X$ is projective, the coefficients of the characteristic
        polynomial of $\phi_t^{\vee\vee}$ are constant. Thus
        $\mathcal E^{\vee\vee}$ admits the generalized eigensubbundle
        decomposition
        \[
        \mathcal E^{\vee\vee}
        =
        \bigoplus_{\lambda}
        \left(\mathcal E^{\vee\vee}\right)_\lambda.
        \]
        The corresponding primary projectors are polynomials in
        $\phi_t^{\vee\vee}$. Since $\phi_t^{\vee\vee}$ preserves the
        subsheaf $\mathcal E\subset\mathcal E^{\vee\vee}$, these projectors
        also preserve $\mathcal E$. Consequently,
        \[
        \mathcal E
        =
        \bigoplus_{\lambda}\mathcal E_\lambda,
        \qquad
        \mathcal E_\lambda
        :=
        \mathcal E\cap
        \left(\mathcal E^{\vee\vee}\right)_\lambda.
        \]
        
        It follows from \eqref{eq:fixed-higgs-isomorphism} that
        \begin{equation}\label{eq:higgs-eigensheaf}
        \theta(\mathcal E_\lambda)
        \subset
        \mathcal E_{t\lambda}\otimes\mathcal L.
        \end{equation}
        Indeed, if $(\phi_t-\lambda)^N$ annihilates a local section $u$, then
        \eqref{eq:fixed-higgs-isomorphism} implies that
        $(\phi_t-t\lambda)^N$ annihilates $\theta(u)$.
        
        Only finitely many eigenvalues of $\phi_t$ occur. Since $t$ has
        infinite order, the eigenvalues in each orbit under multiplication
        by $t$ form a finite chain
        \[
        \lambda,\ t\lambda,\ldots,t^N\lambda.
        \]
        After indexing the generalized eigensheaves along these chains and
        taking direct sums of the terms having the same index, we obtain a
        decomposition
        \[
        \mathcal E=\bigoplus_{p\in\mathbb Z}\mathcal F_p
        \ \text{such that} \
        \theta(\mathcal F_p)
        \subset
        \mathcal F_{p+1}\otimes\mathcal L.
        \]
        Hence $(\mathcal E,\theta)$ is a system of Hodge sheaves.
        
        We next prove the compatibility with the Harder--Narasimhan
        filtration. By the uniqueness of the Harder--Narasimhan filtration, the
        automorphism $\phi_t$ preserves every
        $\mathrm{HN}_i\mathcal E$. It therefore induces an automorphism
        \[
        \phi_{t,i}:\mathcal E_i\rightarrow\mathcal E_i.
        \]
        Since $\mathcal E_i$ is a rank-one torsion-free sheaf,
        $\phi_{t,i}$ is multiplication by a scalar
        $\lambda_i\in k^\times$.
        
        By assumption, the Higgs field induces a nonzero morphism 
        $\overline{\theta}_i:
        \mathcal E_i
        \rightarrow
        \mathcal E_{i+1}\otimes\mathcal L$
        for every $1\leq i<r$. Passing
        \eqref{eq:fixed-higgs-isomorphism} to the corresponding
        Harder--Narasimhan factors gives
        \[
        \lambda_{i+1}\overline{\theta}_i
        =
        t\lambda_i\overline{\theta}_i.
        \]
        Since $\overline{\theta}_i\neq0$, we obtain $\lambda_{i+1}=t\lambda_i$.
        It follows that
        $\lambda_i=t^{i-1}\lambda_1$, $1\leq i\leq r$. 
        Because $t$ has infinite order, the scalars
        $\lambda_1,\ldots,\lambda_r$ are pairwise distinct.
        
        We now apply the primary decomposition of $\phi_t$. Let
        $\mathcal F_i\subset\mathcal E$ be the generalized
        $\lambda_i$-eigensheaf. Since $\phi_t$ preserves the
        Harder--Narasimhan filtration and acts on its $i$-th graded factor
        by the distinct scalar $\lambda_i$, the natural morphism
        $\mathcal F_i\rightarrow\mathcal E_i$ is an isomorphism. Moreover,
        $\mathrm{HN}_j\mathcal E
        =
        \bigoplus_{i=1}^{j}\mathcal F_i$
        for every $1\leq j\leq r$. Hence the Harder--Narasimhan filtration
        splits and
        \[
        \mathcal E
        \cong
        \bigoplus_{i=1}^{r}\mathcal E_i.
        \]
        Under this identification,
        \eqref{eq:higgs-eigensheaf} becomes $\theta(\mathcal E_i)
        \subset
        \mathcal E_{i+1}\otimes\mathcal L$,
        and therefore
        $\theta
        =
        \bigoplus_{i=1}^{r-1}\theta_i$, $\theta_i:
        \mathcal E_i
        \rightarrow
        \mathcal E_{i+1}\otimes\mathcal L$. Thus the system-of-Hodge-sheaves decomposition is compatible with
        the Harder--Narasimhan filtration.
        
        Finally, for the polystable representatives of closed points in
        \[
        \left(
        {}^{GF}{\operatorname{Higgs}}^{\mathrm{tf},ss}_{r,c_1,c_2}
        \right)^{\mathbb G_m},
        \qquad c_2>c_2^{\mathrm{g.bun}},
        \]
        and in
        \[
        \left(
        {\operatorname{Higgs}}^{\mathrm{tf},ss}_{r,c_1,c_2}
        \right)^{\mathbb G_m},
        \qquad c_2=c_2^{\mathrm{g.bun}},
        \]
        the required Harder--Narasimhan type and the nonvanishing of the
        successive morphisms $\overline{\theta}_i$ follow from
        Theorem~\ref{thm:graded piece are simple} and
        Theorem~\ref{thm:small c2}, respectively. This proves the final
        assertions.
        \end{proof}

       For integers
        \[
        n_1\geq n_2\geq\cdots\geq n_r\geq0,
        \]
        let $X^{[\underline n]}$ denote the nested Hilbert scheme
        parametrizing flags
        \[
        Z_r\subseteq Z_{r-1}\subseteq\cdots\subseteq Z_1\subset X,
        \qquad
        \operatorname{length}(Z_i)=n_i.
        \]
        Equivalently, $X^{[\underline n]}$ parametrizes flags of ideal
        sheaves
        \[
        \mathcal I_{Z_1}\subseteq
        \mathcal I_{Z_2}\subseteq\cdots\subseteq
        \mathcal I_{Z_r}\subseteq\mathcal O_X.
        \]
        It is a closed subscheme of
        $\prod_{i=1}^r\operatorname{Hilb}^{n_i}(X)$, and hence is
        projective.
        
        Fix $\delta\in \operatorname{Im}\bigl(c_1:\operatorname{Pic}(X)\longrightarrow H^2(X,\mathbb Q)
        \bigr)$ and let $\operatorname{Pic}_X^\delta$ denote the union of
        components of the Picard scheme parametrizing line bundles whose
        rational first Chern class is $\delta$. We define
        \[
        \mathrm X_{\delta}^{[\underline n]}
        :=
        \operatorname{Pic}_X^\delta\times X^{[\underline n]}.
        \]
        
        A geometric point of $\mathrm X_{\delta}^{[\underline n]}$ is
        therefore a tuple
        \[
        (B;Z_1\supseteq Z_2\supseteq\cdots\supseteq Z_r),
        \qquad
        c_1(B)=\delta.
        \]
        Such a tuple determines $r$ line bundles
        \[
        B_i:=B\otimes\mathcal L^{-(i-1)},
        \qquad
        1\leq i\leq r,
        \]
        and hence rank-one torsion-free sheaves
        \[
        E_i:=B_i\otimes\mathcal I_{Z_i}.
        \]
        Thus $\mathrm X_{\delta}^{[\underline n]}$ equivalently
        parametrizes tuples
        \[
        (B_1,\ldots,B_r;Z_1,\ldots,Z_r)
        \]
        satisfying
        \[
        B_i\cong B_1\otimes\mathcal L^{-(i-1)}
        \quad\text{and}\quad
        Z_r\subseteq\cdots\subseteq Z_1.
        \]
        Since both
        $\operatorname{Pic}_X^\delta$ and $X^{[\underline n]}$ are
        projective, so is $\mathrm X_{\delta}^{[\underline n]}$.
         
        \begin{lemma}
        \label{prop:nested-Hilbert-to-Higgs}
        Suppose that
        \[
            c_1
            =
            r\delta-\frac{r(r-1)}{2}c_1({\mathcal{L}})
        \]
        and
        \[
            c_2
            =
            \sum_{i=1}^r n_i
            +
            \sum_{1\leq i<j\leq r}
            \bigl(\delta-(i-1)c_1({\mathcal{L}})\bigr)
            \bigl(\delta-(j-1)c_1({\mathcal{L}})\bigr).
        \]
        Then there is a natural proper morphism
        \[
            \Phi_{\delta}^{[\underline {n}]}\colon
            \mathrm{X}_{\delta}^{[\underline {n}]}
            \rightarrow
            \left(
                \operatorname{Higgs}^{\mathrm{tf,ss}}_{r,c_1,c_2}
            \right)^{\mathbb G_m}.
        \]
        Every Higgs sheaf in the image of $\Phi_{\delta}^{[\underline {n}]}$ has
        nonzero nilpotent Higgs field and is slope stable, and hence Gieseker
        stable.
        \end{lemma}
        
        \begin{proof}
         Choose a closed point $x_0\in X$, and let $\mathcal B$ be a Poincar\'e line bundle on
        $X\times\operatorname{Pic}^{\delta}(X)$, rigidified along
        $\{x_0\}\times\operatorname{Pic}^{\delta}(X)$. We continue to denote
        by $\mathcal B$ its pullback to $X\times\mathrm{X}_{\delta}^{[\underline {n}]}.
        $
        
        Let
        \[
            \mathcal I_1\subseteq
            \mathcal I_2\subseteq\cdots\subseteq
            \mathcal I_r\subseteq
            \mathcal O_{X\times\mathrm{X}_{\delta}^{[\underline {n}]}}
        \]
        be the universal flag of ideal sheaves which is flat over $X\times \mathrm{X}_{\delta}^{[\underline {n}]}$. Define
        \[
            \mathcal U_i
            :=
            \mathcal B_i
            \otimes\mathcal I_i,
            \qquad
            \mathcal U
            :=
            \bigoplus_{i=1}^r\mathcal U_i,
        \]
        where
        \[
            p_X\colon X\times\mathrm{X}_{\delta}^{[\underline {n}]}\rightarrow X
        \]
        is the projection.
        
        For $1\leq i<r$, the universal inclusion
        $\mathcal I_i\hookrightarrow\mathcal I_{i+1}$ induces a morphism
        \[
        \begin{split}
            \Theta_i\colon\mathcal U_i
            &=
            \mathcal B_i\otimes\mathcal I_i\\
            &\rightarrow
            \mathcal B_i
                \otimes\mathcal I_{i+1}=
            \mathcal U_{i+1}\otimes p_X^*L.
        \end{split}
        \]
        Set $\Theta|_{\mathcal U_i}=\Theta_i$ for $i<r$ and
        $\Theta|_{\mathcal U_r}=0$. We thereby obtain a flat family of
        ${\mathcal{L}}$-valued Higgs sheaves
            $\Theta\colon
            \mathcal U\rightarrow
            \mathcal U\otimes p_X^*L$ and all fibers have the prescribed Chern classes.

        At a geometric point
        $$s=(B,Z_1\supseteq\cdots\supseteq Z_r)
            \in\mathrm{X}_{\delta}^{[\underline {n}]},$$
        the corresponding Higgs sheaf is
        $$   {\mathcal{E}}=\bigoplus_{i=1}^r {\mathcal{E}}_i,\quad{\mathcal{E}}_i=B\otimes L^{-(i-1)}\otimes I_{Z_i}=B_i\otimes I_{Z_i},$$
        with Higgs field
        $    \theta_i\colon {\mathcal{E}}_i\rightarrow {\mathcal{E}}_{i+1}\otimes L$
        induced by $I_{Z_i}\subseteq I_{Z_{i+1}}$. 
        
        In the following, we prove that $({\mathcal{E}},\theta)$ is slope stable. Let
        $0\neq {\mathcal{F}}\subsetneq {\mathcal{E}}$   
        be a $\theta$-invariant satuarated subsheaf, and put
        $s=\operatorname{rk}({\mathcal{F}})$.  Let ${\eta}$ be the generic point of $X$ and let $K=K(X)$. Since the
        $Z_i$ are zero-dimensional, $(I_{Z_i})_{\eta}\cong K$. After trivializing $B_{\eta}$ and ${\mathcal{L}}_{\eta}$, there is a basis
        $e_1,\ldots,e_r$ of ${\mathcal{E}}_{\eta}$ for which
        \[
            \theta_{\eta}(e_i)=e_{i+1}
            \quad (1\leq i<r),
            \qquad
            \theta_{\eta}(e_r)=0.
        \]
        Thus $\theta_{\eta}$ is a regular nilpotent endomorphism with one
        Jordan block.
        
        The subspace ${\mathcal{F}}_{\eta}\subset {\mathcal{E}}_{\eta}$ is
        $\theta_{\eta}$-invariant and has dimension $s$. Since the restriction
        of $\theta_{\eta}$ to ${\mathcal{F}}_{\eta}$ is nilpotent on an
        $s$-dimensional vector space, $\theta_{\eta}^s({\mathcal{F}}_{\eta})=0$.
        Hence ${\mathcal{F}}_{\eta}\subseteq \ker(\theta_{\eta}^s)$.
        Then by the regularity of $\theta_{{\eta}}$, ${\mathcal{F}}_{\eta}$ is spanned by $\{e_{r-s+1},\ldots,e_r\}$.
        We write ${\mathcal{T}}_s:={\mathcal{E}}_{r-s+1}\oplus\cdots\oplus {\mathcal{E}}_r$.
        The composite morphism
        \[
            {\mathcal{F}}\rightarrow {\mathcal{E}}\rightarrow {\mathcal{E}}/{\mathcal{T}}_s
        \]
        vanishes at the generic point. Since ${\mathcal{E}}/{\mathcal{T}}_s$ is torsion-free, it
        vanishes identically, and therefore
        $   {\mathcal{F}}\subseteq {\mathcal{T}}_s.$
        Both sheaves have rank $s$, so ${\mathcal{T}}_s/{\mathcal{F}}$ is a torsion
        sheaf. Moreover, it injects into ${\mathcal{E}}/{\mathcal{F}}$, which is
        torsion-free, and hence ${\mathcal{F}}={\mathcal{T}}_s.$ 
        Thus the only proper saturated $\theta$-invariant subsheaves of $E$
        are the tails
        \[
            E_{r-s+1}\oplus\cdots\oplus E_r,
            \qquad
            1\leq s<r.
        \] 
        Then we have
       \[
        \mu_H(\mathcal T_q)
        =
        \mu_H(\mathcal E)
        -
        \frac{r-q}{2}c_1(\mathcal L)\cdot H
        <
        \mu_H(\mathcal E),
        \]
        so $(\mathcal E,\theta)$ is slope stable. In particular, it is
        $H$-Gieseker stable. Consequently, the flat family $(\mathcal E,\Theta)$ determines a
        classifying morphism
        \[
            \Phi_{\delta}^{[\underline {n}]}\colon
            \mathrm{X}_{\delta}^{[\underline {n}]}
            \rightarrow
            \operatorname{Higgs}^{\mathrm{tf,s}}_{r,c_1,c_2}.
        \]
        which is proper since $\mathrm{X}_{\delta}^{[\underline {n}]}$ is projective and $\operatorname{Higgs}^{\mathrm{tf,s}}_{r,c_1,c_2}$ is separated.
        
        For $t\in\mathbb G_m$, define an automorphism $\rho(t)$ of
        $\mathcal E$ by $\rho(t)|_{\mathcal E_i}=t^{i-1}\operatorname{id}_{\mathcal E_i}$.
        Then
        \[
            (\rho(t)\otimes\operatorname{id}_L)
            \circ\Theta\circ\rho(t)^{-1}
            =
            t\Theta.
        \]
        It follows that each Higgs sheaf in the family is fixed by the
        $\mathbb G_m$-action which scales the Higgs field. Hence
        $\Phi_{\delta}^{[\underline {n}]}$ factors through $(\operatorname{Higgs}^{\mathrm{tf,s}}_{r,c_1,c_2})^{\mathbb G_m}$. 
        Moreover, $\Theta\neq 0$ because every inclusion
        $\mathcal I_i\hookrightarrow\mathcal I_{i+1}$ is generically
        nonzero. Thus the image is contained in the nonzero-Higgs part of
        the fixed locus.\end{proof}

        By definition, monopole branches are irreducible components in $({\operatorname{Higgs}}^{\text{tf},ss}_{r,c_1,c_2})^{\mathbb{G}_m}$ with nozero Higgs fields. Since $^{GF}{\operatorname{Higgs}}^{\text{tf,ss}}_{r,c_1,c_2}$ is a union of irreducible components of ${\operatorname{Higgs}}^{\text{tf,ss}}_{r,c_1,c_2}$. By Theorem~\ref{thm:graded piece are simple}(1) and Lemma~\ref{prop:nested-Hilbert-to-Higgs}, we can talk about the monopole branches lying in the closure of generic fibers. 
   
		\begin{corollary}
			Monopole branches in $({}^{GF}{\operatorname{Higgs}}^{\text{tf},ss}_{r,c_1,c_2})^{\mathbb{G}_m}$ are contained in the images of the unions of some nested Hilbert schemes $\mathrm{X}_{\delta}^{[\underline {n}]}$, satisfying that:
			\begin{align*}
				\sum_{i=1}^{r}n_i&=\frac{r(r^2-1)}{24}c_1({\mathcal{L}})^2+\frac{\Delta({\mathcal{E}})}{2r},n_1\ge n_2\ge\ldots\ge n_r.
			\end{align*}
		\end{corollary}
	
		\begin{remark}
			The geometry of the nested Hilbert scheme of points can be rather complicated. For the criterion of the smoothness of nested Hilbert schemes of points, see Cheah\cite{Ch98}. Ryan and Taylor \cite[Corollary 3.17]{RT22} show that $X^{[n_1,n_2,\ldots,n_k]}$ is reducible when $k>22$.
		\end{remark}

		\section{Connectedness of the moduli spaces}

    	In the last section, we know that the moduli space is reducible when $c_2\gg 0$. In this section, we show the connectedness of the moduli space in a special case. From now on, we consider rank two Higgs sheaves over the base field $\mathbb C$, i.e. $r=2$, on a very general hyper-surface $X$ in $\mathbb P^3$ of degree $d\geq 5$. By the classical Noether--Lefschetz theorem, we have ${\operatorname{Pic}}(X)\cong \mathbb{Z}$. Let $H$ be the hyperplane in $\mathbb P^3$ and ${\mathcal{L}}$ is chosen as $\mathcal O_{\mathbb P^3}(H)|_{X}$. Since $r=2$, we cannot apply Theorem \ref{thm:iso for surfaces} directly. But we still have the following lemma.    
		\begin{lemma}\label{special rank 2}Assume that the ground field is $\mathbb C$. Let
		$X\subset\mathbb P^3$ be a very general smooth hypersurface of degree
		$d\geq 5$ with Picard number $\rho(X)=1$, and let $\mathcal L=\mathcal O_X(1)$. Then, for a very
		general
		\[
		s=(a_1,a_2)\in\boldsymbol A
		=
		H^0(X,\mathcal L)\oplus H^0(X,\mathcal L^{\otimes 2}),
		\]
		such that $X_s$ is smooth, the pullback map
		$\pi_s^*:\operatorname {Pic}(X)\rightarrow\operatorname {Pic}(X_s)$
		is an isomorphism.
		\end{lemma}
		
		\begin{proof} The proof proceeds in two steps. First, we apply the Noether--Lefschetz theorem for double covers \cite{Bu83} to establish that
        \[
        \rho(X_s)=1.
        \] Then, using this result and the primitiveness of the hyperplane class proved in \cite{LW81}, we deduce that the pullback map
        \[
        \pi_s^*:\operatorname{Pic}(X)\rightarrow\operatorname{Pic}(X_s)
        \]induces an isomorphism of Picard groups.
        
        Buium \cite{Bu83} (see also Lanteri--Struppa \cite{LS89}) shows that, if $\mathcal L$ is very ample, then for a double cover branched along a very general smooth divisor in $\left|\mathcal L^{\otimes 2}\right|$, the covering surface and the base surface have the same Picard number if and only if
        \[
        H^0(X,K_X\otimes\mathcal L)\neq 0.
        \]
        In our situation, 
        $K_X\otimes\mathcal L
        \simeq
        \mathcal O_X(d-3)$. Since $d\geq 5$, for a very general $s$, Buium's theorem implies that
        $\rho(X_s)=\rho(X)=1$.

        We next prove that
        \[
        \pi_s^*:\operatorname{Pic}(X)\rightarrow \operatorname{Pic}(X_s)
        \]
        is an isomorphism. We first show that both $X$ and $X_s$ have trivial
        $\operatorname{Pic}^0$ and torsion-free Picard groups.

        Since $X\subset \mathbb P^3$ is a smooth hypersurface, the the Lefschetz
        hyperplane theorem implies that
        \[
        \pi_1(X)\cong \pi_1(\mathbb P^3)=0.
        \]
        Consequently,
        $H_1(X,\mathbb Z)=0$. Thus $H^1(X,\mathcal O_X)=0$ and by the exponential sequence we obtain
        \[
        \operatorname{Pic}^0(X)=0.
        \]

        Moreover, we claim that $\operatorname{Pic}(X)$ is torsion-free. 
        Since $H_1(X,\mathbb Z)=0$, by the universal coefficient theorem for cohomology, there is a short exact sequence
        \[
        0\rightarrow
        \operatorname{Ext}^1(H_1(X,\mathbb Z),\mathbb Z)
        \rightarrow
        H^2(X,\mathbb Z)
        \rightarrow
        \operatorname{Hom}(H_2(X,\mathbb Z),\mathbb Z)
        \rightarrow 0 .
        \]
        Since $H_1(X,\mathbb Z)=0$, 
        \[
        H^2(X,\mathbb Z)
        \cong
        \operatorname{Hom}(H_2(X,\mathbb Z),\mathbb Z),
        \]
        which shows that $H^2(X,\mathbb Z)$ is torsion-free. Since ${\operatorname{Pic}}^0(X)=0$, by the exponential sequence, there is an injection
        \[
        \operatorname{Pic}(X)_{\mathrm{tors}}
        \hookrightarrow
        H^2(X,\mathbb Z)_{\mathrm{tors}} .
        \]
        Therefore, $\operatorname{Pic}(X)$ is torsion-free.
        
        In fact, $\pi_s^*{\mathcal{L}}$ is also very ample on $X_s$, its linear system induces an embedding $X_s\hookrightarrow \mathbb P^4$ with defining equations:
        \[
        X_s=
        \left\{[x_0:x_1:x_2:x_3:\lambda] \left|
        \begin{aligned}
        F_d(x_0,\ldots,x_3)&=0,\\
        \lambda^2+a_1(x)\lambda+a_2(x)&=0,
        \end{aligned}\right.
        \right\}
        \subset \mathbb P^4.
        \]
        where $s=(a_1,a_2)\in H^0(X,{\mathcal{L}})\oplus H^0(X,{\mathcal{L}}^2)$. Therefore, $X_s$ is a smooth complete intersection surface of type
        $(d,2)$ in $\mathbb P^4$. The same argument applies to $X_s$ and we have $\operatorname {Pic}^0(X_s)=0$ and $\operatorname {Pic}(X_s)$ is torsion free.

        By \cite[ Theorem 2.1 ]{LW81}, the hyperplane class of a smooth complete
        intersection is primitive in the middle integral cohomology. Hence the
        hyperplane class
        $\pi_s^*H$ is a primitive element of $H^2(X_s,\mathbb Z)$. Since $X$ is chosen to be very general such that ${\operatorname{Pic}}(X)=\mathbb Z H$, hence $\pi_s^*$ is  an isomorphism.\end{proof}
        
		Recall that in this section, we put $c_1({\mathcal{L}})=H$. The equations in Theorem \ref{thm:criterion of existence} become
		\begin{align}\label{eq: rank2, special1}
			2\delta&=c_1+H,\\
			\label{eq: rank2, special2}	4c_2-c_1^2&=-H^2+4\pi_{s*}[\mathfrak D].
		\end{align}

		In particular, if $c_1=H$, the equation (\ref{eq: rank2, special1}) always has a solution with $\delta=H$ and the generic fiber of the Hitchin map
		$ h:\bigsqcup_{c_2} {\mathrm{Higgs}}_{2,[H],c_{2}}^\textrm{tf,ss}\to \bm A$ 
		is always non-empty. And in particular, the Hitchin map  for $c_2=c_2^\mathrm {g.bun}$ (see Definition \ref{def:minimal c2})
		 $h: {\mathrm{Higgs}}_{2,[H],c_{2}}^\textrm{bun,ss}\to \bm A
		$ 
		has nonempty generic fiber        
		and the discriminant $\Delta(\mathcal E)=4c_2-c_1^2=-H^2<0$. 
		\begin{remark}
            By the above choice of Chern classes, semistability coincides with stability. To keep notations consistent, we still denote the coarse moduli by ${\operatorname{Higgs}}_{r,H,c_2}^{\text{tf,ss}}$. 
        \end{remark}
		\begin{proposition}\label{prop: Gmfixempty} The $\mathbb G_m$ fixed points $({\mathrm{Higgs}}_{2,[H],c_{2}}^\textrm{tf,ss})^{\mathbb G_m}$ is empty if and only if the moduli space ${\mathrm{Higgs}}_{2,[H],c_{2}}^\textrm{tf,ss}$ is empty.
		\end{proposition}
		\begin{proof}
			By the properness of the Hitchin map $$h_{{2,[H],c_{2}}}: {\mathrm{Higgs}}_{2,[H],c_{2}}^\textrm{tf,ss}\to \bm A,$$
			for each Higgs sheaf $[(\mathcal E,\theta)]$ in the moduli space ${\mathrm{Higgs}}_{2,[H],c_{2}}^\textrm{tf,ss}$, the $\mathbb G_m$ orbit 
			\[\begin{tikzcd}
				\mathbb G_m\arrow[r] \arrow[d] 
				&  {\mathrm{Higgs}}_{2,[H],c_{2}}^\textrm{tf,ss}\arrow[d]&t\mapsto [(\mathcal E,t\theta)] \\
				\mathbb A^1 \arrow[r] \arrow[ru,dotted]
				& \bm A &t\mapsto \textrm {char.poly.}(t\theta)
			\end{tikzcd}\]
			extend to $\mathbb A^1$ by adding $\displaystyle\lim_{t\to 0}[(\mathcal E,t\theta)]$ and the separatedness implies that the limit point is unique.
		\end{proof}

		Since $r=2$, Theorem~\ref{thm:only nested hilbert} shows
        that every closed $\mathbb G_m$-fixed point is represented
        either by a Higgs sheaf with zero Higgs field or by a system
        of Hodge sheaves with two rank-one summands and a nonzero
        Higgs morphism between them. 

        We first identify the two rank-one summands more precisely.
        Let $(\mathcal E,\theta)$ represent a closed fixed point with
        $\theta\neq0$. By Theorem~\ref{thm:only nested hilbert}, after
        ordering the two weight spaces we may write
        \[
        \mathcal E=\mathcal E_1\oplus\mathcal E_2,
        \qquad
        0\neq\theta_1:
        \mathcal E_1\longrightarrow
        \mathcal E_2\otimes\mathcal L,
        \qquad
        \theta|_{\mathcal E_2}=0.
        \]
        Since ${\operatorname{Pic}}(X)=\mathbb Z[\mathcal L]$, we can verify that:
        \[
        \det{\mathcal{E}}_1\cong {\mathcal{L}},\qquad \det{\mathcal{E}}_2\cong {\mathcal{O}}
        \]
        It follows that there are zero-dimensional subschemes
        $D_1,D_2\subset X$ such that
        \[
        \mathcal E_1\cong\mathcal L\otimes\mathcal I_{D_1},
        \qquad
        \mathcal E_2\cong\mathcal I_{D_2}.
        \]
        The morphism induced by $\theta_1$ on the double duals is a
        nonzero endomorphism of $\mathcal L$, hence a scalar isomorphism.
        After rescaling one summand, $\theta_1$ is therefore induced by an
        inclusion
        \[
        \mathcal I_{D_1}\hookrightarrow\mathcal I_{D_2},
        \]
        equivalently, $D_2\subseteq D_1$. Finally,
        \[
        c_2(\mathcal E)
        =
        \operatorname{length}(D_1)+
        \operatorname{length}(D_2).
        \]
        Hence every nonzero-Higgs fixed closed point belongs to one of the
        subsets defined below.
        
        For integers $n_1\geq n_2\geq0$ with $n_1+n_2=c_2$, let
        \[
        F^{(n_1,n_2)}
        \subset
        \left(
        {\operatorname{Higgs}}_{2,[H],c_2}^{\mathrm{tf,ss}}
        \right)^{\mathbb G_m}
        \]
        denote the subset of closed points represented by Higgs sheaves
        \[
        \mathcal E
        =
        \mathcal L\otimes\mathcal I_{D_1}
        \oplus
        \mathcal I_{D_2},
        \qquad
        D_2\subseteq D_1,
        \qquad
        \operatorname{length}(D_i)=n_i,
        \]
        whose Higgs field is zero on $\mathcal I_{D_2}$ and whose
        restriction
        \[
        \mathcal L\otimes\mathcal I_{D_1}
        \longrightarrow
        \mathcal I_{D_2}\otimes\mathcal L
        \]
        is induced by the inclusion
        $\mathcal I_{D_1}\subseteq\mathcal I_{D_2}$.
        We regard $F^{(n_1,n_2)}$ only as a subset of the fixed locus.
		\begin{proposition} If $c_2$ is sufficiently large, the moduli space
			${\mathrm{Higgs}}_{2,[H],c_{2}}^\textrm{tf,ss}$ is connected.
		\end{proposition}
	
		\begin{proof} 	
        By Proposition~\ref{prop: Gmfixempty}, every scaling orbit admits a fixed limit as $t\to0$. It therefore suffices to prove that all closed fixed points lie in a single connected component of ${\operatorname{Higgs}}_{2,[H],c_2}^{\mathrm{tf,ss}}$.
        
        We first construct a connected subset of the Higgs moduli space which contains every fixed point with nonzero Higgs field. We then connect this subset, through a point of $F^{(c_2,0)}$, to the instanton locus.

			\medskip
			\noindent
			\textbf{Step 1: a connected subset containing all nonzero-Higgs fixed points.}
			Let $Y_0$ be the spectral surface corresponds to $0$ in the Hitchin base 
            and 
			$\pi:Y_0\to X$ be the projection. For later use, we formally write $Y_0={\operatorname{Spec}}_X\bigl(\mathcal{O}_X\oplus \mathcal{L}^{-1}\varepsilon\bigr), \varepsilon^2=0$.
			For a zero-dimensional subscheme
			$Z\subset Y_0$ of length $c_2$, define
			\[
			\mathcal F_Z:=\mathcal{I}_{Z}\otimes \pi^*\mathcal{L},
			\qquad
			E_Z:=\pi_*\mathcal F_Z.
			\]
			The $\mathcal{O}_{Y_0}$-module structure on $\mathcal F_Z$ induces an
			$\mathcal{L}$-twisted nilpotent Higgs field
			\[
			\theta_Z:E_Z\rightarrow E_Z\otimes\mathcal{L},
			\qquad \theta_Z^2=0.
			\]
			Since $\pi$ is finite, pushing forward the defining sequence of $Z$ gives
			\[
			0\rightarrow E_Z
			\rightarrow \mathcal{L}\oplus\mathcal{O}_X
			\rightarrow Q_Z
			\rightarrow 0,
			\qquad \mathrm{Length}(Q_Z)=c_2.
			\]
			Consequently,
			\[
			\mathrm{rank}(E_Z)=2,
			\qquad \det(E_Z)\simeq\mathcal{L},
			\qquad c_2(E_Z)=c_2.
			\]
			
			Similar as Lemma \ref{prop:nested-Hilbert-to-Higgs},  we have a morphism:
            \[
            Y_0^{[c_2]}\rightarrow {\operatorname{Higgs}}_{2,[H],c_2}^{\mathrm{tf,ss}}.
            \] The proof is almost the same as that in the proof of \ref{prop:nested-Hilbert-to-Higgs}. For completeness, we show that $(E_Z,\theta_Z)$ is slope stable. At the generic point ${\eta}$ of $X$,
			the ideal $\mathcal{I}_Z$ is the unit ideal, so $\theta_{Z,{\eta}}$ is a nonzero regular
			nilpotent endomorphism of a two-dimensional vector space.  Its only invariant
			line is its kernel.  Let $0\neq G\subsetneq E_Z$ be a rank-one
			$\theta_Z$-invariant saturated subsheaf.  Then			$G_{\eta}=\ker(\theta_{Z,{\eta}})$.
			
			Under the inclusion $E_Z\subset \mathcal{L}\oplus\mathcal{O}_X$, the generic kernel is the
			second summand.  Hence the composite
			\[
			G\rightarrow E_Z\rightarrow\mathcal{L}
			\]
			vanishes generically, and therefore vanishes identically because $\mathcal{L}$ is
			torsion-free.  Thus $G$ is a rank-one subsheaf of $\mathcal{O}_X$, and
			\[
			\mu_H(G)\leq 0.
			\]
			On the other hand,
			\[
			\mu_H(E_Z)=\frac{H^2}{2}>0.
			\]
			It follows that $(E_Z,\theta_Z)$ is slope-stable, hence Gieseker-stable as a
			Higgs sheaf.  The construction is compatible with families and therefore
			defines a morphism
			\[
			\Psi_{c_2}:Y_0^{[c_2]}\rightarrow {\operatorname{Higgs}}_{2,[H],c_2}^{\mathrm{tf,ss}}.
			\]
			
			Now let $D_2\subset D_1$ with $\mathrm{Length}(D_i)=n_i$ be two closed subscheme of $X$.  Since
			$\mathcal{I}_{D_1}\subset\mathcal{I}_{D_2}$, the sheaf
			\[
			\mathcal J_{D_1,D_2}
			:=\mathcal{I}_{D_1}\oplus
			\bigl(\mathcal{L}^{-1}\mathcal{I}_{D_2}\bigr)\varepsilon
			\subset \mathcal{O}_{Y_0}
			\]
			is an ideal.  Its quotient has length $n_1+n_2$, and
			\[
			\pi_*\bigl(\mathcal J_{D_1,D_2}\otimes\pi^*\mathcal{L}\bigr)
			=\mathcal{L}\mathcal{I}_{D_1}\oplus\mathcal{I}_{D_2}.
			\]
			The induced Higgs field is
			\[
			\theta_0\big|_{\mathcal{I}_{D_2}}=0,
			\qquad
			\theta_0\big|_{\mathcal{L}\mathcal{I}_{D_1}}:
			\mathcal{L}\mathcal{I}_{D_1}\hookrightarrow \mathcal{L}\mathcal{I}_{D_2}
			=\mathcal{I}_{D_2}\otimes\mathcal{L}.
			\]
			Thus
			\[
			\Psi_{c_2}\bigl([\mathcal J_{D_1,D_2}]\bigr)
			=\bigl(\mathcal{L}\mathcal{I}_{D_1}\oplus\mathcal{I}_{D_2},\theta_0\bigr),
			\]
			which is precisely the fixed point represented by the nested pair
			$D_2\subset D_1$. In particular, we have $F^{(n_1,n_2)}\subset \Psi_{c_2}(Y_0^{[c_2]})$. By \cite[\S 2]{RaSa24}, $Y_0^{[c_2]}$ is connected. Hence all $F^{(n_1,n_2)}$ lie in the same connected component of ${\operatorname{Higgs}}_{2,[H],c_2}^{\mathrm{tf,ss}}$

			\medskip
			\noindent
			\textbf{Step 2: connecting a point of $F^{(c_2,0)}$ to the instanton locus.}
			Assume that $c_2$ is sufficiently large.  Choose a general reduced
			zero-dimensional subscheme $D\subset X$ of length $c_2$ satisfying the
			Cayley--Bacharach property with respect to
			$K_X\otimes\mathcal{L}$.        
            Such a subscheme exists for $c_2\gg0$. Indeed, if $c_2>h^0(X,K_X\otimes\mathcal L)$, then a general set of $c_2$ distinct points has the property that no nonzero section of $K_X\otimes\mathcal L$ vanishes on any subset of $c_2-1$ of these points. It therefore satisfies the Cayley--Bacharach property with respect to $K_X\otimes\mathcal L$. The Serre--Cayley--Bacharach construction
			\cite[Theorem~5.1.1]{HL10} gives a locally free extension
			\[
			0\rightarrow\mathcal{O}_X
			\xrightarrow{i} E
			\xrightarrow{q}\mathcal{L}\mathcal{I}_D
			\rightarrow0.
			\]
			We first show that $E$ is $H$-stable.  Suppose that a saturated rank-one
			subsheaf $M\subset E$ satisfies
			$\mu_H(M)\geq\mu_H(E)=H^2/2$.  Since ${\operatorname{Pic}}(X)=\mathbb Z[\mathcal{L}]$, one has
			$M^{\vee\vee}\simeq\mathcal{L}^k$ for some $k\geq1$.  As $E$ is locally free, the
			inclusion $M\hookrightarrow E$ extends across the finite singular locus of
			$M$ to a nonzero map $\mathcal{L}^k\to E$. Therefore
			\[
		    \mathrm H^0\bigl(X,E\otimes\mathcal{L}^{-k}\bigr)\neq0.
			\]
			Twisting the extension by $\mathcal{L}^{-k}$ gives
			\[
			0\rightarrow\mathcal{L}^{-k}
			\rightarrow E\otimes\mathcal{L}^{-k}
			\rightarrow\mathcal{L}^{1-k}\mathcal{I}_D
			\rightarrow0.
			\]
			For $k\geq2$, both outer terms have no nonzero global sections, while for
			$k=1$ one has $\mathrm H^0(X,\mathcal{I}_D)=0$.  This is a contradiction.  Hence $E$ is
			stable.
			
			Let
			\[
			\xi\in\operatorname{Ext}^1_X(\mathcal{L}\mathcal{I}_D,\mathcal{O}_X)
			\]
			be the extension class.  The line $s\xi$, $s\in{\mathbb {A}}^1$, gives a flat family
			\[
			0\rightarrow\mathcal{O}_X
			\xrightarrow{i_s}E_s
			\xrightarrow{q_s}\mathcal{L}\mathcal{I}_D
			\rightarrow0.
			\]
			Define
			\[
			\theta_s:
			E_s\xrightarrow{q_s}\mathcal{L}\mathcal{I}_D
			\hookrightarrow\mathcal{L}
			\xrightarrow{i_s\otimes1_{\mathcal{L}}}E_s\otimes\mathcal{L}.
			\]
			For $s\neq0$, the sheaf $E_s$ is isomorphic to $E$ and is therefore stable.
			For $s=0$, one obtains
			\[
			(E_0,\theta_0)
			=\bigl(\mathcal{O}_X\oplus\mathcal{L}\mathcal{I}_D,\theta_0\bigr),
			\]
			which is the point of $F^{(c_2,0)}$ corresponding to $D$ and is stable by Step~1.
			Thus the family $(E_s,\theta_s)$ is a curve in ${\operatorname{Higgs}}_{2,[H],c_2}^{\mathrm{tf,ss}}$ joining $F^{(c_2,0)}$ to
			$(E,\theta_1)$.
			
			Finally, since $E$ itself is stable, the family
			\[
			t\mapsto(E,t\theta_1),
			\qquad t\in{\mathbb {A}}^1,
			\]
			consists entirely of stable Higgs sheaves and joins $(E,\theta_1)$ to
			$(E,0)\in F_{\mathrm{ins}}$.  For $c_2\gg0$, the instanton branch is irreducible by
			the results of Gieseker--Li and O'Grady \cite{GL94,O96}; hence the
			entire fixed locus lies in a single connected component of ${\operatorname{Higgs}}_{2,[H],c_2}^{\mathrm{tf,ss}}$.
			
			Let now $(V,\theta)\in{\operatorname{Higgs}}_{2,[H],c_2}^{\mathrm{tf,ss}}$ be arbitrary.  By Proposition~\ref{prop: Gmfixempty}, the scaling orbit admits a limit as $t\to0$, the orbit map
			\[
			\mathbb{G}_m\rightarrow{\operatorname{Higgs}}_{2,[H],c_2}^{\mathrm{tf,ss}},
			\qquad t\mapsto(V,t\theta),
			\]
			extends at $t=0$ to a point of the fixed locus.  The image of the extended
			${\mathbb {A}}^1$ is connected, so $(V,\theta)$ lies in the same connected component as
			its fixed-point limit.  Since all fixed points lie in one connected
			component, ${\operatorname{Higgs}}_{2,[H],c_2}^{\mathrm{tf,ss}}$ is connected.
		\end{proof} 
		\bibliographystyle{alpha}

\begin{thebibliography}{BLM{\etalchar{+}}21}
		
		\bibitem[AK79]{AK79}
		A.~B. Altman and S.~L. Kleiman.
		\newblock Compactifying the {Picard} scheme. {II}.
		\newblock {\em Am. J. Math.}, 101:10--41, 1979.
		
		\bibitem[AK80]{AK80}
		A.~B. Altman and S.~L. Kleiman.
		\newblock Compactifying the {Picard} scheme.
		\newblock {\em Adv. Math.}, 35:50--112, 1980.
		
		\bibitem[Ara12]{Ar12}
		D.~Arapura.
		\newblock {\em Algebraic geometry over the complex numbers}.
		\newblock Universitext. Berlin: Springer, 2012.
		
		\bibitem[Art91]{Artamkin91}
		I.~V. Artamkin.
		\newblock Deforming torsion-free sheaves on an algebraic surface.
		\newblock {\em Math. USSR, Izv.}, 36(3):449--485, 1991.
		
		\bibitem[BLM{\etalchar{+}}21]{BLMNPS21}
		A.~Bayer, M.~Lahoz, E.~Macr\`i, H.~Nuer, A~Perry, and P.~Stellari.
		\newblock Stability conditions in families.
		\newblock {\em Publ. Math. Inst. Hautes \'Etudes Sci.}, 133:157--325, 2021.
		
		\bibitem[BNR89]{BNR}
		A.~Beauville, M.S. Narasimhan, and S.~Ramanan.
		\newblock Spectral curves and the generalised theta divisor.
		\newblock {\em J. Reine Angew. Math.}, 398:169--179, 1989.
		
		\bibitem[ByB73]{Bia73}
		A.~Bia\l~ynicki Birula.
		\newblock Some theorems on actions of algebraic groups.
		\newblock {\em Ann. of Math. (2)}, 98:480--497, 1973.

\bibitem[Bui83]{Bu83}
Alexandru Buium.
\newblock Sur le nombre de {Picard} des rev{\^e}tements doubles des surfaces
  alg{\'e}briques.
\newblock {\em C. R. Acad. Sci., Paris, S{\'e}r. I}, 296:361--364, 1983.


        
		\bibitem[Che98]{Ch98}
		J.~Cheah.
		\newblock Cellular decompositions for nested {H}ilbert schemes of points.
		\newblock {\em Pacific J. Math.}, 183(1):39--90, 1998.

  \bibitem[CN20]{CN20}
T.~H. Chen and B.~C. Ng{\^o}.
\newblock On the {Hitchin} morphism for higher-dimensional varieties.
\newblock {\em Duke Math. J.}, 169(10):1971--2004, 2020.


		
		\bibitem[DI87]{DI87}
		P.~Deligne and L.~Illusie.
		\newblock Rel\`evements modulo {$p^2$} et d\'ecomposition du complexe de {de
			Rham}.
		\newblock {\em Invent. Math.}, 89(2):247--270, 1987.
		
		\bibitem[Don90]{Donaldson90}
		S.~K. Donaldson.
		\newblock Polynomial invariants for smooth four-manifolds.
		\newblock {\em Topology}, 29(3):257--315, 1990.

        \bibitem[Del73]{Del73}
        P.Deligne.
        \newblock {\em Le Théorème de Noether.}
        \newblock{Lect. in Math.}, 40:328--340, 1973
        \newblock Springer-Verlag
		\bibitem[EV92]{EVvanishing}
		H.~Esnault and E.~Viehweg.
		\newblock {\em Lectures on vanishing theorems}, volume~20 of {\em DMV Seminar}.
		\newblock Birkh\"{a}user Verlag, Basel, 1992.
		
		\bibitem[FGI{\etalchar{+}}05]{FGA}
		B.~Fantechi, L.~G{\"o}ttsche, L.~Illusie, S.~Kleiman, N.~Nitsure, and
		A.~Vistoli.
		\newblock {\em Fundamental algebraic geometry: {Grothendieck}'s FGA explained},
		volume 123 of {\em Mathematical Surveys and Monographs}.
		\newblock American Mathematical Society, 2005.
		
		\bibitem[Fri98]{Fri98}
		R.~Friedman.
		\newblock {\em Algebraic surfaces and holomorphic vector bundles}.
		\newblock Universitext. New York, NY: Springer, 1998.
		
		\bibitem[Fuj83]{Fuj83}
		T.~Fujita.
		\newblock Vanishing theorems for semipositive line bundles.
		\newblock In {\em Algebraic geometry ({T}okyo/{K}yoto, 1982)}, volume 1016 of
		{\em Lecture Notes in Math.}, pages 519--528. Springer, Berlin, 1983.


\bibitem[GGPN24]{GGPN24}
Guillermo Gallego, Oscar Garc\'{\i}a-Prada, and M.~S. Narasimhan.
\newblock Higgs bundles twisted by a vector bundle.
\newblock {\em International Journal of Mathematics}, 35(09):2441007, 2024.

  
		\bibitem[GL94]{GL94}
		D.~Gieseker and J.~Li.
		\newblock Irreducibility of moduli of rank-2 vector bundles on algebraic
		surfaces.
		\newblock {\em J. Differ. Geom.}, 40(1):23--104, 1994.



  
		\bibitem[Got94]{Go94}
		P.~B. Gothen.
		\newblock The {Betti} numbers of the moduli space of stable rank 3 {Higgs}
		bundles on a {Riemann} surface.
		\newblock {\em Int. J. Math.}, 5(6):861--875, 1994.
		
		\bibitem[GPH13]{Hein13}
		J.~O. Garc\'{\i}a-Prada and J.~Heinloth.
		\newblock The {$y$}-genus of the moduli space of {${\rm PGL}_n$}-{H}iggs
		bundles on a curve (for degree coprime to {$n$}).
		\newblock {\em Duke Math. J.}, 162(14):2731--2749, 2013.
		
		\bibitem[GPHS14]{Hein14}
		J.~O. Garc\'{\i}a-Prada, J.~Heinloth, and A.~Schmitt.
		\newblock On the motives of moduli of chains and {H}iggs bundles.
		\newblock {\em J. Eur. Math. Soc. (JEMS)}, 16(12):2617--2668, 2014.
		
		\bibitem[Gro61]{EGAIII-1}
		A.~Grothendieck.
		\newblock \'{E}l\'{e}ments de g\'{e}om\'{e}trie alg\'{e}brique. {III}.
		\'{E}tude cohomologique des faisceaux coh\'{e}rents. {I}.
		\newblock {\em Inst. Hautes \'{E}tudes Sci. Publ. Math.}, 11:167, 1961.
		
		\bibitem[Gro68]{SGA2}
		A.~Grothendieck.
		\newblock {\em Cohomologie locale des faisceaux coh\'{e}rents et
			th\'{e}or\`emes de {L}efschetz locaux et globaux {$(SGA$} {$2)$}}.
		\newblock North-Holland Publishing Co., Amsterdam; Masson \& Cie, \'{E}diteur,
		Paris, 1968.
		\newblock Augment\'{e} d'un expos\'{e} par Mich\`ele Raynaud, S\'{e}minaire de
		G\'{e}om\'{e}trie Alg\'{e}brique du Bois-Marie, 1962, Advanced Studies in
		Pure Mathematics, Vol. 2.
		
		\bibitem[GSY20a]{GSY20L}
		A.~Gholampour, A.~Sheshmani, and S.T. Yau.
		\newblock Localized {D}onaldson-{T}homas theory of surfaces.
		\newblock {\em Amer. J. Math.}, 142(2):405--442, 2020.
		
		\bibitem[GSY20b]{GSY20}
		A.~Gholampour, A.~Sheshmani, and S.T. Yau.
		\newblock Nested {H}ilbert schemes on surfaces: virtual fundamental class.
		\newblock {\em Adv. Math.}, 365:107046, 50, 2020.
		
		\bibitem[GT20a]{GT20}
		A.~Gholampour and R.~P. Thomas.
		\newblock Degeneracy loci, virtual cycles and nested {H}ilbert schemes, {I}.
		\newblock {\em Tunis. J. Math.}, 2(3):633--665, 2020.
		
		\bibitem[GT20b]{GT20II}
		A.~Gholampour and R.~P. Thomas.
		\newblock Degeneracy loci, virtual cycles and nested {H}ilbert schemes {II}.
		\newblock {\em Compos. Math.}, 156(8):1623--1663, 2020.
		
		\bibitem[GW20a]{GW20}
		U.~G{\"o}rtz and T.~Wedhorn.
		\newblock {\em Algebraic geometry {I}. {Schemes}. {With} examples and
			exercises}.
		\newblock Springer Stud. Math. -- Master. Wiesbaden: Springer Spektrum, 2nd
		edition, 2020.
		
		\bibitem[GW20b]{GoWe20}
		U.~G{\"o}rtz and T.~Wedhorn.
		\newblock {\em Algebraic geometry {I}. {Schemes}. {With} examples and
			exercises}.
		\newblock Springer Stud. Math. -- Master. Wiesbaden: Springer Spektrum, 2nd
		edition edition, 2020.
		
		\bibitem[Har70]{Har70}
		R.~Hartshorne.
		\newblock {\em Ample subvarieties of algebraic varieties}, volume 156 of {\em
			Lecture Notes in Mathematics}.
		\newblock Springer, Berlin, Heidelberg, 1970.

  \bibitem[HL24]{HL24}
Siqi He and Jie Liu.
\newblock On the spectral variety for rank two higgs bundles.
\newblock {\em Proceedings of the London Mathematical Society}, 129(5):e70004,
  2024.

\bibitem[HLM24]{HLM24}
Siqi He, Jie Liu, and Ngaiming Mok.
\newblock The {Spectral} base and quotients of bounded symmetric domains.
\newblock Preprint, {arXiv}:2401.15852 [math.{AG}] (2024), 2024.


		
		\bibitem[Hei10]{HeinLec}
		J.~Heinloth.
		\newblock Lectures on the moduli stack of vector bundles on a curve.
		\newblock In {\em Affine flag manifolds and principal bundles}, Trends Math.,
		pages 123--153. Birkh\"{a}user/Springer Basel AG, Basel, 2010.
		
		\bibitem[Hei15]{Hein15}
		J.~Heinloth.
		\newblock A conjecture of {H}ausel on the moduli space of {H}iggs bundles on a
		curve.
		\newblock {\em Ast\'{e}risque}, (370):157--175, 2015.
		
		\bibitem[Hit87]{Hit87}
		N.J. Hitchin.
		\newblock The self-duality equations on a {Riemann} surface.
		\newblock {\em Proc. London Math. Soc.}, 55(3):39--126, 1987.
		
		\bibitem[HL10]{HL10}
		D.~Huybrechts and M.~Lehn.
		\newblock {\em The geometry of moduli spaces of sheaves}.
		\newblock Cambridge University Press, 2010.
		
		\bibitem[HPL21]{HPL21}
		V.~Hoskins and S.~Pepin~Lehalleur.
		\newblock On the {V}oevodsky motive of the moduli space of {H}iggs bundles on a
		curve.
		\newblock {\em Selecta Math. (N.S.)}, 27(11):1--37, 2021.
		
		\bibitem[HT03]{HT03r}
		T.~Hausel and M.~Thaddeus.
		\newblock Relations in the cohomology ring of the moduli space of rank 2
		{H}iggs bundles.
		\newblock {\em J. Amer. Math. Soc.}, 16(2):303--329, 2003.
	

        \bibitem[Ji24]{Ji24}
        L.~Ji.
        \newblock{\em The {N}oether-{L}efschetz theorem in arbitrary characteristic}
        \newblock{\em J. Algebraic Geom.}, 33(3):567--600, 2024


		\bibitem[JK22]{JK22}
		Y.~Jiang and M.~Kool.
		\newblock Twisted sheaves and {${\rm SU}(r)/\Bbb Z_r$} {V}afa-{W}itten theory.
		\newblock {\em Math. Ann.}, 382(1-2):719--743, 2022.
		
		\bibitem[Jos95]{Jos95}
		K.~Joshi.
		\newblock A {N}oether-{L}efschetz theorem and applications.
		\newblock {\em J. Algebraic Geom.}, 4(1):105--135, 1995.
		
		\bibitem[KM76]{KM76}
		F.~Knudsen and D.~Mumford.
		\newblock The projectivity of the moduli space of stable curves. {I}:
		{Preliminaries} on ''det'' and ''{Div}''.
		\newblock {\em Math. Scand.}, 39:19--55, 1976.
		
		\bibitem[KM98]{KM98}
		J.~Koll\'{a}r and S.~Mori.
		\newblock {\em Birational geometry of algebraic varieties}, volume 134 of {\em
			Cambridge Tracts in Mathematics}.
		\newblock Cambridge University Press, Cambridge, 1998.
		\newblock With the collaboration of C. H. Clemens and A. Corti, Translated from
		the 1998 Japanese original.
		
		\bibitem[Kol86]{Kol86I}
		J.~Koll\'{a}r.
		\newblock Higher direct images of dualizing sheaves. {I}.
		\newblock {\em Ann. of Math. (2)}, 123(1):11--42, 1986.
		
		\bibitem[Kol90]{Kol90}
		J.~Koll{\'a}r.
		\newblock Projectivity of complete moduli.
		\newblock {\em J. Differ. Geom.}, 32(1):235--268, 1990.
		
		\bibitem[Laa20]{Laa20}
		T.~Laarakker.
		\newblock Monopole contributions to refined {V}afa-{W}itten invariants.
		\newblock {\em Geom. Topol.}, 24(6):2781--2828, 2020.
		
		\bibitem[Laa21]{Laa21}
		T.~Laarakker.
		\newblock Vertical {V}afa-{W}itten invariants.
		\newblock {\em Selecta Math. (N.S.)}, 27(4):Paper No. 56, 28, 2021.


\bibitem[LS89]{LS89}
Antonio Lanteri and Daniele~C. Struppa.
\newblock Topological properties of cyclic coverings branched along an ample
  divisor.
\newblock {\em Can. J. Math.}, 41(3):462--479, 1989.

        
		\bibitem[Laz04]{Laz04}
		R.~Lazarsfeld.
		\newblock {\em Positivity in Algebraic Geometry {I,II}}, volume 48,49 of {\em
			Ergebnisse der Mathematik und ihrer Grenzgebiete. 3. Folge}.
		\newblock Springer-Verlag Berlin Heidelberg, 2004.

\bibitem[LW81]{LW81}
A.~S. Libgober and J.~W. Wood.
\newblock On the topological structure of even-dimensional complete
  intersections.
\newblock {\em Trans. Am. Math. Soc.}, 267:637--660, 1981.



        
		\bibitem[Mar77]{Mar77}
		M.~Maruyama.
		\newblock Moduli of stable sheaves. {I}.
		\newblock {\em J. Math. Kyoto Univ.}, 17(1):91--126, 1977.
		
		\bibitem[Mar78]{Mar78}
		M.~Maruyama.
		\newblock Moduli of stable sheaves. {II}.
		\newblock {\em J. Math. Kyoto Univ.}, 18(3):557--614, 1978.
		
		\bibitem[MFK94]{MFK94}
		D.~Mumford, J.~Fogarty, and F.~Kirwan.
		\newblock {\em Geometric invariant theory.}, volume~34 of {\em Ergeb. Math.
			Grenzgeb.}
		\newblock Berlin: Springer-Verlag, 3rd enl. ed. edition, 1994.
		
		\bibitem[Mil80]{Mil80}
		J.~S. Milne.
		\newblock {\em \'{E}tale cohomology}, volume~33 of {\em Princeton Mathematical
			Series}.
		\newblock Princeton University Press, Princeton, N.J., 1980.
		
		\bibitem[Nit91]{N91}
		N.~Nitsure.
		\newblock Moduli space of semistable pairs on a curve.
		\newblock {\em Proc. London Math. Soc.}, 62(3):275--300,
		1991.
		
		\bibitem[O'G93]{O93}
		K.~G. O'Grady.
		\newblock The irreducible components of moduli spaces of vector bundles on
		surfaces.
		\newblock {\em Invent. Math.}, 112(3):585--613, 1993.
		
		\bibitem[O'G96]{O96}
		K.~G. O'Grady.
		\newblock Moduli of vector bundles on projective surfaces: {Some} basic
		results.
		\newblock {\em Invent. Math.}, 123(1):141--207, 1996.
		
		\bibitem[RS06]{RS06}
		G.~V. Ravindra and V~Srinivas.
		\newblock The grothendieck-lefschetz theorem for normal projective varieties.
		\newblock {\em J. Alg. Geom.}, 15(3):563--590, 2006.
		
		\bibitem[RS09]{RS09}
		G.~V. Ravindra and V.~Srinivas.
		\newblock The {N}oether-{L}efschetz theorem for the divisor class group.
		\newblock {\em J. Algebra}, 322(9):3373--3391, 2009.
		
		\bibitem[RT22]{RT22}
		T.~Ryan and G.~Taylor.
		\newblock Irreducibility and singularities of some nested {H}ilbert schemes.
		\newblock {\em J. Algebra}, 609:380--406, 2022.
		
		\bibitem[Sas]{Sasha11}
		Sasha.
		\newblock {Chern} classes in flat families.
		\newblock mathoverflow.net.
		\newblock
		\url{https://mathoverflow.net/questions/52667/chern-classes-in-flat-families}.
		
		\bibitem[Sha77]{Sha77}
		S.~S. Shatz.
		\newblock The decomposition and specialization of algebraic families of vector
		bundles.
		\newblock {\em Compos. Math.}, 35:163--187, 1977.
		
		\bibitem[{Sim}92]{Sim92}
		C.~T. {Simpson}.
		\newblock Higgs bundles and local systems.
		\newblock {\em Publ. Math. Inst. Hautes \'Etudes Sci.}, 75:5--95, 1992.
		
		\bibitem[Sim94]{Sim94II}
		C.~T. Simpson.
		\newblock Moduli of representations of the fundamental group of a smooth
		projective variety. {II}.
		\newblock {\em Publ. Math. Inst. Hautes \'Etudes Sci.}, 80:5--79, 1994.

  \bibitem[SS24]{SS24}
Lei Song and Hao Sun.
\newblock On the image of {Hitchin} morphism for algebraic surfaces: the case
  {$\mathrm{GL}_n$}.
\newblock {\em Int. Math. Res. Not.}, 2024(1):492--514, 2024.


		
		\bibitem[{Sta}19]{stacks-project}
		The {Stacks project authors}.
		\newblock The stacks project.
		\newblock \url{https://stacks.math.columbia.edu}, 2019.
		
		\bibitem[Sum74]{Su74}
		H.~Sumihiro.
		\newblock Equivariant completion.
		\newblock {\em J. Math. Kyoto Univ.}, 14:1--28, 1974.
		
		\bibitem[Tho20]{Tho20}
		R.~P. Thomas.
		\newblock Equivariant {$K$}-theory and refined {V}afa-{W}itten invariants.
		\newblock {\em Comm. Math. Phys.}, 378(2):1451--1500, 2020.
		
		\bibitem[TT18]{TT18}
		Y.~{Tanaka} and R.~P. {Thomas}.
		\newblock {Vafa-Witten invariants for projective surfaces. II: Semistable
			case}.
		\newblock {\em {Pure Appl. Math. Q.}}, 13(3):517--562, 2018.
		
		\bibitem[TT20]{TT20}
		Y.~Tanaka and R.~P. Thomas.
		\newblock Vafa-{W}itten invariants for projective surfaces {I}: stable case.
		\newblock {\em J. Algebraic Geom.}, 29(4):603--668, 2020.
		
		\bibitem[Web17]{We17}
		A.~Weber.
		\newblock Hirzebruch class and {Bia{{\l}}ynicki}-{Birula} decomposition.
		\newblock {\em Transform. Groups}, 22(2):537--557, 2017.
		
		\bibitem[Yok91]{Yo91}
		K.~Yokogawa.
		\newblock Moduli of stable pairs.
		\newblock {\em J. Math. Kyoto Univ.}, 31(1):311--327, 1991.
		
		\bibitem[Yok93]{Yo93C}
		K.~Yokogawa.
		\newblock Compactification of moduli of parabolic sheaves and moduli of
		parabolic {Higgs} sheaves.
		\newblock {\em J. Math. Kyoto Univ.}, 33(2):451--504,
		1993.
		\bibitem[RaSa24]{RaSa24}
         R.~Ramkumar and A.~Sammartano.
        \newblock Rational singularities of nested {Hilbert} schemes.
        \newblock {\em Int. Math. Res. Not.}, 2024(2):1061--1122, 2024.
	\end{thebibliography}

	\newcommand{\etalchar}[1]{$^{#1}$}

	\end{document}